\documentclass[english]{article} 
\usepackage[absolute,overlay]{textpos} 
\usepackage{amsthm}
\usepackage{amsmath}
\usepackage{amssymb}
\usepackage[utf8]{inputenc}
\usepackage[english]{babel}
\usepackage{graphicx}
\usepackage{geometry}
\usepackage{caption}
\usepackage{subcaption}
\usepackage{pdfpages}
\usepackage{comment}
\usepackage{bm}
\usepackage{etoolbox}
\usepackage{mathtools}
\usepackage{tikz}

\usepackage{algpseudocode}
\usepackage{algorithm}
\usepackage{mathrsfs} 
\usepackage{multirow}

\usepackage{pgfplots}
\usepackage{pgfplotstable}
\usepackage{stmaryrd} 

\usepackage{hyperref}

\usepackage{booktabs}

\usepackage[style=numeric, giveninits=true, sorting=none, maxnames=99]{biblatex}
\usepackage[noabbrev,capitalise]{cleveref}
\crefformat{appendix}{#2#1#3}

\usetikzlibrary{fillbetween}
\usepgfplotslibrary{fillbetween}

\usetikzlibrary{shapes.misc}
\usetikzlibrary{math} 

\usetikzlibrary{arrows}

\definecolor{darkspringgreen}{rgb}{0., 0.55, 0.3}
\definecolor{dartmouthgreen}{rgb}{0.05, 0.5, 0.06}
\definecolor{etonblue}{rgb}{0.59, 0.78, 0.64}
\definecolor{airforceblue}{rgb}{0., 0.4, 0.66}
\definecolor{arylideyellow}{rgb}{0.91, 0.84, 0.42}
\definecolor{emerald}{rgb}{0.31, 0.78, 0.47}
\definecolor{uclagold}{rgb}{1.0, 0.7, 0.0}
\definecolor{cadmiumorange}{rgb}{0.93, 0.53, 0.18}

\theoremstyle{thmstyleone}

\theoremstyle{thmstyletwo}

\newtheorem{remark}{Remark}

\theoremstyle{thmstylethree}

\newcommand{\uvec}[2][3]{\boldsymbol{#2\mkern-#1mu}\mkern#1mu}

\newcommand\abs[1]{\left\lvert#1\right\rvert}

\newcommand{\R}{\mathbb{R}}

\newcommand{\bu}{\uvec{u}}

\newcommand{\lSSP}{\ell \text{SSP}}

\newcommand{\RIcolor}[1]{{\leavevmode\color{black} #1}}

\newcommand{\bF}{\uvec{f}}
\newcommand{\diff}[1]{{\mathrm{d}{#1}}}

\newcommand{\ubar}{\overline{\uvec{u}}}

\newcommand{\xip}{x_{i+\frac{1}{2}}}
\newcommand{\xin}{x_{i-\frac{1}{2}}}

\newcommand{\iip}{i+\frac{1}{2}}
\newcommand{\iin}{i-\frac{1}{2}}

\newcommand{\iph}{{i+\frac{1}{2}}}
\newcommand{\imh}{{i-\frac{1}{2}}}

\newcommand{\qbar}{\overline{q}}

\begin{document}
\author{Lorenzo Micalizzi\footnote{Affiliation: Department of Mathematics, North Carolina State University, SAS Hall, 2108, 2311 Stinson Dr, Raleigh, NC 27607, United States. Email: lmicali@ncsu.edu}, and Eleuterio F. Toro\footnote{Affiliation: Laboratory of Applied Mathematics, DICAM, University of Trento, Via Mesiano 77, 38123 Trento, Italy. Email: eleuterio.toro@unitn.it}}
\title{Algorithms of very high space--time orders of accuracy for hyperbolic equations in the semidiscrete WENO--DeC framework}

\maketitle

\abstract{
	In this work, we provide a deep investigation of a family of arbitrary high order numerical methods for hyperbolic partial differential equations (PDEs), with particular emphasis on very high order versions, i.e., with order higher than 5. 
	More in detail, within the context of a generic Finite Volume (FV) semidiscretization, we consider Weighted Essentially Non--Oscillatory (WENO) spatial reconstruction and Deferred Correction (DeC) time discretization.
	The goal of this paper is twofold.
	On the one hand, we want to demonstrate the possibility of utilizing very high order schemes in concrete situations and highlight the related advantages.
	On the other one, we want to debunk the myth according to which, 
	in the context of numerical resolution of hyperbolic PDEs with very high order spatial discretizations, 
	the adoption of lower order time discretizations, e.g., strong stability preserving (SSP) or linearly strong stability preserving ($\lSSP$) Runge--Kutta (RK) schemes, does not affect the overall accuracy of the resulting approach and consequently its computational efficiency.
	Numerical results are reported for the linear advection equation (LAE) and for the Euler equations of fluid dynamics, showing the advantages and the critical aspects of the adoption of very high order numerical methods. 
	Overall, the results indicate the potential for their use in real--life applications, offering advantages in terms of efficiency, such as requiring shorter computational times to achieve a prescribed error, even in problems involving discontinuities.
	Furthermore, the results confirm order degradation and efficiency loss when coupling very high order space discretizations with lower order SSPRK time discretizations.
}

\section{Introduction}
In recent years, the whole scientific community has been putting a lot of effort in the design and in the investigation of high order methods, in particular but not exclusively, in the context of arbitrary high order frameworks~\cite{Decremi,toro2001towards,dumbser2006arbitrary,toro2024ader,dumbser2008unified,ciallella2022arbitrary,abgrall2023extensions,ciallella2023arbitrary,lore_phd_thesis,micalizzi2023efficient,veiga2024improving,velasco2023spectral,jund2007arbitrary,ADERNSE,titarev2002ader,dumbser2009very,boscheri2019high,micalizzi2023new,micalizzi2024novel,ciallella2024high}, i.e., theoretically allowing for the construction of schemes of ``any'' order with no accuracy barrier.
The main motivation lies in the ability of high order methods to achieve, in regions in which the solution is smooth, smaller errors within coarser discretizations with respect to lower order ones.
Moreover, their reduced error allows to better capture the details of the solution in long--time simulations.
On the other hand, they do not lack critical aspects: their implementation is usually more involved with respect to the one of low order methods and, when discontinuities are present, they are more prone to spurious oscillations, which can lead to simulation blow--ups.
Therefore, generally speaking, high order methods are considered less safe to use in the context of tough simulations, where robustness is strictly required, and practitioners tend to use low order discretizations in real applications.
In particular, the adoption of second order methods~\cite{chakravarthy1983high,vanleer1977towardsiii,roe1986characteristic,sweby1984high,van1982comparative,van1974towards} in problems involving discontinuities is nowadays well--established, with higher order methods being primarily relegated, for now, to the academic context.
Even though convergence analyses on smooth problems attaining the formal order of accuracy have been reported for very high orders~\cite{evstigneev2016construction,Evstigneev2016OnTC,gerolymos2009very,balsara2000monotonicity,shi2003resolution,hermes2012linear,gao2020seventh,dumbser2014posteriori,velasco2023spectral,veiga2024improving}, tests on problems involving discontinuities are rarely conducted with spatial discretizations with order higher than~5.
Moreover, in references considering very high order space discretizations, strong stability preserving (SSP)~\cite{shu1988total,shu1988efficient} or linearly strong stability preserving ($\lSSP$)~\cite{gottlieb2001strong} Runge--Kutta (RK) methods are usually adopted for the time integration.
The accuracy barrier for such schemes, if one wants to have non-negative RK coefficients~\cite{shu1988efficient}, is order 4~\cite{ruuth2002two} on nonlinear problems. 
With the main goal of preventing loss in accuracy due to the mismatch between temporal and spatial order, in some works, a reduction of the time step is performed. However, this causes an increase in the diffusion and makes the resulting approach too computationally costly for real applications.

Here follows a collection of the references that we have found in literature addressing spatial discretizations of order higher than~5.
In~\cite{evstigneev2016construction,Evstigneev2016OnTC}, the author studies Weighted Essentially Non--Oscillatory (WENO) reconstructions~\cite{liu1994weighted,shu1998essentially,shu1989efficient} up to order 13. In particular, in~\cite{evstigneev2016construction} a theoretical analysis is presented for the Burgers equation, while, in \cite{Evstigneev2016OnTC} a numerical investigation on the linear advection equation (LAE), the Burgers equation and the Euler equations is reported. The time discretizations considered are 3-rd and 4-th order accurate, and the time step is chosen, in smooth tests, to be $1/350$ of the cell length in order to avoid loss of accuracy. Although this may work for coarse meshes, and despite the rates corresponding to the spatial accuracy being attained, this makes the order of accuracy of the resulting scheme formally limited by the order of the time discretization.
In~\cite{gerolymos2009very}, WENO reconstructions up to order 17 are investigated on the LAE, on the Burgers equation and on the Euler equations. $\lSSP$ RK methods are adopted for simulations involving the LAE, including a convergence analysis attaining the expected rates, however, SSPRK schemes of order 3 are used for all tests involving nonlinear equations.
In~\cite{balsara2000monotonicity}, the authors consider WENO reconstructions up to order 13 and test them in the context of the numerical resolution of the LAE, of the Burgers equation, of the Euler equations, and of the magnetohydrodynamics equations. Also there, $\lSSP$ RK methods are adopted for simulations involving the LAE, and a 3-rd order SSPRK scheme is adopted for simulations involving nonlinear problems.
In the latter case, the time step is, in most of the tests, suitably chosen to be proportional to a power of the cell length in order not to spoil the accuracy. As already underlined, this practice is not suitable for real applications and results unfeasible for multidimensional problems. This is the reason why, in the same reference, two--dimensional convergence analyses are performed with a fixed Courant--Friedrichs--Lewy (CFL) number.
In~\cite{shi2003resolution}, the authors investigate the capability of WENO reconstructions up to order 9 to capture complex flow structures in the context of the Euler equations. The time discretization is again performed through a 3-rd order accurate SSPRK scheme.
In~\cite{hermes2012linear}, a stability analysis of WENO reconstructions up to order 11 is performed for the LAE. The time discretizations considered are of order 3 and 4, without any adjustment of the time step to make up for the lower order accuracy in time. The reported convergence analyses confirm the fact that the order of the resulting schemes is affected in a negative way by the choice of lower order time discretizations.
In~\cite{gao2020seventh}, WENO reconstructions of order 5, 7 and 9 are considered in the context of the Euler equations and the time discretization is performed through a 3-rd order accurate SSPRK scheme, adjusting the time step size in the convergence analysis in order not to spoil the accuracy.
The only exceptions we found, adopting very high order in space and time, are given by~\cite{toro2001towards,dumbser2006arbitrary,dumbser2014posteriori,popov2024space,velasco2023spectral,veiga2024improving,micalizzi2024impact}.
The ADER--FV framework was introduced in~\cite{toro2001towards}, where results from implementations of nonlinear versions were reported up to order 10 for the LAE.
The same setting is considered in~\cite{dumbser2006arbitrary}, where a theoretical analysis is conducted on the LAE with a simple central linear spatial reconstruction not suitable for discontinuous solutions, and convergence analyses were reported up to order 16.
%
%
In~\cite{dumbser2014posteriori,popov2024space,velasco2023spectral,veiga2024improving}, a discontinuous finite element setting is considered:
in~\cite{dumbser2014posteriori,popov2024space} a fully--discrete ADER--Discontinuous Galerkin framework is assumed, while, in \cite{velasco2023spectral,veiga2024improving} the authors adopt a semidiscrete ADER--Spectral Difference framework. 
In~\cite{dumbser2014posteriori}, simulations up to order 10 in space and time are reported both on smooth problems, with related convergence analysis attaining the expected order, and on tests involving discontinuities in the context of the Euler equations. 
In~\cite{popov2024space}, convergence analyses achieving the expected order are shown up to order 12 and very tough simulations are performed with methods accurate up to order 10 for the reactive Euler equations.  
In~\cite{velasco2023spectral},  convergence analyses are reported up to order 7 and tests on shocks are performed up to order 10 in the context of the Euler equations.
In~\cite{veiga2024improving}, convergence analyses and tests on discontinuities are reported up to order 9 for the Euler equations.
In~\cite{dumbser2014posteriori,popov2024space,velasco2023spectral,veiga2024improving}, the oscillatory character of the finite element setting on non--smooth solutions is tackled via ``a posteriori limiting'', resorting to lower order schemes when the numerical solution violates positivity of density and pressure and a local discrete maximum principle.
In~\cite{micalizzi2024impact}, instead, in the context of the Euler equations, the same setting under investigation here is adopted, namely, a family of ``truly'' (i.e., both in space and time) arbitrary high order methods for hyperbolic partial differential equations (PDEs), based on a finite volume (FV) semidiscretization with WENO space reconstruction and Deferred Correction (DeC)~\cite{micalizzi2023new,ciallella2022arbitrary,Decremi,minion2003semi,Decoriginal} time integration.
The performance of the approach is systematically assessed up to order 7.

The semidiscrete WENO--DeC framework had been already investigated, up to order 5, in~\cite{ciallella2022arbitrary,ciallella2023arbitrary,ciallella2024high}.
In particular, in such references, structure--preserving modifications of the approach for the Shallow Water equations have been proposed. 
In~\cite{micalizzi2024impact}, instead, an investigation has been conducted to assess the role of the numerical flux in high order semidiscrete FV schemes for the Euler equations. 
More in detail, it was shown how the importance of the adopted numerical flux tends to diminish as the order of accuracy increases (without vanishing for the range of investigate orders), and how the order of accuracy of the scheme plays a crucial role in improving the quality of the results, both on smooth and discontinuous solutions.
This work is a follow-up of~\cite{micalizzi2024impact}; here, we push the order of accuracy up to 13, and we report results for the LAE and for the Euler equations of fluid dynamics in one space dimension.
Numerical results show the advantages in the adoption of a truly arbitrary high order framework. Furthermore, reported comparisons with results obtained coupling WENO reconstructions with classical SSPRK schemes confirm the advantages of high order time integration against accuracy and efficiency losses of $3$-rd and $4$-th order SSPRK time discretizations.
We also confirm what conjectured in~\cite{micalizzi2024impact}: in the context of the numerical solution of the Euler equations, the differences between results obtained with different numerical fluxes become smaller as the order of accuracy increases. Nonetheless, depending on the test and on the mesh refinement, such differences may not disappear and still be evident also for order 13.
We leave for future works a systematic investigation of this aspect.
Simulations on tough tests, characterized by strong discontinuities, confirm that the proposed approach is a good baseline framework for the construction of robust schemes. However, the presence of spurious oscillations and some simulation crashes show that further limiting is required to make it fail--safe.
Due to the investigative character of this work, we do not couple the basic framework with any extra strategy, such as a posteriori limiting or fixes in the reconstruction ensuring positivity of density and pressure, which would indeed improve the quality of the results and avoid simulation crashes, but we rather focus on the basic approach obtained by simply coupling WENO and DeC in the FV framework.
Further studies on suitable limiting of the framework are planned for future works.

The paper is structured as follows.
In Section \ref{sec:basic_notions}, we recall the basic notions on semidiscrete FV methods for hyperbolic PDEs.
We describe, in Section~\ref{sec:WENO}, the WENO reconstruction in space.
The DeC time discretization is introduced in Section~\ref{sec:DeC}.
Numerical results are reported in Section~\ref{sec:numerical_results}.
Conclusions and future perspectives follow in Section~\ref{sec:conclusions}.

\section{Semidiscrete Finite Volume schemes}\label{sec:basic_notions}
In this section, we describe the basic structure and background of semidiscrete FV methods for hyperbolic systems of conservation laws in one space dimension.
For multidimensional generalizations of the presented notions, the reader is referred to~\cite{ToroBook,toro2024computational,balsara2000monotonicity,ciallella2022arbitrary,ciallella2024high,lore_phd_thesis,micalizzi2024impact,abgrall1994essentially,wolf2007high,zhu2008runge} and references therein.
One--dimensional hyperbolic systems of conservation laws appear in the form
\begin{equation}\label{eq:sys}
	\frac{\partial}{\partial t} \uvec{u}(x,t) + \frac{\partial}{\partial x}\uvec{f}(\uvec{u}(x,t)) = \uvec{0}, \quad (x,t)\in\Omega\times\mathbb{R}^+_0,
\end{equation}
with $\uvec{u}:\Omega\times\mathbb{R}^+_0\longrightarrow \mathbb{R}^{N_c}$ being the unknown vector of conserved variables, where $N_c\in \mathbb{N}^+$ is the number components of the system, $\uvec{f}:\mathbb{R}^{N_c}\longrightarrow\mathbb{R}^{N_c}$ being flux, and $\Omega:=[x_L,x_R]$ being the space domain.
System~\eqref{eq:sys} is hyperbolic, namely, the Jacobian of the flux, $J(\uvec{u}):=\frac{\partial\uvec{f}}{\partial \uvec{u}}(\uvec{u})$, is real--diagonalizable: it has $N_c$ real eigenvalues and a corresponding set of linearly independent eigenvectors.
In particular, its eigenvalues represent the wave speeds of the system.

In this work, we focus on semidiscrete FV schemes for solving system~\eqref{eq:sys}.
The FV method, originally introduced by Godunov~\cite{Godunov}, is based on a spatial integral formulation of the governing equations.
Due to its favorable numerical and analytical properties, it has gained significant popularity over the years, leading to the development of numerous FV schemes. For a comprehensive overview, we refer readers to classical works on the subject~\cite{hirsch2007numerical,godlewski2021numerical,toro2024computational,shu1988efficient,shu1989efficient,shu1998essentially,leveque2002finite,ToroBook,abgrall1994essentially}.
The basic idea behind semidiscrete schemes is to decouple the discretizations in space and time via the so--called ``method of lines''.

Let us introduce a tessellation of $\Omega$, a family of non--overlapping control volumes $C_{i}:=[x_{\imh},x_{\iph}]$ covering $\Omega$ exactly, which we assume for simplicity to be uniform, namely, we assume $\xip-\xin=\Delta x$ for any $i$.
Making use of the divergence theorem, the spatial integral of the governing equations over the generic cell $C_{i}:=[x_{\imh},x_{\iph}]$ yields
\begin{equation}\label{eq:FV_semidiscretization_exact}
	\frac{d}{dt}\ubar_{i}(t) + \frac{1}{\Delta x}(\uvec{f}_{\iip}-\uvec{f}_{\iin}) = \uvec{0},
\end{equation}
in which $\ubar_{i}$ is the cell average of the solution in $C_{i}$, while, $\uvec{f}_{\iip}$ and $\uvec{f}_{\iin}$ represent the flux at cell boundaries $x_{\iip}$ and $x_{\iin}$ at time $t$
\begin{align}
	\uvec{f}_{\iip} &:= \bF(\bu(\xip,t)). 
	\label{eq:flux F} 
\end{align}

So far, no approximation has been performed and Equation~\eqref{eq:FV_semidiscretization_exact} is, as a matter of fact, exact.
The discretization in space is completed by specifying a suitable approximation for the flux $\uvec{f}_{\iip}$.
This is achieved via a space reconstruction and a numerical flux function.
More in detail, one considers the following discretization 
\begin{equation}\label{eq:FV_semidiscretization}
	\frac{d}{dt}\ubar_{i}(t) + \frac{1}{\Delta x}(\widehat\bF_{\iip}-\widehat\bF_{\iin}) = \uvec{0},
\end{equation}
where, for any $i$, the term  $\widehat\bF_{\iip}$ is defined as
\begin{equation}\label{eq:flux}
	\widehat\bF_{\iip} = \widehat\bF(\bu_{\iip}^L(t),\bu_{\iip}^R(t)),
\end{equation}
where $\widehat\bF(\cdot,\cdot)$ is the numerical flux function, while,  $\bu_{\iip}^L(t)$ and $\bu_{\iip}^R(t)$ are the reconstructions of the numerical solution at the interface $x_{\iip}$ from left and right side respectively at generic time $t$.
The choice of the spatial reconstruction plays a crucial role in the development of a robust (high order) scheme.
In particular, the reconstruction must be able to handle discontinuities or steep gradients reducing spurious oscillations, which could cause simulation blow--ups when numerically solving systems of equations in which positivity is strictly required, such as the Euler equations or the Shallow Water equations.
To this end, the nonlinear character of the reconstruction results essential, as underlined in Godunov's theorem~\cite{Godunov}, see also~\cite{ToroBook,toro2024computational}.
For this reason, here, we consider a WENO space reconstruction locally performed in each cell, which is detailed in Section~\ref{sec:WENO}.
The order of accuracy of the resulting spatial discretization equals the one of the adopted spatial reconstruction.
The choice of the numerical flux does not influence the order of accuracy, however, it is well--known to influence stability and diffusion of the scheme. Plenty of numerical fluxes are available in literature, e.g., Lax--Friedrichs~\cite{lax1954weak}, First--Order Centred~\cite{Toro1996,toro2000centred,chen2003centred}, Rusanov~\cite{Rusanov1961}, Harten--Lax--van Leer~\cite{harten1983upstream}, Central--Upwind~\cite{kurganov2001semidiscrete,kurganov2000new}, Low--Dissipation Central--Upwind~\cite{kurganov2023new}, HLLC~\cite{toro1992restoration,toro1994restoration}, the Godunov flux from the exact Riemann solver~\cite{Godunov}.
The interested reader is referred to~\cite{toro2016riemann,qiu2006numerical,qiu2007numerical,qiu2008development,san2015evaluation,hongxia2020numerical,leidi2024performance,micalizzi2024impact} for more information and comparisons.

The system of ordinary differential equations (ODEs) obtained considering the semidiscretization~\eqref{eq:FV_semidiscretization} for all cell averages must be numerically solved in time through a suitable scheme.
The order of accuracy of the resulting method is $P$ if both the space and time discretizations are at least $P$-th order accurate.
In the numerical solution of hyperbolic PDEs with very high order spatial reconstructions, it is a very common practice to adopt lower order time discretizations, see~\cite{evstigneev2016construction,Evstigneev2016OnTC,gerolymos2009very,balsara2000monotonicity,shi2003resolution,hermes2012linear,gao2020seventh}, where SSP or $\lSSP$ RK time integration schemes are considered.
On general nonlinear problems, the order of accuracy of such schemes is well--known to be limited to 4~\cite{ruuth2002two}, although $\lSSP$ RK schemes could reach arbitrary high order of accuracy on linear problems.
In principle, one could modify RK methods through downwind computation of numerical fluxes in some stages to achieve, at the same time, arbitrary high order and the SSP property~\cite{shu1988efficient}. However, this modification is never considered in the aforementioned works dealing with very high order spatial reconstructions, thus, the order of accuracy of the schemes is limited from the one of the adopted time discretization.
In some of those references~\cite{balsara2000monotonicity,gao2020seventh}, the time step is adjusted to be proportional to a power of the cell length to make up for the choice of a lower order time discretization and match the orders of accuracy in space and time.
The resulting schemes are however ``computationally impracticable''~\cite{balsara2000monotonicity} especially on multidimensional tests, due to the huge computational time required, which is also associated to higher levels of diffusion and to accumulation of rounding errors.
Let us also notice that, as stated in~\cite{balsara2000monotonicity}, when coupling a high order space discretization with a lower order time integration, a strong time step reduction is actually required to keep the non--oscillatory property.
In order to avoid these issues, here we consider an arbitrary high order time integration strategy, namely, a DeC method~\cite{micalizzi2023new,ciallella2022arbitrary,Decremi,micalizzi2024impact,lore_phd_thesis}, detailed in Section~\ref{sec:DeC}.

\section{Space reconstruction}\label{sec:WENO}
The WENO reconstruction has been introduced in 1994~\cite{liu1994weighted,jiang1996efficient} as an improvement of the Essentially Non-Oscillatory (ENO)~\cite{harten1987uniformly1,harten1987uniformly2} reconstruction.
Since then, a huge number of WENO--based methods has been proposed due to its favorable features.
The interested reader is referred to~\cite{shu1998essentially,shu1988efficient,shu1989efficient,abgrall1994essentially} and references therein for a non-exhaustive overview of the methodology and related applications.
In this section, we will describe the WENO reconstruction for a scalar component $q(x)$ in one space dimension. 
For the extension to higher dimensions, the interested reader is referred to~\cite{lore_phd_thesis,ciallella2022arbitrary,micalizzi2024novel} in a Cartesian framework and to~\cite{abgrall1994essentially,wolf2007high,zhu2008runge} in an unstructured framework.

Let us considered, for $q(x)$, the same discretization setting described in the previous section:
we assume to know its cell average $\qbar_i$ in each cell $C_i:=[x_{\imh},x_{\iph}]$ of a uniform tessellation of the one--dimensional domain $\Omega$.
We will now describe how to perform the local reconstruction of order $2r-1$ in a generic point $x^*$ of the generic cell $C_i$.
The main idea is to reconstruct, in such a point, a high order approximation $q_h^{HO}(x^*)\approx q(x^*)$, obtained from the cell averages of a ``big'' stencil containing $C_i$, and $r$ low order approximations $q_h^{\ell,LO}(x^*)\approx q(x^*)$ $\ell=0,\dots,r-1$, obtained from the cell averages of $r$ ``small'' stencils containing $C_i$.
At this point, assuming that we have some coefficients $d_\ell^{x^*}$, called ``linear weights'', which allow to express the high order approximation $q_h^{HO}(x^*)$ as a linear combination of the low order approximations $q_h^{\ell,LO}(x^*)$, i.e.,
\begin{equation}
	q_h^{HO}(x^*)=\sum_{\ell=0}^{r-1}d_\ell^{x^*} q_h^{\ell,LO}(x^*),
	\label{eq:linear_weights}
\end{equation}
the WENO reconstruction in $x^*$ is obtained by replacing the linear weights $d_\ell^{x^*}$ in \eqref{eq:linear_weights} with modified ones $\omega_\ell^{x^*}$, called ``nonlinear weights'', leading to 
\begin{equation}
	q_h^{WENO}(x^*)=\sum_{\ell=0}^{r-1}\omega_\ell^{x^*} q_h^{\ell,LO}(x^*).
	\label{eq:WENO_approximation}
\end{equation}
In particular, the nonlinear weights $\omega_\ell^{x^*}$ are designed in such a way to recover the linear ones, $\omega_\ell^{x^*}\approx d_\ell^{x^*}$, in smooth cases and to privilege the approximations associated to the smoothest stencils involved in the linear combination~\eqref{eq:WENO_approximation} in case of discontinuous solutions.

In order to complete the description of the method, we need to specify how to construct the approximations $q_h^{HO}(x^*)$ and $q_h^{\ell,LO}(x^*)$ for $\ell=0,\dots,r-1$, and how to define the linear and nonlinear weights $d_\ell^{x^*}$ and $\omega_\ell^{x^*}$ for $\ell=0,\dots,r-1$.
The two problems are strongly connected, as will be explained in the following.

The required approximations, $q_h^{HO}(x^*)$ and $q_h^{\ell,LO}(x^*)$ for $\ell=0,\dots,r-1$, are obtained through the same strategy applied to different stencils.
More in detail, the approximation $q_h^{HO}(x^*)$ of order $2r-1$ is obtained from a ``big'' stencil of $2r-1$ cells
\begin{equation}
	\label{eq:HO_stencil}
	\mathcal{S}^{HO}:=\left\lbrace C_{i-(r-1)},\dots,C_{i+(r-1)} \right\rbrace,
\end{equation}
while, each of the approximations $q_h^{\ell,LO}(x^*)$ of order $r$  is obtained from one of the ``small'' stencils of $r$ cells
\begin{equation}
	\label{eq:LO_sentcil}
	\mathcal{S}^{LO}_\ell:=\left\lbrace C_{i-(r-1)+\ell},\dots,C_{i+\ell} \right\rbrace, \quad \ell=0,\dots,r-1.
\end{equation}

In both cases, starting from the cell averages of a stencil, $\mathcal{S}:=\left\lbrace C_{s},\dots,C_{f} \right\rbrace$ with $N_{\mathcal{S}}:=f-s+1$ cells, we need to reconstruct an approximation $q_h(x^*)$.
To this end, we consider the interpolation of a primitive $Q$ of $q$. 
Let us notice that the point values of $Q$ are available at $N_{\mathcal{S}}+1$ cell interfaces
\begin{align}
	Q(x_{s-\frac{1}{2}})&=0, \label{eq:arbitrary_initial_value}\\
	Q(x_{s-\frac{1}{2}+k})&=\sum_{m=0}^{k-1} \qbar_{s+m}\RIcolor{\Delta x}, \quad k=1,\dots,N_{\mathcal{S}}.
\end{align}
This allows to construct the interpolation polynomial $Q_h(x)$ of degree $N_{\mathcal{S}}$ associated to such point values and defined as
\begin{align}
	Q_h(x):=\sum_{k=0}^{N_{\mathcal{S}}}Q(x_{s-\frac{1}{2}+k})\varphi_k(x),
\end{align}
with $\varphi_k$ being the Lagrange polynomials involved in the interpolation.
The desired approximation of $q(x^*)$ is readily obtained by differentiating the previous expression
\begin{align}
	q_h(x^*):=\sum_{k=0}^{N_{\mathcal{S}}}Q(x_{s-\frac{1}{2}+k})\frac{d}{d x}\varphi_k(x^*).
\end{align}
Let us notice that, since $x^*\in[x_{\imh},x_{\iph}]$ is fixed, the values $\frac{d}{d x}\varphi_k(x^*)$ are fixed as well, therefore,  $q_h(x^*)$ is nothing but a linear combination of the cell averages of the stencil $\mathcal{S}$. 
This holds for the high order stencil~\eqref{eq:HO_stencil} and for all the low order stencils~\eqref{eq:LO_sentcil}, thus, leading to $q_h^{HO}(x^*)$ with order of accuracy $2r-1$ and $q_h^{\ell,LO}(x^*)$ $\ell=0,\dots,r-1$ with order of accuracy $r$.
The linear weights $d_\ell^{x^*}$ are seeked such that the linear combination of the low order approximations $q_h^{\ell,LO}(x^*)$ through such coefficients gives the high order one according to Equation~\ref{eq:linear_weights}, constituting an over--determined linear system to be solved.
Concerning the nonlinear weights, the most popular definition is the one from~\cite{shu1998essentially}, namely
\begin{equation}
\omega_\ell^{x^*} := \frac{\alpha_\ell^{x^*}}{\sum^{r-1}_{k=0}\alpha_k^{x^*}},\quad \alpha_\ell^{x^*} := \frac{d_\ell^{x^*}}{(\beta_\ell+\epsilon_{{\small \text{WENO}}})^2},
\end{equation}
with $\epsilon_{{\small \text{WENO}}}$ being a small constant adopted to avoid divisions by zero, and $\beta_m$ being the smoothness indicators associated to the stencils and defined as
\begin{equation}
	\beta_\ell := \sum_{k=1}^{r-1} \int_{x_{i-1/2}}^{x_{i+1/2}} \left(\frac{d^k}{dx^k} q_h^{\ell,LO}(x)\right)^2 \Delta x^{2k-1}\diff{x},\quad \ell = 0,\ldots,r-1.
\end{equation}
Here, like in~\cite{shu1998essentially,jiang1996efficient,ciallella2022arbitrary,ciallella2023arbitrary,ciallella2024high,micalizzi2024impact}, we consider $\epsilon_{{\small \text{WENO}}}:=10^{-6}$.

\begin{remark}[Reconstruction of characteristic variables]\label{rmk:characteristic_variables}
In the case of systems, the described scalar reconstruction must be applied to each component.
However, as pointed out in~\cite{qiu2002construction,miyoshi2020short,peng2019adaptive,ghosh2012compact}, performing the reconstruction on conserved variables leads to spurious oscillations for high order.
In order to prevent such an issue, the reconstruction should be applied to characteristic variables, i.e., to the scalar components of $L\ubar_{i}$, where $L$ is the matrix of the left eigenvectors of the flux Jacobian $J$.
The reconstructed variables must be then multiplied by $R:=L^{-1}$, where $R$ is the matrix of the right eigenvectors of the flux Jacobian $J$.
From a practical point of view, the reconstruction in the generic cell $C_{i}$ is performed with ``frozen'' matrices $L:=L(\ubar_{i})$ and $R:=R(\ubar_{i})$.
\end{remark}

\RIcolor{
\begin{remark}[On some critical aspects of the computation of the linear weights]
The computation of the linear weights is a critical step of WENO schemes.
For example, as underlined in~\cite{shi2002technique}, a special treatment is needed in presence of negative linear weights in order to achieve stability.
On the other hand, see same reference, negative linear weights do not occur in one--dimensional FV WENO methods for conservation laws for any order of accuracy, and in finite difference WENO methods for conservation laws in any spatial dimension for any order of accuracy.
Concerning existence, unfortunately linear weights do not always exist for arbitrary $x^*$ and order as mentioned in many references, see~\cite{wang2025extremum,zhang2024application,zhang2016eno,zhu2017new,huang2014re}.
Such a pathology is not restricted to the context of non--uniform or unstructured grids only.
A notable example, see~\cite{huang2014re}, is given by WENO3, for which linear weights do not exist at cell center, i.e., for $x^*:=x_{i}$.
Let us remark that, in the one--dimensional framework under investigation, we have found no issues under this point of view for $x^*:=x_{i\pm\frac{1}{2}}$ up to order 31, even though simulations are reported up to order 13 only.
\end{remark}
}

\section{Time discretization}\label{sec:DeC}
In this section, we describe the time discretization adopted in the context of this work.
It consists in a DeC scheme referred as ``bDeC'' in~\cite{micalizzi2023new,lore_phd_thesis}, also used in~\cite{han2021dec,ciallella2022arbitrary,ciallella2023arbitrary,micalizzi2024impact}, whose construction is based on the DeC formalism introduced by Abgrall in~\cite{Decremi} in 2017.

The DeC methodology has a long history, and it has been firstly introduced in~\cite{fox1949some} as an iterative strategy for initial value problems. However, the approach became popular much later in a modern formulation proposed in~\cite{Decoriginal} for the arbitrary high order numerical resolution of ODEs.
This new formulation led to many developments and applications, see inter alia~\cite{minion2003semi,minion2004semi,layton2005implications,huang2006accelerating,minion2011hybrid,ketcheson2014comparison,speck2015multi,boscarino2016error,boscarino2018implicit}.
In particular, a further abstract formulation was proposed by Abgrall in~\cite{Decremi}, leading to numerous other works based on it~\cite{abgrall2019high,abgrall2020high,michel2021spectral,michel2022spectral,ciallella2022arbitrary,abgrall2021relaxation,ciallella2023arbitrary,ciallella2024high,abgrall2024staggered,abgrall2020multidimensional,offner2020arbitrary,han2021dec,micalizzi2023new,micalizzi2023efficient,veiga2024improving,micalizzi2024novel,micalizzi2024impact}.

Let us consider the Cauchy problem
\begin{equation}
	\label{eq:ODE}
	\begin{cases}
		\frac{d}{dt}\uvec{u}(t) = \uvec{G}(t,\uvec{u}(t)),\quad t\in[0,T_f], \\
		\uvec{u}(0)=\uvec{z},
	\end{cases}
\end{equation}
where $\uvec{u}:[0,T_f] \rightarrow \mathbb{R}^{N_c}$ represents the unknown solution, with $N_c\in \mathbb{N}^+$ being the number of equations of ODEs system, $T_f\in \mathbb{R}^+$ is the final time, $\uvec{z} \in \R^{N_c}$ is the known initial condition, and $\uvec{G}: [0,T_f] \times \R^{N_c} \to \R^{N_c}$ represents the right-hand side function, that we assume to satisfy the classical regularity hypotheses guaranteeing well--posedness of problem~\eqref{eq:ODE}, namely, $\uvec{G}$ is Lipschitz-continuous with respect to $\uvec{u}$ uniformly with respect to $t$.

The DeC method considered here is a one--step method; therefore, we assume the classical setting of a generic time interval $[t_n,t_{n+1}]$, where $t_{n+1}-t_n=:\Delta t$, and our goal is to obtain $\uvec{u}_{n+1}\approx \uvec{u}(t_{n+1})$ given $\uvec{u}_{n}\approx \uvec{u}(t_{n})$. 

To this end, we introduce, over $[t_n,t_{n+1}]$, $M+1$ subtimenodes $t^m$ $m=0,\dots,M$ with  
\begin{equation}
	t_n=:t^0<t^1<\dots<t^M:=t_{n+1}.
\end{equation}
The number and the distribution of subtimenodes influences the order of the method~\cite{micalizzi2023new,lore_phd_thesis}. For example, with $M+1$ equispaced subtimenodes one can obtain order of accuracy $M+1$, while, the same number of Gauss--Lobatto subtimenodes leads to accuracy $2M$.
Here, like in~\cite{micalizzi2024novel,ciallella2022arbitrary,ciallella2023arbitrary,ciallella2024high}, we assume Gauss--Lobatto subtimenodes to get higher accuracy with the same computational cost.
Thus, we denote by $\uvec{u}^m$ an approximation of the exact solution $\uvec{u}(t^m)$ of~\eqref{eq:ODE} in the subtimenode $t^m$.
The vectors $\uvec{u}^m$ are to be determined for $m=1,\dots,M$, while, in the first subtimenode the assumption of a one--step method setting provides a sufficiently accurate approximation and we fix $\uvec{u}^0:=\uvec{u}_n$.
We discretize the integral version of the analytical differential problem,
\begin{equation}
	\label{exint}
	\uvec{u}(t^m)-\uvec{u}(t_n)-\int_{t^0}^{t^m}\uvec{G}(t,\uvec{u}(t))dt=\uvec{0}, \quad m=1,\dots,M,
\end{equation}
in a high order fashion into 
\begin{equation}
	\label{eq:apprint}
	\uvec{u}^m-\uvec{u}_n-\Delta t \sum_{\ell=0}^{M} \theta^m_\ell \uvec{G}(t^\ell,\uvec{u}^\ell)=\uvec{0}, \quad m=1,\dots,M,
\end{equation}
where the coefficients $\theta^m_\ell$ are the normalized quadrature weights of the quadrature formula associated to the introduced subtimenodes over $[t^0,t^m]$.
The obtained discretization is a $2M$-th order accurate implicit nonlinear system in the unknown coefficients $\uvec{u}^m$ for $m=1,\dots,M$.
We do not directly solve such a nonlinear system, rather, we consider the following explicit iterative procedure
\begin{align}\label{eq:DeCODE_Remi}
	\uvec{u}^{m,(p)} = \uvec{u}_n+\Delta t \sum_{\ell=0}^{M} \theta^m_\ell \uvec{G}(t^\ell,\uvec{u}^{\ell,(p-1)}), \quad m=1,\dots,M,\quad p>0,
\end{align}
over the approximated values $\uvec{u}^{m,(p)}\approx \uvec{u}(t^m)$.
In this context, we set $\uvec{u}^{m,(p)}:=\uvec{u}_n$ for $m=0$ or $p=0$, and we consider $P$ iterations, with $P$ being the order of accuracy of the scheme. In the end, we set $\uvec{u}_{n+1}:= \uvec{u}^{M,(P)}$.
More in detail, the number of adopted subtimenodes is $M+1$ with $M=\left \lceil \frac{P}{2}\right \rceil$.
It can be proven~\cite{micalizzi2023new,lore_phd_thesis} that the previously defined iterative procedure gains one order of accuracy at each iteration towards the exact solution of the nonlinear system~\eqref{eq:apprint}.
Connections of the approach with RK schemes are put in evidence in~\cite{micalizzi2023new,lore_phd_thesis} and the reader is referred to such references for further information.

\section{Numerical results}\label{sec:numerical_results}
In this section, we test the numerical approach up to order 13 on several benchmarks with different nature.
In particular, in Section~\ref{sec:LAE_1d}, we provide results for the LAE, while, in Section~\ref{sec:Euler_1d}, we provide results for the Euler equations.
We will now state some general guidelines followed in the setup of the simulations.

As in \cite{micalizzi2024impact}, the initialization of each test is performed through numerical computation of the cell averages at the beginning of the simulation using the Gauss--Legendre quadrature with the minimal number of points to attain the order of accuracy of the employed spatial reconstruction.

Concerning the numerical flux, for the LAE, we have adopted a simple and natural upwind numerical flux~\cite{hirsch2007numerical,ToroBook,leveque2002finite,godlewski2021numerical,toro2024computational,AbgrallMishranotes}, corresponding to the exact Riemann solver for such a problem.
Instead, for the Euler equations, we have systematically investigated two possible options, namely:
\begin{itemize}
	\item Rusanov numerical flux~\cite{Rusanov1961}
	\begin{equation}
		\widehat{\bF}^{\text{Rus}}({\bu}^L,{\bu}^R) := \frac{1}{2}\left(\bF({\bu}^R) + \bF({\bu}^L)\right) - \frac{1}{2}s\left({\bu}^R - {\bu}^L\right),
	\end{equation}
	where $s$ is the maximum in absolute value of the local wave speeds associated to the states ${\bu}^L$ and ${\bu}^R$;
	
	\item exact Riemann solver numerical flux~\cite{Godunov}
	\begin{equation}
		\widehat{\bF}({\bu}^L,{\bu}^R) = \bF({\bu}^*(0)),
	\end{equation}
	where ${\bu}^*(\xi)$ is the exact solution of the Riemann problem between the states ${\bu}^L$ and ${\bu}^R$ for $\xi:=\frac{x}{t}$.
\end{itemize}
In particular, the Riemann problem is solved through the strategy proposed in~\cite{toro1989fast}, also detailed in \cite[Chapter 4]{ToroBook}.
The Rusanov numerical flux is an incomplete upwind one, while, the exact Riemann solver numerical flux is a complete upwind one. For more information, the interested reader is referred to~\cite{ToroBook,micalizzi2024impact}.
In order to save space, in the labels, Rusanov numerical flux and exact Riemann solver numerical flux will be compactly denoted respectively as~``Rus'' and ``Ex.RS''.

For what concerns the reconstruction, as already stated, we have adopted a standard WENO reconstruction. However, in the context of the Euler equations, we have considered two options:
\begin{itemize}
	\item reconstruction of conserved variables, i.e., the reconstruction is applied directly to the cell averages $\ubar_{i}$;
	\item reconstruction of characteristic variables, i.e., the reconstruction is applied to the quantities of the vector $L\ubar_{i}$ as detailed in Remark~\ref{rmk:characteristic_variables}.
\end{itemize}
In order to save space, in the labels, reconstruction of conserved and characteristic variables will be compactly denoted respectively as~``cons.'' and ``char.''.

Concerning the computation of the time step, we have considered here the following expression
\begin{equation}
	\Delta t:= C_{CFL} \frac{\Delta x}{\max{(s)}},
	\label{eq:CFL}
\end{equation}
where $C_{CFL}$ is a constant, and $s$ represents an estimate of the maximum wave speed in absolute value.
For the LAE, the wave speed is a constant and a datum of the problem, while, for nonlinear systems wave speeds must be estimated in some way.
In particular, in the context of the Euler equations, we have considered a direct estimate with $s:=\abs{u}+\sqrt{\gamma \frac{p}{\rho}}$, with $\rho,u,p$ being obtained from the cell averages at time $t_n$.
Let us, however, notice that, despite this being a very simple and broadly adopted estimate, it could fail in bounding the actual wave speeds leading to simulation crashes due to violation of the CFL condition, e.g., in Riemann problems with zero initial velocity. See~\cite{toro2020bounds} for more information.

We tried to run each test starting with $C_{CFL}=0.95$, resorting to lower values in case of simulation crashes. One can prove that~\cite{ToroBook}, for the first order versions of the FV schemes considered herein,
the linear CFL stability constraint is $C_{CFL}^{max}=1$ provided that rigorous speed estimates, bounding the actual wave speeds, are adopted for the time step computation.

In all WENO--DeC simulations, the order of the time discretization matches the one of the space discretization.
In order to demonstrate that coupling a high order space discretization with a lower order time discretization determines efficiency losses, we will also provide comparisons between WENO--DeC and the schemes obtained coupling the WENO space reconstruction with SSPRK(3,3) and SSPRK(5,4) time discretizations from~\cite{spiteri2002new}, compactly denoted respectively as ``SSPRK3'' and ``SSPRK4''.
Moreover, along with WENO--SSPRK3 and WENO--SSPRK4, we will consider also modified versions in which the time step is reduced, as in~\cite{balsara2000monotonicity,gao2020seventh}, 
to be proportional to a power of the cell length to match the orders of accuracy in space and time, 
for spatial discretization with order higher than the one of the time discretization. 
We denote such modified versions by ``WENO--mSSPRK3'' and ``WENO--mSSPRK4'' respectively, where the additional letter ``m'' compactly indicates the described time step modification.
In practice, for mSSPRK3 and mSSPRK4, we consider
\begin{equation}
	\Delta t:= C_{CFL} \left(\frac{\Delta x}{\max{(s)}}\right)^{\frac{P}{R}},
	\label{eq:modified_CFL}
\end{equation}
where $R$ is the order of the SSPRK scheme, i.e., $R:=3$ for SSPRK3 and $R:=4$ for SSPRK4, where $P>R$ is the order of the space discretization.
Notice that we do not investigate mSSPRK3 and mSSPRK4 with WENO3, as this would lead to $\frac{P}{R}>1$, causing instabilities in the mesh refinement.
Hence, WENO--mSSPRK3 and WENO--mSSPRK4 are investigated for spatial discretizations with order from 5 on. 
The comparisons carried out herein are based on the concept of ``efficiency'' meant as the error with respect to the computational time~\cite{micalizzi2024impact}.


\begin{remark}[On the accuracy of the coefficients of SSPRK(5,4)]
	The coefficients of SSPRK(5,4) reported in the original work~\cite{spiteri2002new} by Spiteri and Ruuth are accurate up to $10^{-10}$. 
	Therefore, implementations using such coefficients suffer from error stagnation far above machine precision, which is rather unsatisfactory in the context of investigations of high order discretizations. 
	In order to avoid such an issue, it was necessary to recompute the coefficients with higher accuracy.
	They were kindly provided to us by Prof. David Ketcheson from King Abdullah University of Science and Technology, and we report them in Appendix~\ref{app:SSPRK4}, as we are not aware of any reference reporting the coefficients of SSPRK(5,4) with accuracy higher than the one in~\cite{spiteri2002new}.	
\end{remark}

\subsection{Linear advection equation}\label{sec:LAE_1d}
We start with the simplest hyperbolic PDE, the scalar LAE, which reads
\begin{equation}
	\frac{\partial}{\partial t}u+\frac{\partial}{\partial x}\left(au\right)=0,
	\label{eq:LAE_1d}
\end{equation}
with $a$ being a constant that we assume here to be equal to $1$. 
Its unique analytical solution, for given initial condition $u(x,0)=u_0(x)$ over a periodic spatial domain (or over the whole real axis), can be explicitly computed through the method of characteristics, and it reads $u(x,t)=u_0(x-at)$ in the whole space--time domain.
For this model equation, we consider two problems:
a smooth test in Section~\ref{sec:LAE_1d_sin4} in order to verify the order of accuracy, 
and the advection of a composite wave in Section~\ref{sec:LAE_1d_composite_wave} in order to assess the ability of the high order discretizations under investigation to preserve the solution profile, even in very long--time simulations.

The results confirm the advantages of adopting very high order methods over lower order ones, as well as order degradation and efficiency loss when the order of accuracy of the time discretization does not match the one of the space discretization, even if the time step is modified to make up for such mismatch.

\subsubsection{Test 1: Advection of smooth profile}\label{sec:LAE_1d_sin4}
We consider, on the computational domain $\Omega:=[-1,1]$, with periodic boundary conditions, the smooth initial condition 
\begin{equation}
	u_0(x):=\sin^4{\left(\pi x\right)},
\end{equation}
and we run simulations until final time $T_f:=1$ with $C_{CFL}:=0.95.$
The test has been taken from~\cite{shu1998essentially,balsara2000monotonicity}. It was originally introduced in~\cite{rogerson1990numerical,shu1990numerical} as an example for which ENO schemes suffer from order degradation.

Convergence tables for WENO--DeC, including errors in $L^1$-, $L^2$-, and $L^\infty$-norms, along with computational times, are reported in Table~\ref{tab:LAE_1d_convergence_tables_WENO_DeC}, where $N$ represents the number of elements used for the simulations. As one can clearly see, the expected rate of convergence has been obtained in all cases from order 5 on. The trends of the errors in the three norms are very similar. This is not trivial: many approaches struggle or fail to achieve the expected rate in the $L^\infty$-norm.
The only unexpected behavior has been obtained for WENO3, which shows a superconvergent trend, achieving asymptotical experimental orders approximately equal to 4, 4.5 and 5 respectively in the $L^1$-, $L^2$- and $L^\infty$-norms.
Despite this, the average convergence rate is still around 3.


\begin{table}[htbp]
	\centering
	\caption{LAE, Test 1: convergence tables for WENO--DeC}
	\label{tab:LAE_1d_convergence_tables_WENO_DeC}
	\scalebox{0.65}{ 
		\begin{tabular}{c c c c c c c c}
			\toprule
			\multirow{2}{*}{$N$} & \multicolumn{2}{c}{$L^1$ error} & \multicolumn{2}{c}{$L^2$ error} & \multicolumn{2}{c}{$L^{\infty}$ error} & \multirow{2}{*}{CPU Time} \\
			\cmidrule(lr){2-3} \cmidrule(lr){4-5} \cmidrule(lr){6-7}
			& Error & Order & Error & Order & Error & Order & \\
			\midrule
			
			\multicolumn{8}{c}{\textbf{WENO3--DeC3}} \\ 
			\midrule
			160  &   1.425e-02  &  $-$  &   1.844e-02  &  $-$  &   4.155e-02  &  $-$  &   2.838e-02 \\ 
			320  &   2.127e-03  &  2.744  &   3.609e-03  &  2.353  &   1.135e-02  &  1.872  &   8.120e-02 \\ 
			640  &   2.534e-04  &  3.070  &   5.321e-04  &  2.762  &   2.282e-03  &  2.314  &   2.740e-01 \\ 
			1280  &   1.898e-05  &  3.739  &   3.932e-05  &  3.758  &   2.097e-04  &  3.444  &   1.046e+00 \\ 
			2560  &   1.187e-06  &  3.999  &   1.869e-06  &  4.395  &   8.136e-06  &  4.688  &   4.097e+00 \\ 
			5120  &   7.186e-08  &  4.046  &   8.737e-08  &  4.419  &   2.763e-07  &  4.880  &   1.633e+01 \\ 
			Average order    &      &  3.520  &      &  3.537  &      &  3.440  &    \\ 
			\midrule

			\multicolumn{8}{c}{\textbf{WENO5--DeC5}} \\ 
			\midrule
			80  &   9.335e-04  &  $-$  &   8.607e-04  &  $-$  &   1.507e-03  &  $-$  &   2.171e-02 \\ 
			160  &   2.950e-05  &  4.984  &   2.954e-05  &  4.865  &   6.754e-05  &  4.480  &   6.499e-02 \\ 
			320  &   7.685e-07  &  5.263  &   7.526e-07  &  5.295  &   1.760e-06  &  5.262  &   2.156e-01 \\ 
			640  &   1.809e-08  &  5.409  &   1.569e-08  &  5.584  &   2.873e-08  &  5.937  &   7.819e-01 \\ 
			1280  &   4.128e-10  &  5.453  &   3.386e-10  &  5.534  &   4.569e-10  &  5.975  &   3.020e+00 \\ 
			2560  &   9.740e-12  &  5.406  &   8.526e-12  &  5.312  &   1.275e-11  &  5.163  &   1.211e+01 \\ 
			5120  &   2.342e-13  &  5.378  &   2.165e-13  &  5.299  &   3.351e-13  &  5.250  &   4.762e+01 \\ 
			Average order    &      &  5.315  &      &  5.315  &      &  5.344  &    \\ 
			\midrule

			\multicolumn{8}{c}{\textbf{WENO7--DeC7}} \\ 
			\midrule
			80  &   1.027e-04  &  $-$  &   1.189e-04  &  $-$  &   2.194e-04  &  $-$  &   6.193e-02 \\ 
			160  &   3.150e-07  &  8.348  &   3.966e-07  &  8.228  &   9.739e-07  &  7.815  &   2.563e-01 \\ 
			320  &   2.211e-09  &  7.155  &   4.326e-09  &  6.518  &   2.346e-08  &  5.376  &   9.595e-01 \\ 
			640  &   1.441e-11  &  7.261  &   3.715e-11  &  6.864  &   2.627e-10  &  6.480  &   3.750e+00 \\ 
			1280  &   7.825e-14  &  7.525  &   1.266e-13  &  8.196  &   7.583e-13  &  8.436  &   1.476e+01 \\ 
			Average order    &      &  7.572  &      &  7.452  &      &  7.027  &    \\ 
			\midrule

			\multicolumn{8}{c}{\textbf{WENO9--DeC9}} \\ 
			\midrule
			40  &   5.342e-04  &  $-$  &   6.957e-04  &  $-$  &   1.640e-03  &  $-$  &   4.913e-02 \\ 
			80  &   1.071e-06  &  8.962  &   1.374e-06  &  8.984  &   3.683e-06  &  8.798  &   2.045e-01 \\ 
			160  &   1.072e-09  &  9.965  &   1.225e-09  &  10.131  &   3.376e-09  &  10.092  &   7.376e-01 \\ 
			320  &   2.046e-12  &  9.033  &   2.335e-12  &  9.036  &   6.077e-12  &  9.118  &   2.899e+00 \\ 
			Average order    &      &  9.320  &      &  9.384  &      &  9.336  &    \\ 
			\midrule

			\multicolumn{8}{c}{\textbf{WENO11--DeC11}} \\ 
			\midrule
			40  &   7.303e-05  &  $-$  &   8.076e-05  &  $-$  &   1.607e-04  &  $-$  &   1.249e-01 \\ 
			80  &   8.862e-08  &  9.687  &   1.238e-07  &  9.349  &   3.275e-07  &  8.939  &   4.935e-01 \\ 
			160  &   1.604e-11  &  12.431  &   2.150e-11  &  12.492  &   4.997e-11  &  12.678  &   1.855e+00 \\ 
			Average order    &      &  11.059  &      &  10.920  &      &  10.808  &    \\ 
			\midrule

			\multicolumn{8}{c}{\textbf{WENO13--DeC13}} \\ 
			\midrule
			40  &   7.820e-06  &  $-$  &   8.099e-06  &  $-$  &   1.556e-05  &  $-$  &   2.650e-01 \\ 
			80  &   1.621e-09  &  12.236  &   2.404e-09  &  11.718  &   6.710e-09  &  11.179  &   9.964e-01 \\ 
			160  &   4.408e-14  &  15.166  &   5.082e-14  &  15.530  &   1.308e-13  &  15.647  &   3.806e+00 \\ 
			Average order    &      &  13.701  &      &  13.624  &      &  13.413  &    \\ 
			\midrule

			\bottomrule
	\end{tabular}}
\end{table}

\begin{table}[htbp]
	\centering
	\caption{LAE, Test 1: convergence tables for WENO--SSPRK3}
	\label{tab:LAE_1d_convergence_tables_WENO_SSPRK3}
	\scalebox{0.65}{ 
		\begin{tabular}{c c c c c c c c}
			\toprule
			\multirow{2}{*}{$N$} & \multicolumn{2}{c}{$L^1$ error} & \multicolumn{2}{c}{$L^2$ error} & \multicolumn{2}{c}{$L^{\infty}$ error} & \multirow{2}{*}{CPU Time} \\
			\cmidrule(lr){2-3} \cmidrule(lr){4-5} \cmidrule(lr){6-7}
			& Error & Order & Error & Order & Error & Order & \\
			\midrule
			
			\multicolumn{8}{c}{\textbf{WENO3--SSPRK3}} \\ 
			\midrule
			160  &   1.511e-02  &  $-$  &   1.939e-02  &  $-$  &   4.256e-02  &  $-$  &   1.746e-02 \\ 
			320  &   2.292e-03  &  2.721  &   3.855e-03  &  2.331  &   1.190e-02  &  1.839  &   5.784e-02 \\ 
			640  &   2.534e-04  &  3.177  &   5.321e-04  &  2.857  &   2.282e-03  &  2.383  &   1.652e-01 \\ 
			1280  &   1.898e-05  &  3.739  &   3.932e-05  &  3.758  &   2.097e-04  &  3.444  &   6.013e-01 \\ 
			2560  &   1.187e-06  &  3.999  &   1.869e-06  &  4.395  &   8.136e-06  &  4.688  &   2.342e+00 \\ 
			5120  &   7.186e-08  &  4.046  &   8.737e-08  &  4.419  &   2.763e-07  &  4.880  &   9.350e+00 \\ 
			Average order    &      &  3.536  &      &  3.552  &      &  3.447  &    \\ 
			\midrule

			\multicolumn{8}{c}{\textbf{WENO5--SSPRK3}} \\ 
			\midrule
			80  &   2.341e-03  &  $-$  &   1.929e-03  &  $-$  &   2.485e-03  &  $-$  &   5.646e-03 \\ 
			160  &   2.777e-04  &  3.076  &   2.264e-04  &  3.091  &   2.814e-04  &  3.143  &   1.286e-02 \\ 
			320  &   3.500e-05  &  2.988  &   2.804e-05  &  3.014  &   3.421e-05  &  3.040  &   6.353e-02 \\ 
			640  &   4.385e-06  &  2.997  &   3.499e-06  &  3.002  &   4.246e-06  &  3.010  &   1.924e-01 \\ 
			1280  &   5.490e-07  &  2.998  &   4.376e-07  &  2.999  &   5.305e-07  &  3.000  &   6.781e-01 \\ 
			2560  &   6.863e-08  &  3.000  &   5.470e-08  &  3.000  &   6.633e-08  &  3.000  &   2.663e+00 \\ 
			5120  &   8.582e-09  &  2.999  &   6.840e-09  &  3.000  &   8.294e-09  &  2.999  &   1.083e+01 \\ 
			Average order    &      &  3.010  &      &  3.018  &      &  3.032  &    \\ 
			\midrule

			\multicolumn{8}{c}{\textbf{WENO7--SSPRK3}} \\ 
			\midrule
			80  &   2.202e-03  &  $-$  &   1.759e-03  &  $-$  &   2.155e-03  &  $-$  &   1.207e-02 \\ 
			160  &   2.800e-04  &  2.975  &   2.232e-04  &  2.978  &   2.708e-04  &  2.992  &   3.067e-02 \\ 
			320  &   3.506e-05  &  2.997  &   2.794e-05  &  2.998  &   3.389e-05  &  2.999  &   1.309e-01 \\ 
			640  &   4.384e-06  &  2.999  &   3.494e-06  &  2.999  &   4.238e-06  &  2.999  &   4.753e-01 \\ 
			1280  &   5.488e-07  &  2.998  &   4.374e-07  &  2.998  &   5.304e-07  &  2.998  &   1.874e+00 \\ 
			2560  &   6.862e-08  &  3.000  &   5.469e-08  &  3.000  &   6.633e-08  &  3.000  &   7.292e+00 \\ 
			5120  &   8.582e-09  &  2.999  &   6.840e-09  &  2.999  &   8.294e-09  &  2.999  &   2.956e+01 \\ 
			Average order    &      &  2.995  &      &  2.995  &      &  2.998  &    \\ 
			\midrule

			\multicolumn{8}{c}{\textbf{WENO9--SSPRK3}} \\ 
			\midrule
			40  &   1.608e-02  &  $-$  &   1.292e-02  &  $-$  &   1.569e-02  &  $-$  &   4.138e-03 \\ 
			80  &   2.220e-03  &  2.856  &   1.769e-03  &  2.869  &   2.142e-03  &  2.872  &   2.274e-02 \\ 
			160  &   2.800e-04  &  2.987  &   2.232e-04  &  2.987  &   2.705e-04  &  2.986  &   5.450e-02 \\ 
			320  &   3.506e-05  &  2.997  &   2.794e-05  &  2.998  &   3.388e-05  &  2.997  &   2.196e-01 \\ 
			640  &   4.384e-06  &  2.999  &   3.494e-06  &  2.999  &   4.237e-06  &  2.999  &   8.340e-01 \\ 
			1280  &   5.488e-07  &  2.998  &   4.374e-07  &  2.998  &   5.304e-07  &  2.998  &   3.278e+00 \\ 
			2560  &   6.862e-08  &  3.000  &   5.469e-08  &  3.000  &   6.633e-08  &  3.000  &   1.290e+01 \\ 
			5120  &   8.582e-09  &  2.999  &   6.840e-09  &  2.999  &   8.294e-09  &  2.999  &   5.185e+01 \\ 
			Average order    &      &  2.977  &      &  2.978  &      &  2.979  &    \\ 
			\midrule

			\multicolumn{8}{c}{\textbf{WENO11--SSPRK3}} \\ 
			\midrule
			40  &   1.606e-02  &  $-$  &   1.291e-02  &  $-$  &   1.565e-02  &  $-$  &   9.742e-03 \\ 
			80  &   2.220e-03  &  2.855  &   1.769e-03  &  2.868  &   2.142e-03  &  2.869  &   3.570e-02 \\ 
			160  &   2.800e-04  &  2.987  &   2.232e-04  &  2.987  &   2.705e-04  &  2.986  &   9.908e-02 \\ 
			320  &   3.506e-05  &  2.997  &   2.794e-05  &  2.998  &   3.388e-05  &  2.997  &   3.499e-01 \\ 
			640  &   4.384e-06  &  2.999  &   3.494e-06  &  2.999  &   4.237e-06  &  2.999  &   1.352e+00 \\ 
			1280  &   5.488e-07  &  2.998  &   4.374e-07  &  2.998  &   5.304e-07  &  2.998  &   5.355e+00 \\ 
			2560  &   6.862e-08  &  3.000  &   5.469e-08  &  3.000  &   6.633e-08  &  3.000  &   2.116e+01 \\ 
			5120  &   8.582e-09  &  2.999  &   6.840e-09  &  2.999  &   8.294e-09  &  2.999  &   8.444e+01 \\ 
			Average order    &      &  2.977  &      &  2.978  &      &  2.978  &    \\ 
			\midrule

			\multicolumn{8}{c}{\textbf{WENO13--SSPRK3}} \\ 
			\midrule
			40  &   1.606e-02  &  $-$  &   1.290e-02  &  $-$  &   1.564e-02  &  $-$  &   1.473e-02 \\ 
			80  &   2.220e-03  &  2.855  &   1.769e-03  &  2.867  &   2.142e-03  &  2.868  &   3.529e-02 \\ 
			160  &   2.800e-04  &  2.987  &   2.232e-04  &  2.987  &   2.705e-04  &  2.985  &   1.473e-01 \\ 
			320  &   3.506e-05  &  2.997  &   2.794e-05  &  2.998  &   3.388e-05  &  2.997  &   5.327e-01 \\ 
			640  &   4.384e-06  &  2.999  &   3.494e-06  &  2.999  &   4.237e-06  &  2.999  &   2.064e+00 \\ 
			1280  &   5.488e-07  &  2.998  &   4.374e-07  &  2.998  &   5.304e-07  &  2.998  &   8.139e+00 \\ 
			2560  &   6.862e-08  &  3.000  &   5.469e-08  &  3.000  &   6.633e-08  &  3.000  &   3.243e+01 \\ 
			5120  &   8.582e-09  &  2.999  &   6.840e-09  &  2.999  &   8.294e-09  &  2.999  &   1.294e+02 \\ 
			Average order    &      &  2.976  &      &  2.978  &      &  2.978  &    \\ 
			\midrule

			\bottomrule
	\end{tabular}}
\end{table}

\begin{table}[htbp]
	\centering
	\caption{LAE, Test 1: convergence tables for WENO--SSPRK4}
	\label{tab:LAE_1d_convergence_tables_WENO_SSPRK4}
	\scalebox{0.65}{ 
		\begin{tabular}{c c c c c c c c}
			\toprule
			\multirow{2}{*}{$N$} & \multicolumn{2}{c}{$L^1$ error} & \multicolumn{2}{c}{$L^2$ error} & \multicolumn{2}{c}{$L^{\infty}$ error} & \multirow{2}{*}{CPU Time} \\
			\cmidrule(lr){2-3} \cmidrule(lr){4-5} \cmidrule(lr){6-7}
			& Error & Order & Error & Order & Error & Order & \\
			\midrule
			
			\multicolumn{8}{c}{\textbf{WENO3--SSPRK4}} \\ 
			\midrule
			160  &   1.064e-02  &  $-$  &   1.374e-02  &  $-$  &   3.331e-02  &  $-$  &   1.791e-02 \\ 
			320  &   2.057e-03  &  2.371  &   3.494e-03  &  1.975  &   1.098e-02  &  1.602  &   7.795e-02 \\ 
			640  &   2.647e-04  &  2.958  &   5.552e-04  &  2.654  &   2.353e-03  &  2.222  &   2.509e-01 \\ 
			1280  &   1.927e-05  &  3.780  &   3.970e-05  &  3.806  &   2.131e-04  &  3.465  &   9.934e-01 \\ 
			2560  &   1.197e-06  &  4.009  &   1.861e-06  &  4.415  &   8.078e-06  &  4.721  &   3.913e+00 \\ 
			5120  &   7.011e-08  &  4.093  &   8.492e-08  &  4.454  &   2.680e-07  &  4.913  &   1.546e+01 \\ 
			Average order    &      &  3.442  &      &  3.461  &      &  3.385  &    \\ 
			\midrule

			\multicolumn{8}{c}{\textbf{WENO5--SSPRK4}} \\ 
			\midrule
			80  &   9.223e-04  &  $-$  &   8.338e-04  &  $-$  &   1.378e-03  &  $-$  &   8.686e-03 \\ 
			160  &   2.986e-05  &  4.949  &   2.930e-05  &  4.831  &   6.409e-05  &  4.426  &   2.039e-02 \\ 
			320  &   8.166e-07  &  5.193  &   7.764e-07  &  5.238  &   1.830e-06  &  5.130  &   8.554e-02 \\ 
			640  &   2.268e-08  &  5.170  &   1.975e-08  &  5.297  &   3.711e-08  &  5.624  &   2.970e-01 \\ 
			1280  &   9.976e-10  &  4.507  &   8.163e-10  &  4.597  &   9.736e-10  &  5.253  &   1.119e+00 \\ 
			2560  &   5.895e-11  &  4.081  &   4.699e-11  &  4.119  &   5.196e-11  &  4.228  &   4.393e+00 \\ 
			5120  &   3.640e-12  &  4.017  &   2.873e-12  &  4.032  &   3.134e-12  &  4.051  &   1.771e+01 \\ 
			Average order    &      &  4.653  &      &  4.685  &      &  4.785  &    \\ 
			\midrule

			\multicolumn{8}{c}{\textbf{WENO7--SSPRK4}} \\ 
			\midrule
			80  &   1.337e-04  &  $-$  &   1.309e-04  &  $-$  &   2.626e-04  &  $-$  &   1.340e-02 \\ 
			160  &   3.869e-06  &  5.111  &   3.038e-06  &  5.429  &   3.248e-06  &  6.337  &   6.106e-02 \\ 
			320  &   2.399e-07  &  4.011  &   1.887e-07  &  4.009  &   2.039e-07  &  3.994  &   2.076e-01 \\ 
			640  &   1.496e-08  &  4.003  &   1.179e-08  &  4.000  &   1.275e-08  &  3.999  &   7.456e-01 \\ 
			1280  &   9.360e-10  &  3.999  &   7.380e-10  &  3.998  &   7.977e-10  &  3.998  &   2.931e+00 \\ 
			2560  &   5.851e-11  &  4.000  &   4.613e-11  &  4.000  &   4.984e-11  &  4.000  &   1.137e+01 \\ 
			5120  &   3.635e-12  &  4.008  &   2.864e-12  &  4.009  &   3.078e-12  &  4.017  &   4.557e+01 \\ 
			Average order    &      &  4.189  &      &  4.241  &      &  4.391  &    \\ 
			\midrule

			\multicolumn{8}{c}{\textbf{WENO9--SSPRK4}} \\ 
			\midrule
			40  &   1.064e-03  &  $-$  &   8.501e-04  &  $-$  &   1.131e-03  &  $-$  &   6.474e-03 \\ 
			80  &   6.081e-05  &  4.130  &   4.816e-05  &  4.142  &   5.217e-05  &  4.438  &   3.621e-02 \\ 
			160  &   3.824e-06  &  3.991  &   3.017e-06  &  3.997  &   3.260e-06  &  4.000  &   1.010e-01 \\ 
			320  &   2.393e-07  &  3.998  &   1.887e-07  &  3.999  &   2.039e-07  &  3.999  &   3.583e-01 \\ 
			640  &   1.496e-08  &  4.000  &   1.179e-08  &  4.000  &   1.275e-08  &  3.999  &   1.369e+00 \\ 
			1280  &   9.360e-10  &  3.998  &   7.380e-10  &  3.998  &   7.977e-10  &  3.998  &   5.439e+00 \\ 
			2560  &   5.851e-11  &  4.000  &   4.613e-11  &  4.000  &   4.984e-11  &  4.000  &   2.167e+01 \\ 
			5120  &   3.635e-12  &  4.008  &   2.864e-12  &  4.009  &   3.078e-12  &  4.017  &   8.558e+01 \\ 
			Average order    &      &  4.018  &      &  4.021  &      &  4.065  &    \\ 
			\midrule

			\multicolumn{8}{c}{\textbf{WENO11--SSPRK4}} \\ 
			\midrule
			40  &   9.390e-04  &  $-$  &   7.408e-04  &  $-$  &   7.741e-04  &  $-$  &   1.561e-02 \\ 
			80  &   6.110e-05  &  3.942  &   4.816e-05  &  3.943  &   5.206e-05  &  3.894  &   5.439e-02 \\ 
			160  &   3.824e-06  &  3.998  &   3.017e-06  &  3.997  &   3.260e-06  &  3.997  &   1.575e-01 \\ 
			320  &   2.393e-07  &  3.998  &   1.887e-07  &  3.999  &   2.039e-07  &  3.999  &   5.784e-01 \\ 
			640  &   1.496e-08  &  4.000  &   1.179e-08  &  4.000  &   1.275e-08  &  3.999  &   2.261e+00 \\ 
			1280  &   9.360e-10  &  3.998  &   7.380e-10  &  3.998  &   7.977e-10  &  3.998  &   8.962e+00 \\ 
			2560  &   5.851e-11  &  4.000  &   4.613e-11  &  4.000  &   4.984e-11  &  4.000  &   3.585e+01 \\ 
			5120  &   3.635e-12  &  4.008  &   2.864e-12  &  4.009  &   3.078e-12  &  4.017  &   1.433e+02 \\ 
			Average order    &      &  3.992  &      &  3.992  &      &  3.987  &    \\ 
			\midrule

			\multicolumn{8}{c}{\textbf{WENO13--SSPRK4}} \\ 
			\midrule
			40  &   9.160e-04  &  $-$  &   7.294e-04  &  $-$  &   7.747e-04  &  $-$  &   2.454e-02 \\ 
			80  &   6.107e-05  &  3.907  &   4.816e-05  &  3.921  &   5.207e-05  &  3.895  &   7.745e-02 \\ 
			160  &   3.824e-06  &  3.997  &   3.017e-06  &  3.997  &   3.260e-06  &  3.997  &   2.421e-01 \\ 
			320  &   2.393e-07  &  3.998  &   1.887e-07  &  3.999  &   2.039e-07  &  3.999  &   9.047e-01 \\ 
			640  &   1.496e-08  &  4.000  &   1.179e-08  &  4.000  &   1.275e-08  &  3.999  &   3.555e+00 \\ 
			1280  &   9.360e-10  &  3.998  &   7.380e-10  &  3.998  &   7.977e-10  &  3.998  &   1.407e+01 \\ 
			2560  &   5.851e-11  &  4.000  &   4.613e-11  &  4.000  &   4.984e-11  &  4.000  &   5.625e+01 \\ 
			5120  &   3.636e-12  &  4.008  &   2.864e-12  &  4.009  &   3.079e-12  &  4.017  &   2.276e+02 \\ 
			Average order    &      &  3.987  &      &  3.989  &      &  3.987  &    \\ 
			\midrule

			\bottomrule
	\end{tabular}}
\end{table}

\begin{table}[htbp]
	\centering
	\caption{LAE, Test 1: convergence tables for WENO--mSSPRK3}
	\label{tab:LAE_1d_convergence_tables_WENO_mSSPRK3}
	\scalebox{0.65}{ 
		\begin{tabular}{c c c c c c c c}
			\toprule
			\multirow{2}{*}{$N$} & \multicolumn{2}{c}{$L^1$ error} & \multicolumn{2}{c}{$L^2$ error} & \multicolumn{2}{c}{$L^{\infty}$ error} & \multirow{2}{*}{CPU Time} \\
			\cmidrule(lr){2-3} \cmidrule(lr){4-5} \cmidrule(lr){6-7}
			& Error & Order & Error & Order & Error & Order & \\
			\midrule
			
			\multicolumn{8}{c}{\textbf{WENO5--mSSPRK3}} \\ 
			\midrule
			80  &   9.297e-04  &  $-$  &   8.538e-04  &  $-$  &   1.524e-03  &  $-$  &   3.390e-02 \\ 
			160  &   2.955e-05  &  4.975  &   2.945e-05  &  4.858  &   6.656e-05  &  4.517  &   2.187e-01 \\ 
			320  &   7.724e-07  &  5.258  &   7.528e-07  &  5.290  &   1.754e-06  &  5.246  &   1.268e+00 \\ 
			640  &   1.821e-08  &  5.407  &   1.574e-08  &  5.580  &   2.864e-08  &  5.936  &   7.856e+00 \\ 
			1280  &   4.171e-10  &  5.448  &   3.418e-10  &  5.525  &   4.629e-10  &  5.951  &   4.968e+01 \\ 
			2560  &   1.398e-11  &  4.899  &   1.322e-11  &  4.692  &   2.024e-11  &  4.516  &   3.168e+02 \\ 
			5120  &   3.892e-11  &  -1.477  &   3.417e-11  &  -1.370  &   3.990e-11  &  -0.979  &   1.954e+03 \\ 
			Average order    &      &  4.085  &      &  4.096  &      &  4.198  &    \\ 
			\midrule

			\multicolumn{8}{c}{\textbf{WENO7--mSSPRK3}} \\ 
			\midrule
			80  &   1.027e-04  &  $-$  &   1.188e-04  &  $-$  &   2.185e-04  &  $-$  &   9.909e-01 \\ 
			160  &   3.151e-07  &  8.348  &   3.967e-07  &  8.227  &   9.737e-07  &  7.810  &   9.597e+00 \\ 
			320  &   2.211e-09  &  7.155  &   4.327e-09  &  6.518  &   2.346e-08  &  5.376  &   9.351e+01 \\ 
			640  &   5.900e-11  &  5.228  &   5.761e-11  &  6.231  &   2.620e-10  &  6.484  &   9.426e+02 \\ 
			1280  &   2.582e-10  &  -2.129  &   2.267e-10  &  -1.976  &   2.635e-10  &  -0.008  &   9.549e+03 \\ 
			Average order    &      &  4.650  &      &  4.750  &      &  4.915  &    \\ 
			\midrule

			\multicolumn{8}{c}{\textbf{WENO9--mSSPRK3}} \\ 
			\midrule
			40  &   5.391e-04  &  $-$  &   7.010e-04  &  $-$  &   1.649e-03  &  $-$  &   1.345e+00 \\ 
			80  &   1.071e-06  &  8.975  &   1.374e-06  &  8.995  &   3.684e-06  &  8.806  &   2.075e+01 \\ 
			160  &   1.072e-09  &  9.965  &   1.226e-09  &  10.131  &   3.376e-09  &  10.092  &   3.286e+02 \\ 
			320  &   3.142e-10  &  1.771  &   2.758e-10  &  2.152  &   3.208e-10  &  3.396  &   5.268e+03 \\ 
			640  &   2.571e-09  &  -3.033  &   2.258e-09  &  -3.033  &   2.621e-09  &  -3.031  &   8.259e+04 \\ 
			Average order    &      &  4.419  &      &  4.561  &      &  4.816  &    \\ 
			\midrule

			\multicolumn{8}{c}{\textbf{WENO11--mSSPRK3}} \\ 
			\midrule
			40  &   7.303e-05  &  $-$  &   8.076e-05  &  $-$  &   1.607e-04  &  $-$  &   1.735e+01 \\ 
			80  &   8.862e-08  &  9.687  &   1.238e-07  &  9.349  &   3.275e-07  &  8.939  &   4.263e+02 \\ 
			160  &   5.851e-10  &  7.243  &   5.082e-10  &  7.929  &   5.904e-10  &  9.116  &   1.073e+04 \\ 
			320  &   5.076e-09  &  -3.117  &   4.457e-09  &  -3.133  &   5.187e-09  &  -3.135  &   2.502e+05 \\ 
			Average order    &      &  4.604  &      &  4.715  &      &  4.973  &    \\ 
			\midrule

			\multicolumn{8}{c}{\textbf{WENO13--mSSPRK3}} \\ 
			\midrule
			40  &   7.820e-06  &  $-$  &   8.099e-06  &  $-$  &   1.556e-05  &  $-$  &   1.809e+02 \\ 
			80  &   1.800e-09  &  12.085  &   2.421e-09  &  11.708  &   6.710e-09  &  11.179  &   7.095e+03 \\ 
			160  &   2.000e-08  &  -3.474  &   1.755e-08  &  -2.858  &   2.039e-08  &  -1.604  &   2.813e+05 \\ 
			Average order    &      &  4.306  &      &  4.425  &      &  4.787  &    \\ 
			\midrule

			\bottomrule
	\end{tabular}}
\end{table}

\begin{table}[htbp]
	\centering
	\caption{LAE, Test 1: convergence tables for WENO--mSSPRK4}
	\label{tab:LAE_1d_convergence_tables_WENO_mSSPRK4}
	\scalebox{0.65}{ 
		\begin{tabular}{c c c c c c c c}
			\toprule
			\multirow{2}{*}{$N$} & \multicolumn{2}{c}{$L^1$ error} & \multicolumn{2}{c}{$L^2$ error} & \multicolumn{2}{c}{$L^{\infty}$ error} & \multirow{2}{*}{CPU Time} \\
			\cmidrule(lr){2-3} \cmidrule(lr){4-5} \cmidrule(lr){6-7}
			& Error & Order & Error & Order & Error & Order & \\
			\midrule
			
			\multicolumn{8}{c}{\textbf{WENO5--mSSPRK4}} \\ 
			\midrule
			80  &   9.291e-04  &  $-$  &   8.537e-04  &  $-$  &   1.521e-03  &  $-$  &   1.936e-02 \\ 
			160  &   2.953e-05  &  4.976  &   2.945e-05  &  4.857  &   6.657e-05  &  4.514  &   5.633e-02 \\ 
			320  &   7.717e-07  &  5.258  &   7.527e-07  &  5.290  &   1.753e-06  &  5.247  &   2.707e-01 \\ 
			640  &   1.819e-08  &  5.407  &   1.573e-08  &  5.580  &   2.863e-08  &  5.937  &   1.180e+00 \\ 
			1280  &   4.163e-10  &  5.449  &   3.412e-10  &  5.527  &   4.622e-10  &  5.953  &   5.552e+00 \\ 
			2560  &   9.870e-12  &  5.398  &   8.642e-12  &  5.303  &   1.299e-11  &  5.153  &   2.594e+01 \\ 
			5120  &   1.623e-12  &  2.604  &   1.431e-12  &  2.595  &   1.844e-12  &  2.817  &   1.235e+02 \\ 
			Average order    &      &  4.849  &      &  4.859  &      &  4.937  &    \\ 
			\midrule

			\multicolumn{8}{c}{\textbf{WENO7--mSSPRK4}} \\ 
			\midrule
			80  &   1.027e-04  &  $-$  &   1.188e-04  &  $-$  &   2.185e-04  &  $-$  &   1.964e-01 \\ 
			160  &   3.151e-07  &  8.348  &   3.967e-07  &  8.227  &   9.737e-07  &  7.810  &   1.253e+00 \\ 
			320  &   2.211e-09  &  7.155  &   4.327e-09  &  6.518  &   2.346e-08  &  5.376  &   8.207e+00 \\ 
			640  &   1.458e-11  &  7.245  &   3.717e-11  &  6.863  &   2.627e-10  &  6.480  &   5.463e+01 \\ 
			1280  &   3.661e-12  &  1.994  &   3.203e-12  &  3.536  &   3.728e-12  &  6.139  &   3.645e+02 \\ 
			Average order    &      &  6.185  &      &  6.286  &      &  6.451  &    \\ 
			\midrule

			\multicolumn{8}{c}{\textbf{WENO9--mSSPRK4}} \\ 
			\midrule
			40  &   5.391e-04  &  $-$  &   7.010e-04  &  $-$  &   1.649e-03  &  $-$  &   2.537e-01 \\ 
			80  &   1.071e-06  &  8.975  &   1.374e-06  &  8.995  &   3.684e-06  &  8.806  &   2.196e+00 \\ 
			160  &   1.072e-09  &  9.965  &   1.225e-09  &  10.131  &   3.375e-09  &  10.092  &   2.059e+01 \\ 
			320  &   1.096e-11  &  6.612  &   9.109e-12  &  7.072  &   1.066e-11  &  8.307  &   1.934e+02 \\ 
			640  &   3.944e-11  &  -1.847  &   3.463e-11  &  -1.927  &   4.024e-11  &  -1.917  &   1.824e+03 \\ 
			Average order    &      &  5.926  &      &  6.068  &      &  6.322  &    \\ 
			\midrule

			\multicolumn{8}{c}{\textbf{WENO11--mSSPRK4}} \\ 
			\midrule
			40  &   7.303e-05  &  $-$  &   8.076e-05  &  $-$  &   1.607e-04  &  $-$  &   1.721e+00 \\ 
			80  &   8.862e-08  &  9.687  &   1.238e-07  &  9.349  &   3.275e-07  &  8.939  &   2.250e+01 \\ 
			160  &   1.652e-11  &  12.389  &   2.153e-11  &  12.490  &   4.998e-11  &  12.678  &   2.978e+02 \\ 
			320  &   7.282e-11  &  -2.140  &   6.394e-11  &  -1.570  &   7.430e-11  &  -0.572  &   3.967e+03 \\ 
			Average order    &      &  6.645  &      &  6.756  &      &  7.015  &    \\ 
			\midrule

			\multicolumn{8}{c}{\textbf{WENO13--mSSPRK4}} \\ 
			\midrule
			40  &   7.820e-06  &  $-$  &   8.099e-06  &  $-$  &   1.556e-05  &  $-$  &   1.225e+01 \\ 
			80  &   1.621e-09  &  12.236  &   2.404e-09  &  11.718  &   6.709e-09  &  11.179  &   2.259e+02 \\ 
			160  &   1.281e-10  &  3.661  &   1.124e-10  &  4.418  &   1.306e-10  &  5.683  &   4.233e+03 \\ 
			Average order    &      &  7.949  &      &  8.068  &      &  8.431  &    \\ 
			\midrule

			\bottomrule
	\end{tabular}}
\end{table}

\begin{figure}[htbp]
	\centering
	\begin{subfigure}[b]{0.45\textwidth}
		\centering
		\includegraphics[width=\textwidth]{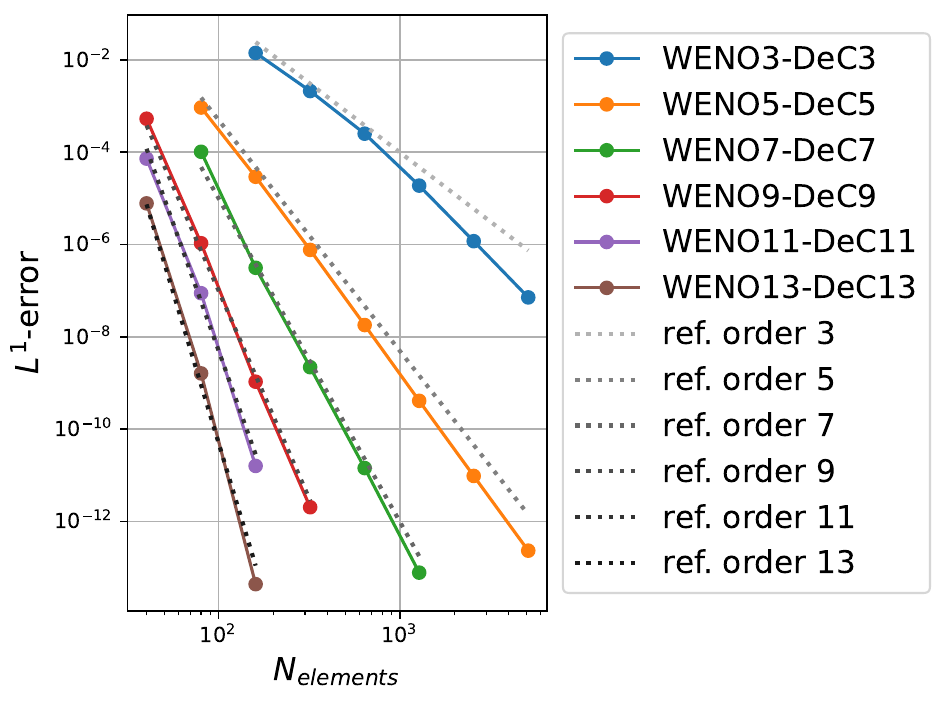}
		\caption{Convergence analysis for WENO--DeC}
	\end{subfigure}
	\quad
	\begin{subfigure}[b]{0.45\textwidth}
		\centering
		\includegraphics[width=\textwidth]{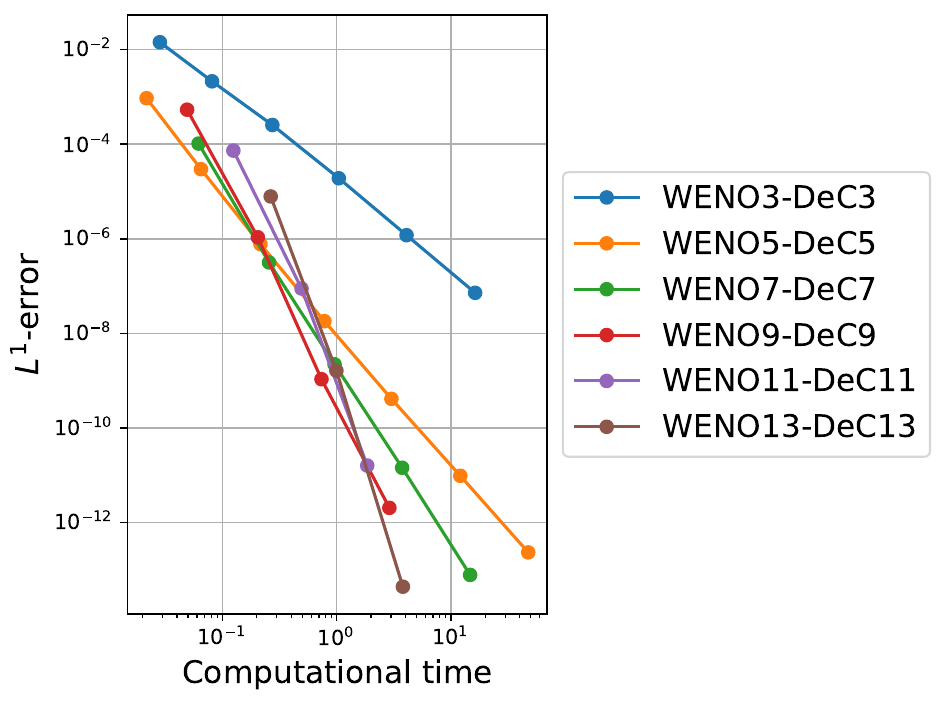}
		\caption{Efficiency analysis for WENO--DeC}
	\end{subfigure}
	\caption{LAE, Test 1: Results obtained with WENO--DeC}
	\label{fig:LAE_1d_sin4_WENODeC}
\end{figure}

A visual representation of the behavior of the error in the $L^1$-norm is reported in Figure~\ref{fig:LAE_1d_sin4_WENODeC}.
In particular, in the left plot, we report the error with respect to the number of elements, while, in the right one we report the error with respect to the computational time. In both cases, data are represented on logarithmic scales.
From the left plot, one can see how the expected rate of convergence is always achieved.
The expected trend is also recovered in the right plot, with high order methods being asymptotically more convenient.
The errors in the $L^2$- and $L^\infty$-norms behave similarly, hence, the associated plots are omitted.

Convergence tables for WENO--SSPRK3 and WENO--SSPRK4 are reported in Tables~\ref{tab:LAE_1d_convergence_tables_WENO_SSPRK3} and \ref{tab:LAE_1d_convergence_tables_WENO_SSPRK4} respectively.
Also for such settings, the errors in the three different norms display similar behaviors amongst themselves. From the convergence trends, it is immediately evident how the employment of lower order time discretizations, for very high order space discretizations, results in order degradation.
As expected, all WENO space discretizations from order 5 on asymptotically converge with order 3 and 4 when combined with SSPRK3 and SSPRK4 respectively.
Instead, also in this case, WENO3 is characterized by the same asymptotical superconvergence observed with DeC time integration.
\begin{remark}[On the asymptotic character of convergence]\label{rmk:asymptotic_convergence}
	The formal order of accuracy is related to the asymptotic behavior of the error with respect to the mesh refinement.
	When coupling a high order spatial discretization with a lower order time discretization, it may take some levels of refinement for the error in time to dominate with respect to the one in space and, hence, for the resulting approach to achieve the formal convergence rate determined by the time discretization.
	This is the case for WENO5--SSPRK4, which is characterized by a 5-th order convergence trend until a very high level of refinement, see Table~\ref{tab:LAE_1d_convergence_tables_WENO_SSPRK4}.
	From 2560 elements on, the expected 4-th order is finally achieved.
	Depending on the specific problem, one may not see order 4 for such a method before hitting machine precision.
\end{remark}

\begin{figure}[htbp]
	\centering
	\begin{subfigure}[b]{0.45\textwidth}
		\centering
		\includegraphics[width=\textwidth]{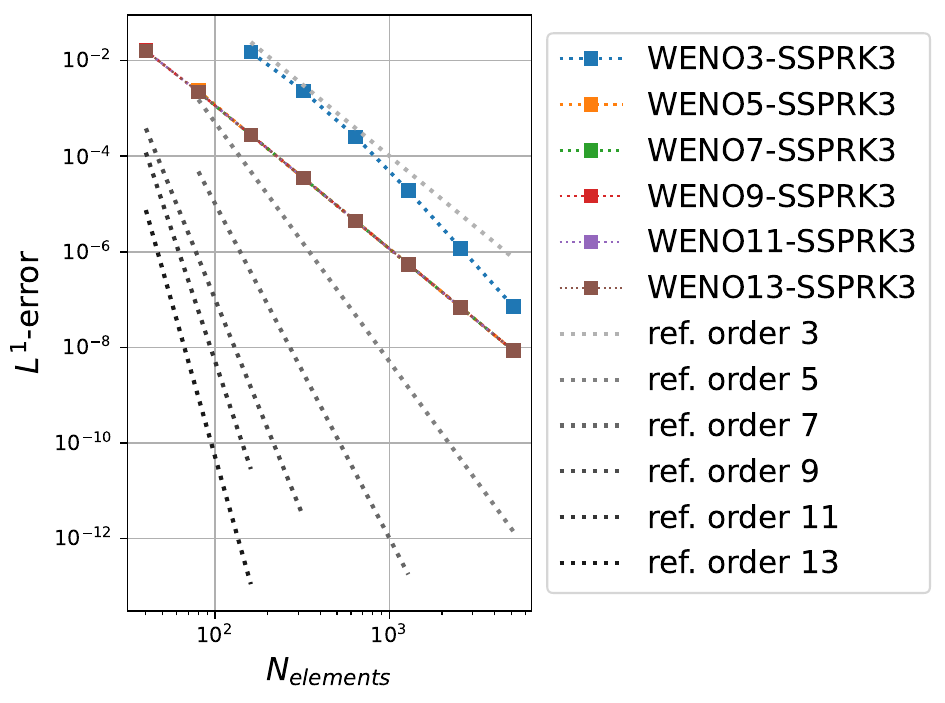}
		\caption{Convergence analysis for WENO--SSPRK3}
	\end{subfigure}
	\quad
	\begin{subfigure}[b]{0.45\textwidth}
		\centering
		\includegraphics[width=\textwidth]{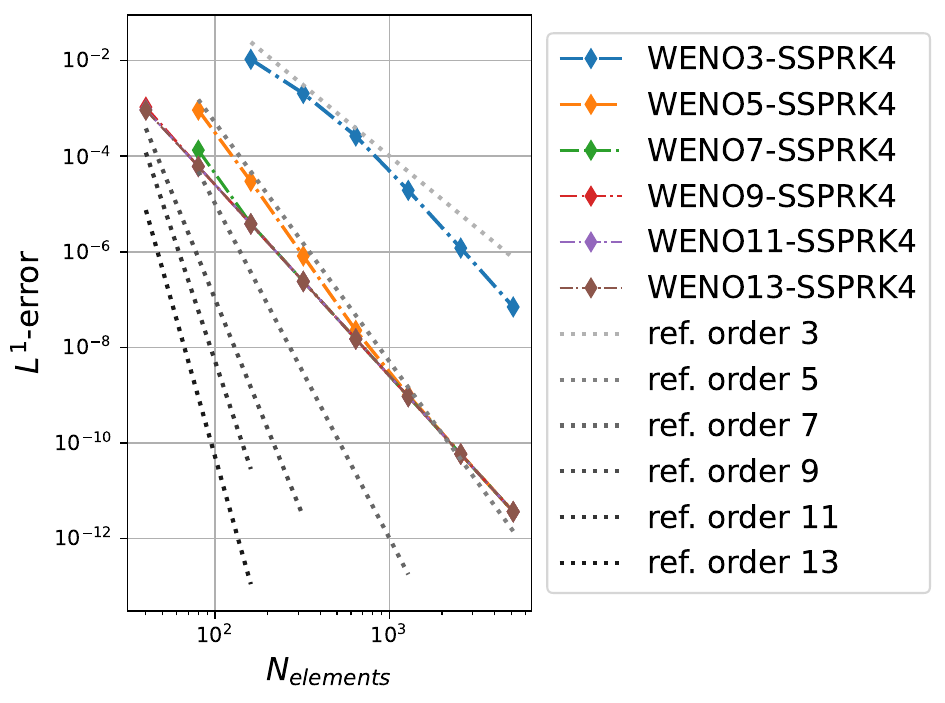}
		\caption{Convergence analysis for WENO--SSPRK4}
	\end{subfigure}\\
	\begin{subfigure}[b]{0.45\textwidth}
		\centering
		\includegraphics[width=\textwidth]{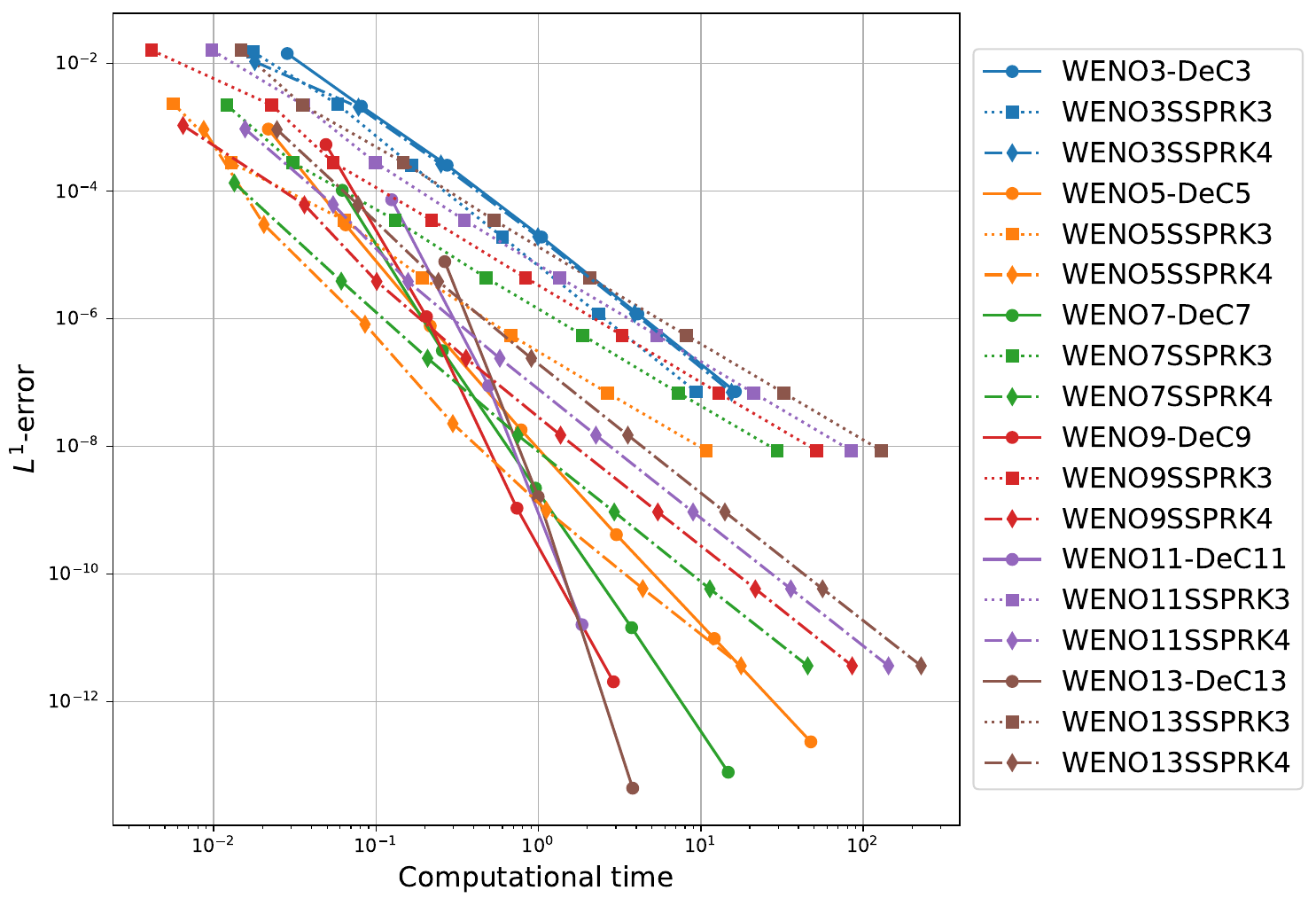}
		\caption{Efficiency analysis for WENO--SSPRK3, WENO--SSPRK4 and WENO--DeC}
	\end{subfigure}
	\caption{LAE, Test 1: Results obtained with WENO--SSPRK3 and WENO--SSPRK4 on top and efficiency comparison with WENO--DeC on the bottom}
	\label{fig:LAE_1d_sin4_WENOSSPRK}
\end{figure}

In Figure~\ref{fig:LAE_1d_sin4_WENOSSPRK}, we report the results of convergence analyses of WENO--SSPRK3 and WENO--SSPRK4, along with the efficiency comparison between such couplings and WENO--DeC, in the $L^1$-norm.
The efficiency comparison on the bottom makes clear that, although for coarse meshes and, hence, large errors SSPRK methods can be better than DeC, for very high order methods and higher level of refinements adopting lower order time discretizations implies losses in computational efficiency.

\begin{figure}[htbp]
	\centering
	\begin{subfigure}[b]{0.55\textwidth}
		\centering
		\includegraphics[width=1.3\linewidth, trim={0 710 0 0}, clip]{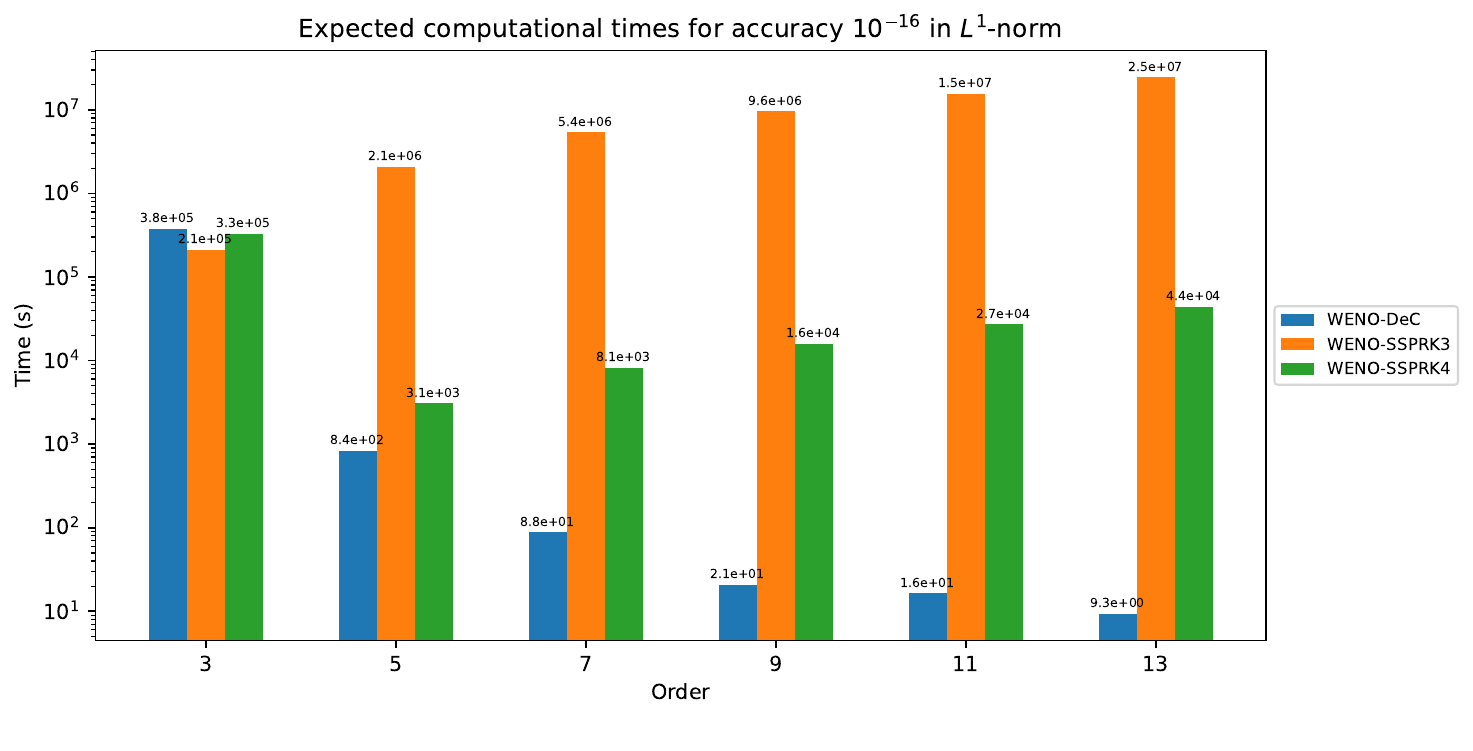}
		\caption{Results for $L^1$-norm}
	\end{subfigure}
	\\
	\begin{subfigure}[b]{0.55\textwidth}
		\centering
		\includegraphics[width=1.3\linewidth, trim={0 360 0 360}, clip]{figures/LAE1D/sin4/advection_sin4_longer_domain_paper_LAE_1D_bar_plots_expected_times.pdf}
		\caption{Results for $L^2$-norm}
	\end{subfigure}
	\\
	\begin{subfigure}[b]{0.55\textwidth}
		\centering
		\includegraphics[width=1.3\linewidth, trim={0 0 0 700}, clip]{figures/LAE1D/sin4/advection_sin4_longer_domain_paper_LAE_1D_bar_plots_expected_times.pdf}
		\caption{Results for $L^{\infty}$-norm}
	\end{subfigure}
	\caption{LAE, Test 1: Expected computational times in seconds to reach an accuracy level equal to $10^{-16}$ in the $L^1$-, $L^2$-, and $L^{\infty}$-norms.}
	\label{fig:expected_time_LAE_1d}
\end{figure}

In order to give a quantitative meaning to the efficiency loss obtained adopting SSPRK3 and SSPRK4 with very high order space discretizations, we provide in Figure~\ref{fig:expected_time_LAE_1d} the estimated computational times needed to reach an accuracy tolerance equal to $10^{-16}$ with WENO--DeC, WENO--SSPRK3 and WENO--SSPRK4, in the three norms $L^1$, $L^2$-and $L^{\infty}$.
They are computed considering the linear regression of the last three points of each curve error versus time in logarithmic scale in order to only focus on the asymptotic behavior and exclude possible underestimates of the required computational time determined by initial superconvergences occurring for coarse meshes, for which the expected rate has not been achieved yet, see for example Remark~\ref{rmk:asymptotic_convergence}.
We see that, while in the context of WENO--DeC increasing the order is always advantageous resulting in smaller computational times, for WENO--SSPRK3 and WENO--SSPRK4 this is not the case. 
In fact, since the error in time dominates for very high spatial orders, increasing the order of the space discretization essentially results in an increase of the needed computational time and consequently in a waste of resources.
The only exception is represented by WENO3--SSPRK4, which requires more computational time than WENO5--SSPRK4.
However, this is not completely unexpected, as WENO5--SSPRK4 is formally 4-th order accurate, while, WENO3--SSPRK4 has formal order of accuracy equal to 3 only.

\begin{figure}[htbp]
	\centering
	\begin{subfigure}[b]{0.45\textwidth}
		\centering
		\includegraphics[width=\textwidth]{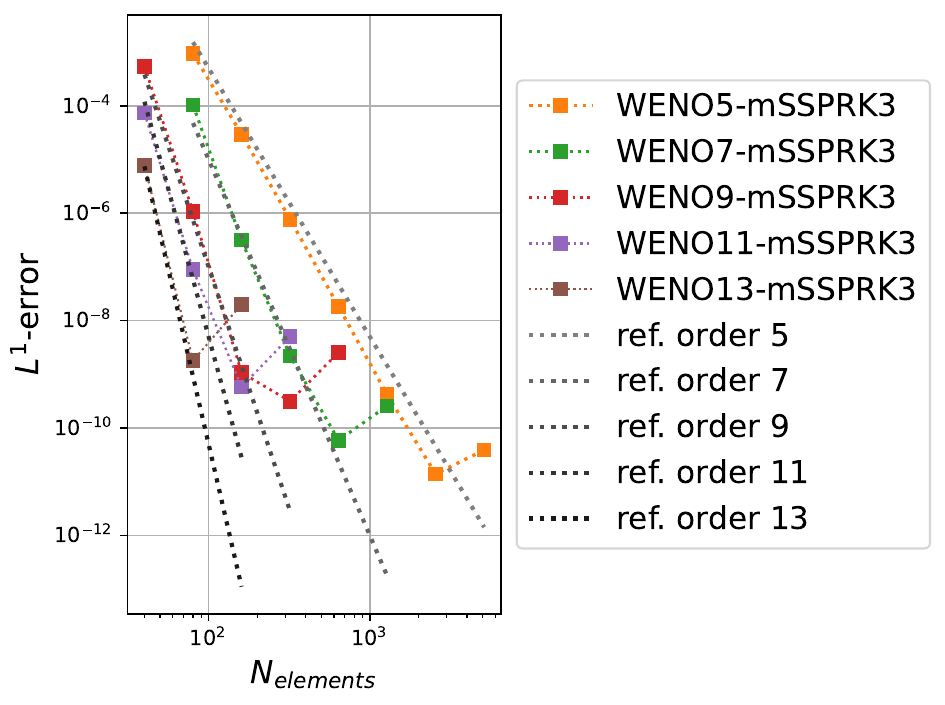}
		\caption{Convergence analysis for WENO--mSSPRK3}
	\end{subfigure}
	\quad
	\begin{subfigure}[b]{0.45\textwidth}
		\centering
		\includegraphics[width=\textwidth]{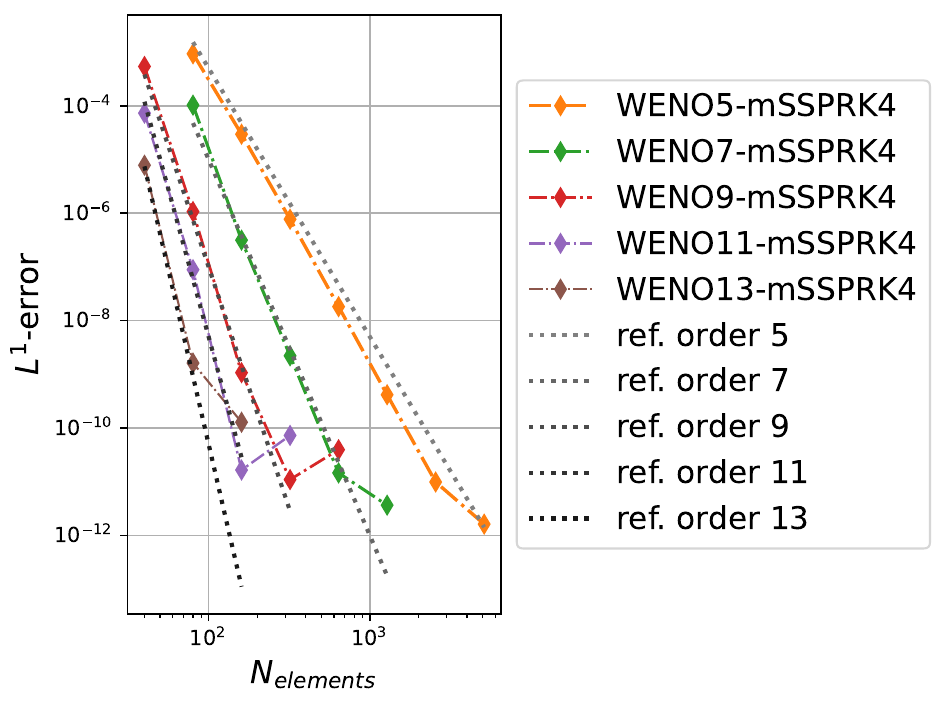}
		\caption{Convergence analysis for WENO--mSSPRK4}
	\end{subfigure}\\
	\begin{subfigure}[b]{0.45\textwidth}
		\centering
		\includegraphics[width=\textwidth]{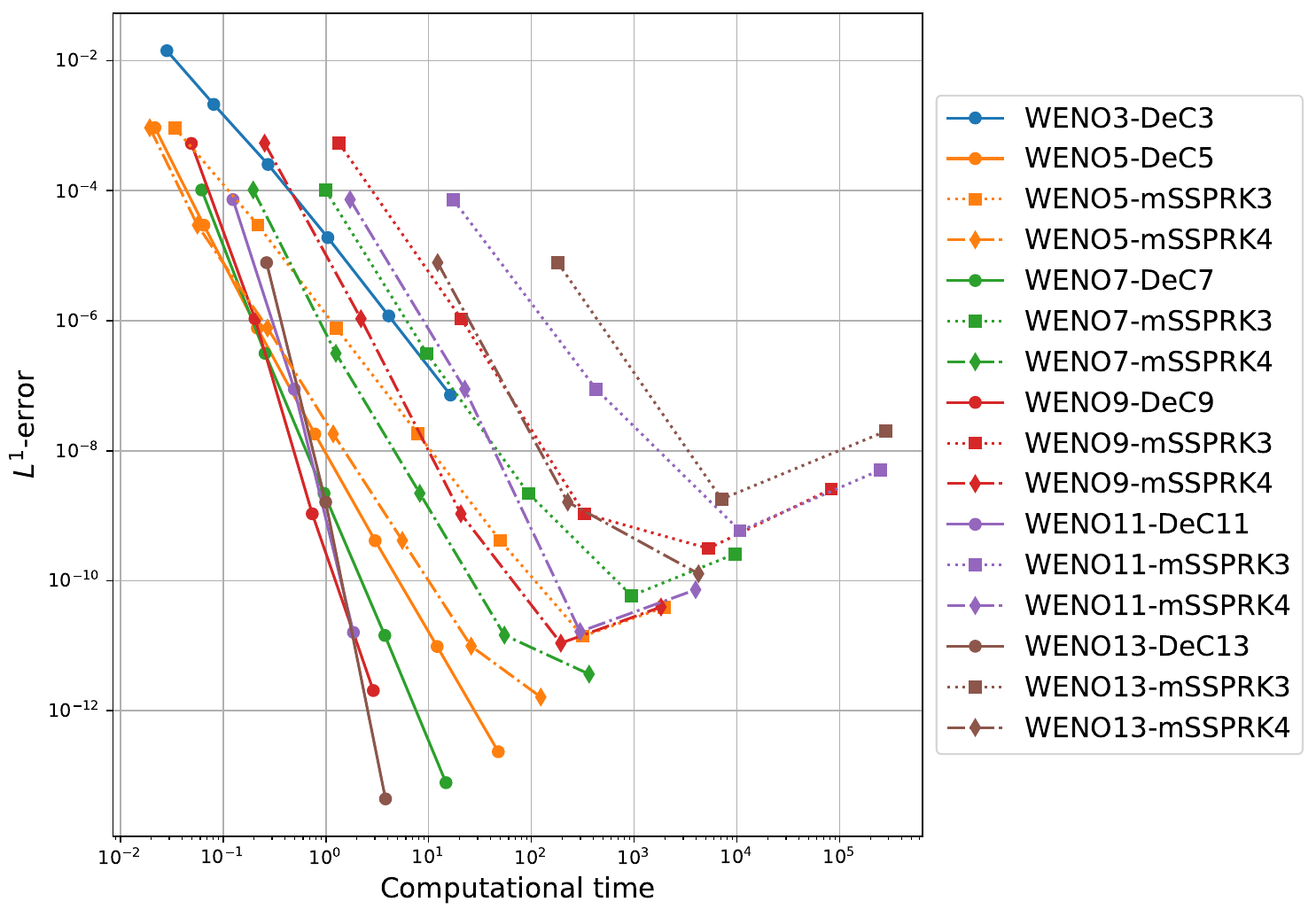}
		\caption{Efficiency analysis for WENO--mSSPRK3, WENO--mSSPRK4 and WENO--DeC}
	\end{subfigure}
	\caption{LAE, Test 1: Results obtained with WENO--mSSPRK3 and WENO--mSSPRK4 on top and efficiency comparison with WENO--DeC on the bottom}
	\label{fig:LAE_1d_sin4_WENOmSSPRK}
\end{figure}

In Tables~\ref{tab:LAE_1d_convergence_tables_WENO_mSSPRK3} and~\ref{tab:LAE_1d_convergence_tables_WENO_mSSPRK4}, the results of the convergence analyses obtained with WENO--mSSPRK3 and WENO--mSSPRK4 are reported. 
As one can notice, not always the expected rate of convergence has been achieved, and there is a strong tendency to order degradation in the mesh refinement. This is due to the time step restrictions, which characterize the methods. The excessive (but necessary in order not to spoil the spatial order of accuracy) reduction in the time step determines a huge number of computations to be performed, which are associated with excessive accumulation of rounding errors, leading to machine precision effects arising far above the usual values.
Graphical representations of the convergence and efficiency analyses of WENO--mSSPRK3 and WENO--mSSPRK4 in $L^1$--norm are reported in Figure~\ref{fig:LAE_1d_sin4_WENOmSSPRK}.
In particular, from the efficiency comparison at the bottom, one can observe that such methods are never more convenient than WENO--DeC. It is also evident that they require a huge amount of computational time in the mesh refinement, much longer than what required by WENO--DeC.

In conclusion, the best results among the considered schemes are the ones obtained with WENO--DeC.

\subsubsection{Test 2: Composite wave}\label{sec:LAE_1d_composite_wave}
Here, we consider a more challenging test consisting in the propagation of a wave made by signals with various shapes~\cite{jiang1996efficient}.
More in detail, we consider the computational domain $\Omega:=[-1,1]$ with periodic boundary conditions and initial condition given by
\begin{align}
	u_0(x):=\begin{cases}
				\frac{1}{6}[G(x,\beta,z-\delta)+4 G(x,\beta,z)+G(x,\beta,z+\delta)],\quad &\text{if} \quad -0.8 \leq x \leq -0.6,\\
				1,\quad &\text{if} \quad -0.4 \leq x \leq -0.2,\\
				1.0-\abs{10(x-0.1)},\quad &\text{if} \quad 0 \leq x \leq 0.2,\\
				\frac{1}{6}[F(x,\alpha,a-\delta)+4F(x,\alpha,a)+F(x,\alpha,a+\delta)],\quad &\text{if} \quad 0.4 \leq x \leq 0.6,\\
				0, \quad &\text{otherwise},
			\end{cases}
\end{align}
with $\delta:=0.005$, $z:=-0.7$, $a:=0.5$, $\alpha:=10$, $\beta:=\frac{\log{(2)}}{36\delta^2}$ and functions $G$ and $F$ defined as
\begin{align}
	G(x,\beta,z):=e^{-\beta(x-z)^2}, \quad F(x,\alpha,a):=\sqrt{\max{\left(1-\alpha^2(x-a)^2,0\right)}}.
\end{align}

\begin{figure}[htbp]
	\centering
	\begin{subfigure}[b]{0.45\textwidth}
		\centering
		\includegraphics[width=\textwidth]{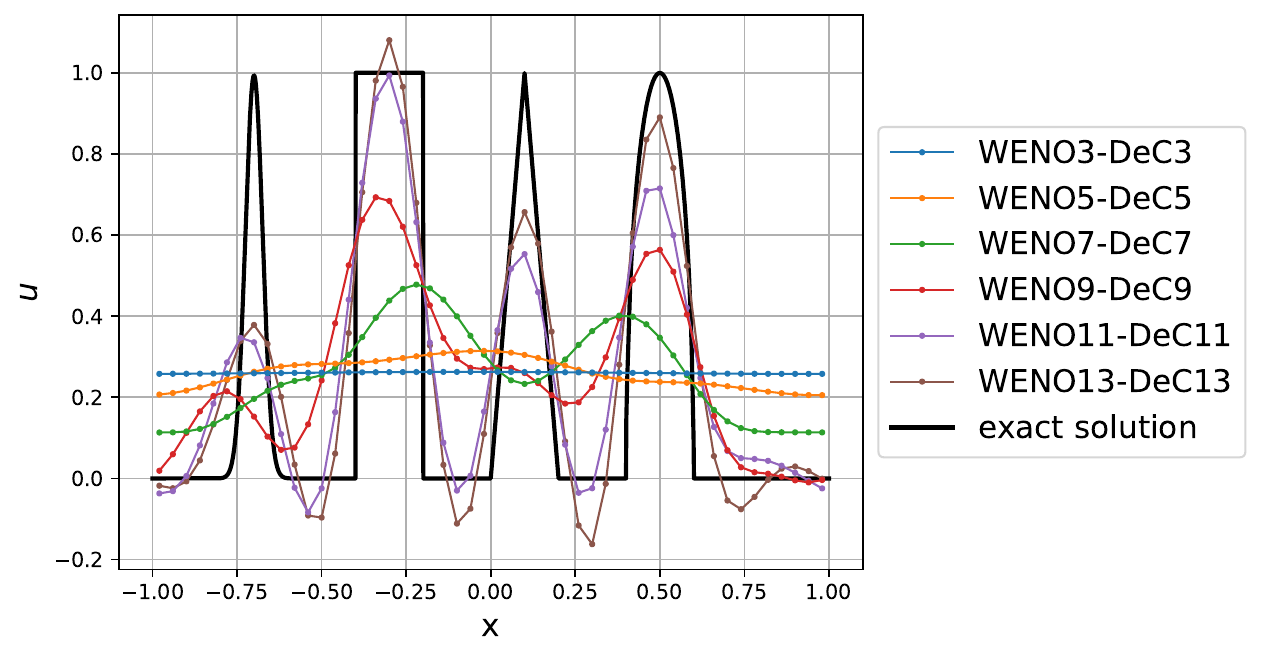}
		\caption{50 elements}
	\end{subfigure}
	\begin{subfigure}[b]{0.45\textwidth}
		\centering
		\includegraphics[width=\textwidth]{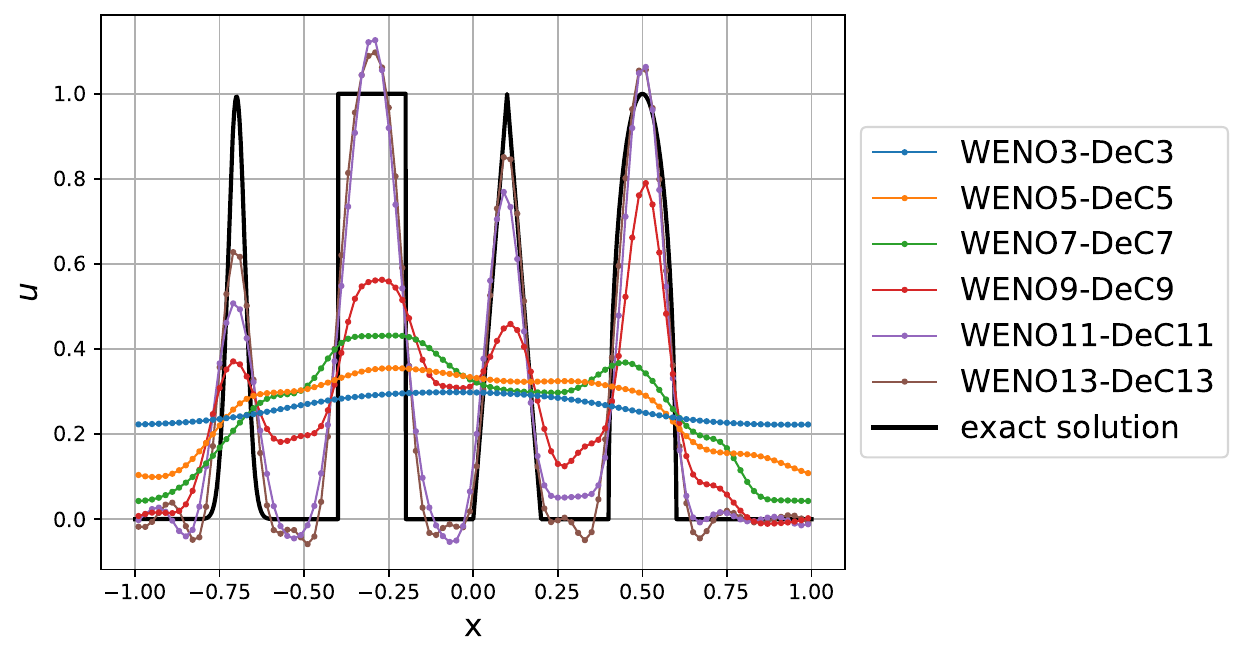}
		\caption{100 elements}
	\end{subfigure}\\
	\begin{subfigure}[b]{0.8\textwidth}
		\centering
		\includegraphics[width=\textwidth]{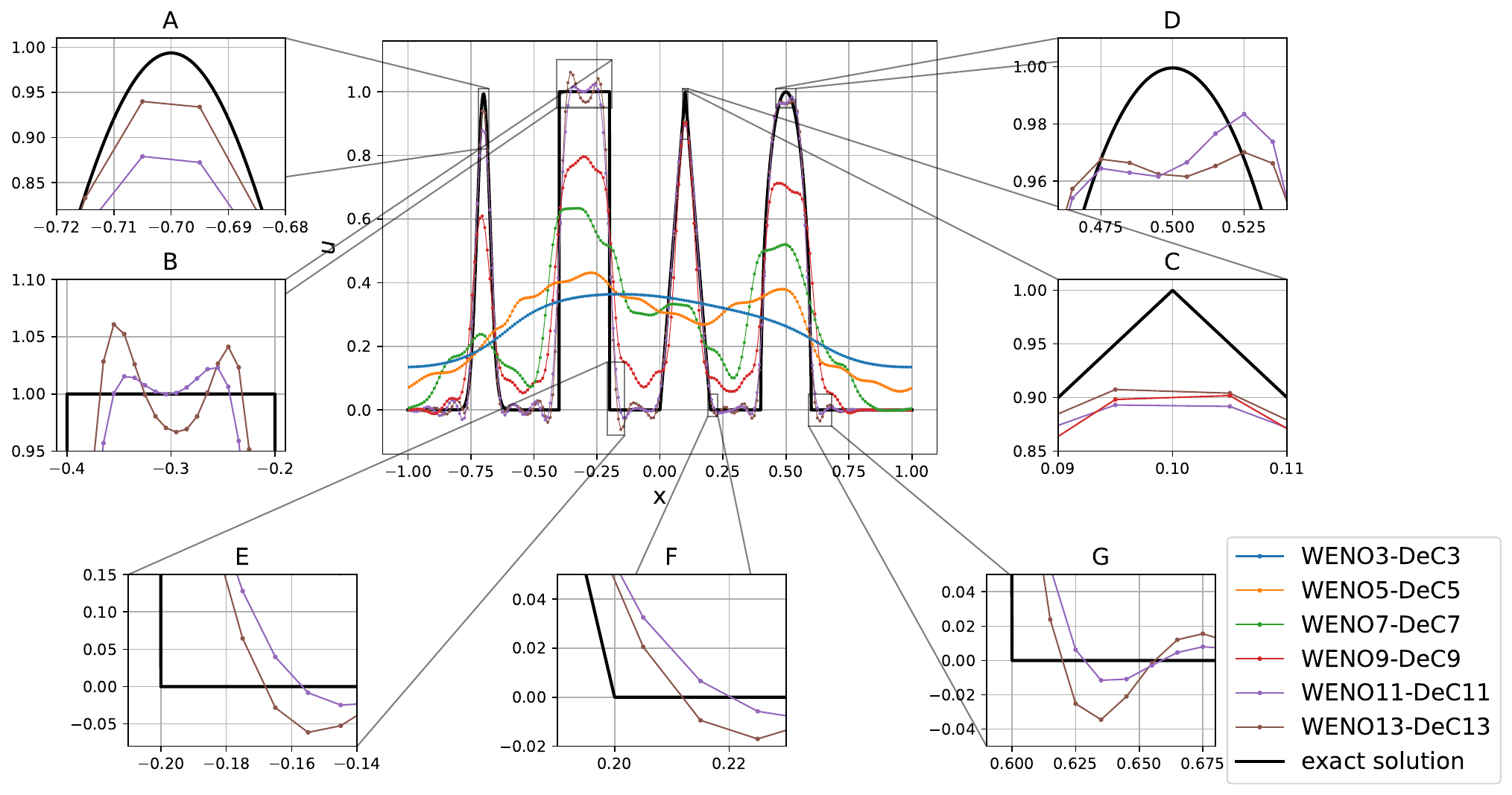}
		\caption{200 elements}
	\end{subfigure}\\
	\begin{subfigure}[b]{0.8\textwidth}
		\centering
		\includegraphics[width=\textwidth]{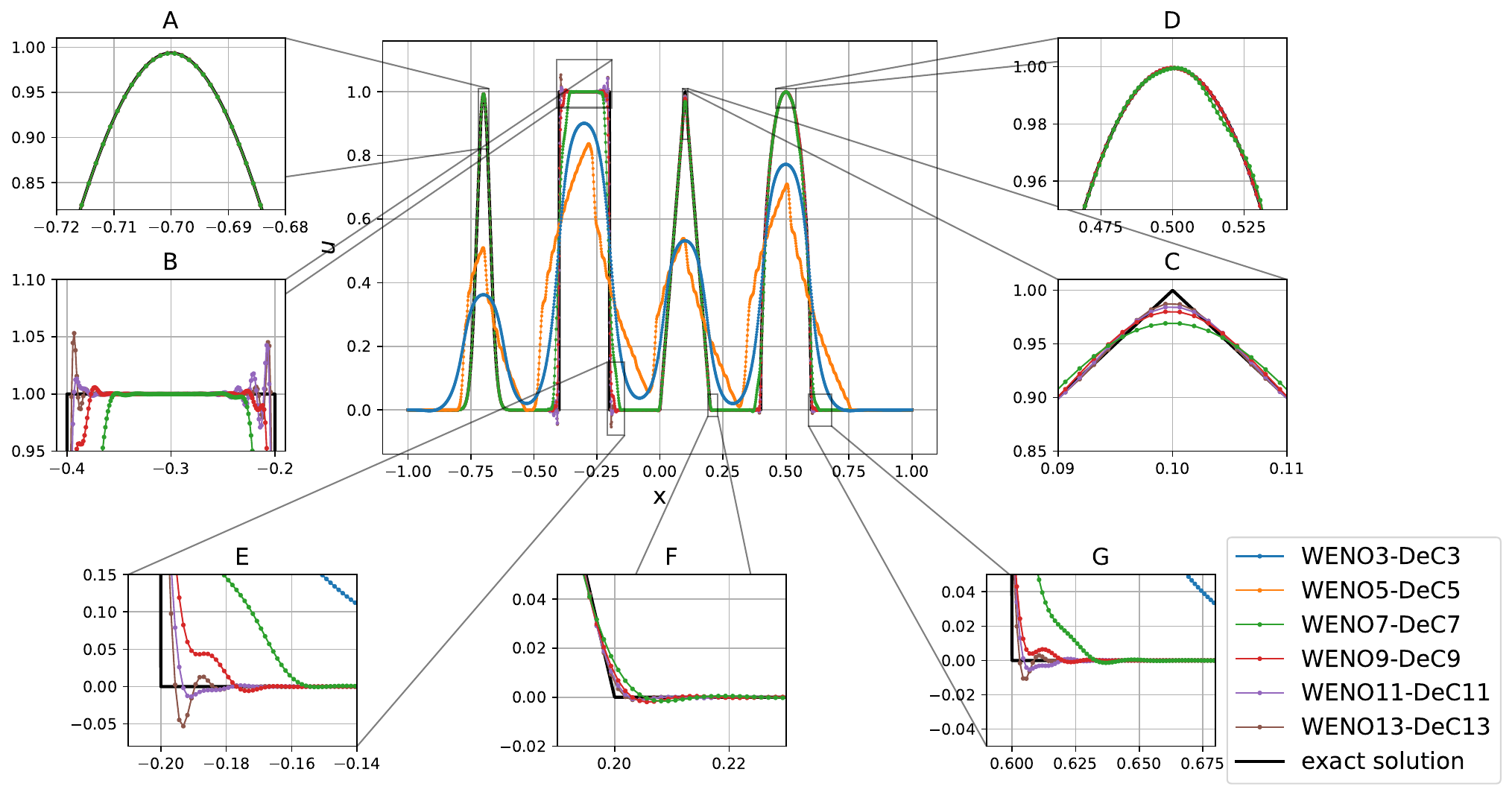}
		\caption{1600 elements}
	\end{subfigure}\\
	\caption{LAE, Test 2: Results obtained with WENO--DeC with final time $T_f:=2000$ and several levels of mesh refinement}
	\label{fig:LAE_1d_composite_wave_longer_time_WENO_DeC}
\end{figure}

Being the solution only $C^{0}$, we do not expect to achieve high order convergence, yet, we would like higher order discretizations to better preserve the solution profile, especially for very long simulation times.
%
%
%
We have considered $T_f:=2000$ as in~\cite{toro2005tvd} and, again, $C_{CFL}:=0.95$.
In Figure~\ref{fig:LAE_1d_composite_wave_longer_time_WENO_DeC}, we report the results obtained with WENO--DeC for 50, 100, 200 and 1600 elements.
We systematically see that, despite being affected by little over- and undershoots, high order methods better preserve the solution profile.
In particular, up to 200 cells, the results obtained with all methods but the ones of order 11 and 13 are affected by a severe smearing. 
For 1600 cells, the results obtained with orders 7 and 9 get comparable with the ones obtained for higher order, while, orders 3 and 5 are still far from convergence.
Let us remark that the coupling of WENO--DeC with extra adaptive limiting strategies would definitely improve the results preventing over- and undershoots. This is planned for future works.

\begin{figure}[htbp]
	\centering
	\begin{subfigure}[b]{0.65\textwidth}
		\centering
		\includegraphics[width=\textwidth]{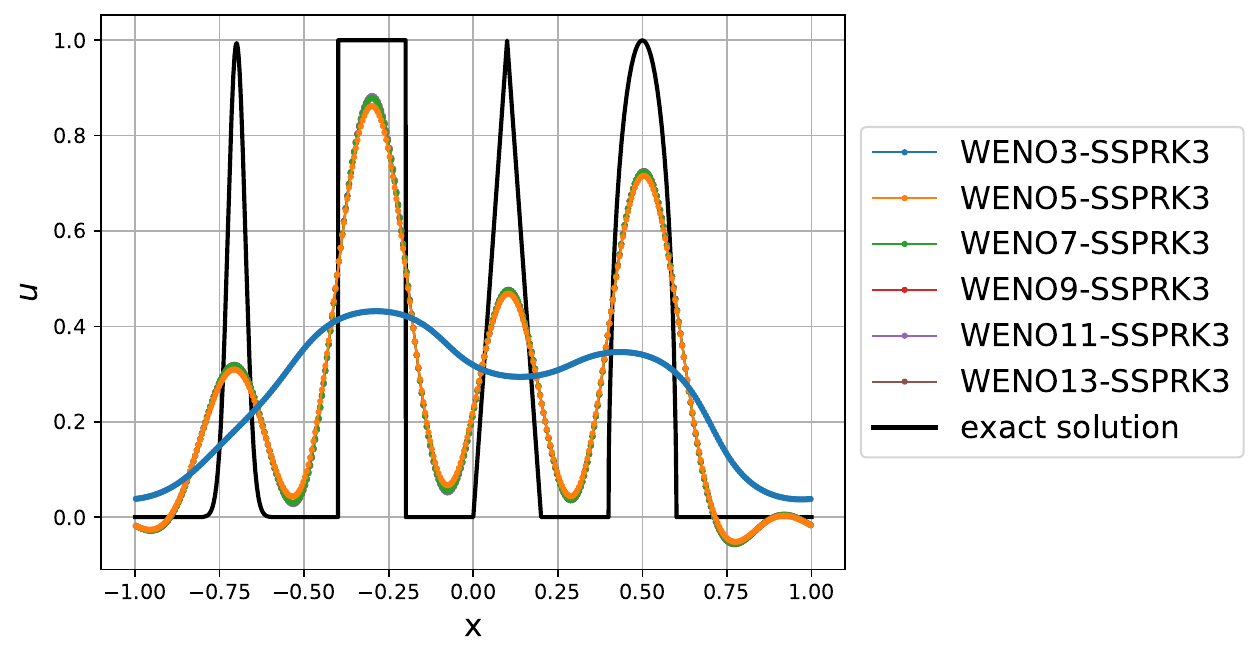}
		\caption{WENO--SSPRK3}
	\end{subfigure}\\
	\quad
	\begin{subfigure}[b]{0.65\textwidth}
		\centering
		\includegraphics[width=\textwidth]{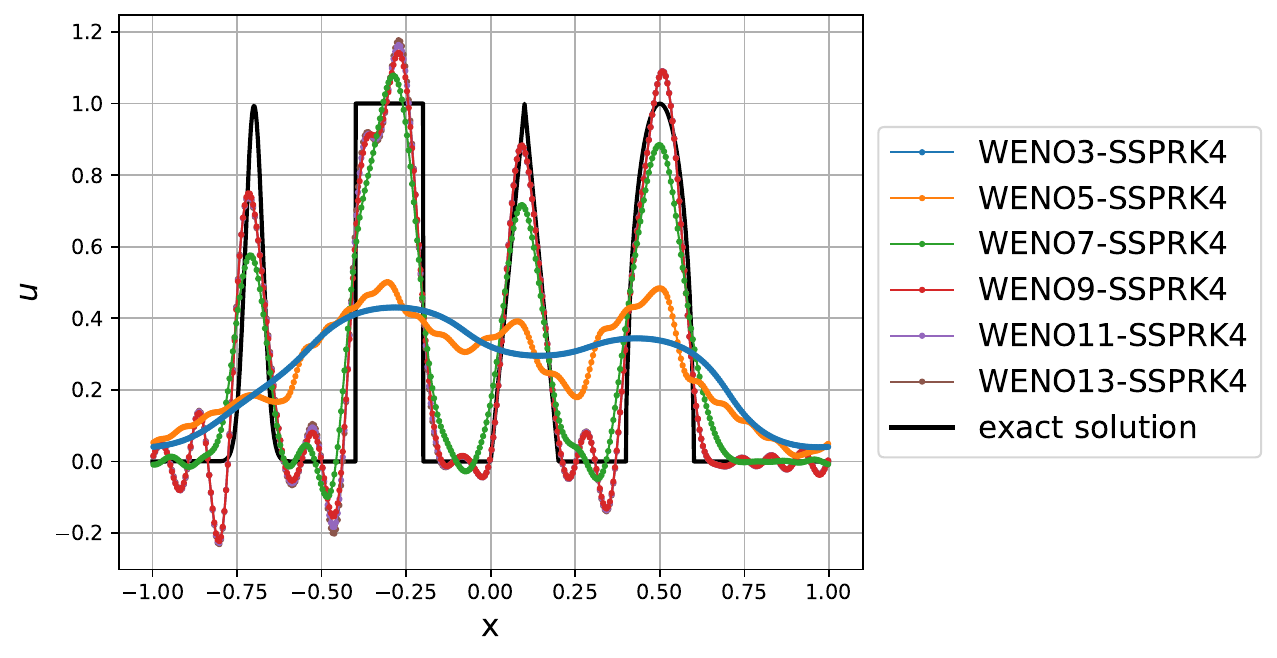}
		\caption{WENO--SSPRK4}
	\end{subfigure}\\
	\begin{subfigure}[b]{0.65\textwidth}
		\centering
		\includegraphics[width=\textwidth]{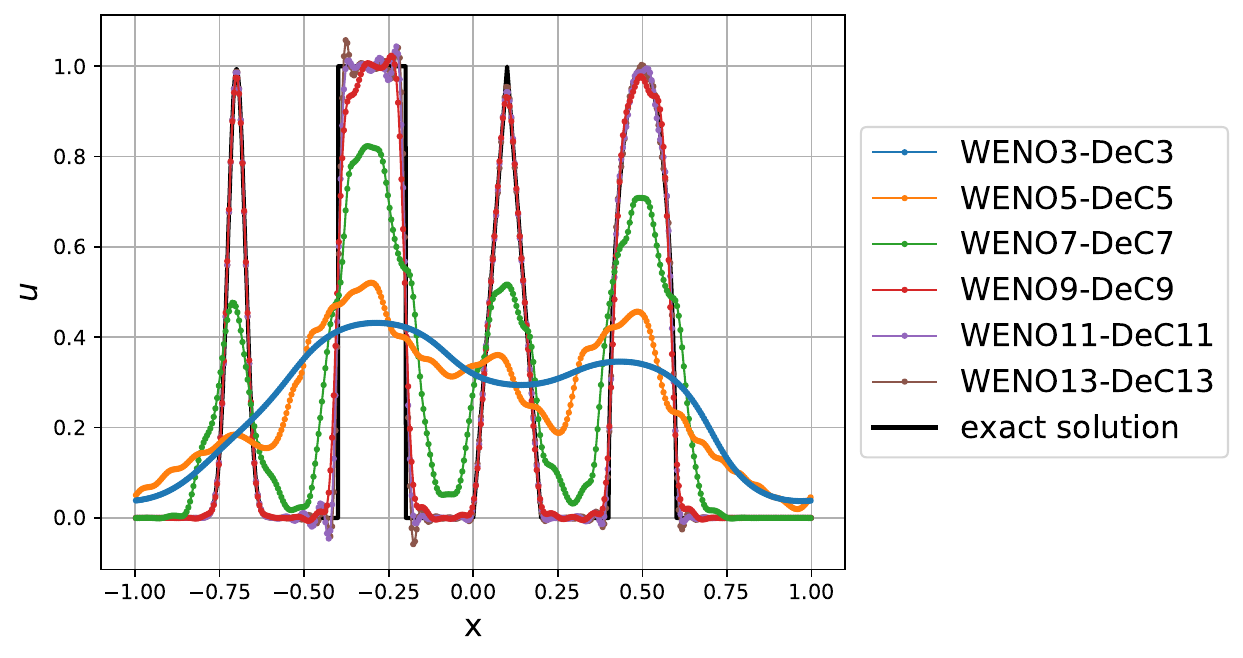}
		\caption{WENO--DeC}
	\end{subfigure}
	\caption{LAE, Test 2: Results obtained with WENO--SSPRK3, WENO--SSPRK4 and WENO--DeC with final time $T_f:=2000$ and $400$ elements}
	\label{fig:LAE_1d_composite_wave_longer_time_WENO_SSPRK}
\end{figure}

In Figure~\ref{fig:LAE_1d_composite_wave_longer_time_WENO_SSPRK}, we compare the results obtained  with WENO--SSPRK3, WENO--SSPRK4 and WENO--DeC on a mesh of 400 elements for same $C_{CFL}$ and final time.
We see a complete smearing of the numerical solution when SSPRK3 is adopted. Moreover, all curves from order 5 on coincide, meaning the the error in time is dominating.
WENO--SSPRK4 seems more able to preserve the shape of the solution, however, the associated results are much less accurate than the ones obtained with WENO--DeC from order 5 on.

Aiming at a quantitative description, we report the errors in the $L^1$-norm obtained with WENO--DeC for several levels of mesh refinement on this problem in Table~\ref{tab:LAE_1d_composite_wave_longer_time_WENO_DeC}. As expected, we see the error decreasing for increasing number of cells and order with no exception.
In Table~\ref{tab:LAE_1d_composite_wave_longer_time_WENO_SSPRK}, we report the errors in the $L^1$-norm obtained with WENO--DeC, WENO--SSPRK3 and WENO--SSPRK4 for 400 cells. It is interesting to notice that, passing from spatial order 7 to 13, there is no substantial decrease in the error for WENO--SSPRK3 and WENO--SSPRK4, meaning that the error in time is dominating. On the other hand, for WENO--DeC, the error always decreases. Moreover, the errors obtained with WENO--DeC from order 7 on are smaller than the ones obtained with WENO--SSPRK3 and WENO--SSPRK4. The difference is one order of magnitude from order 9 on.

\begin{table}[h]
	\centering
	\begin{tabular}{|c|c|c|c|c|c|c|}
		\hline
		Order & 50 cells & 100 cells & 200 cells & 400 cells & 800 cells & 1600 cells \\
		\hline
		3 & 6.726e-01  & 6.681e-01 & 6.346e-01 & 5.942e-01 & 4.871e-01 & 2.404e-01 \  \\
		\hline
		5 & 6.580e-01 & 6.214e-01 & 5.917e-01 & 5.526e-01 & 4.860e-01 & 3.416e-01 \\
		\hline
		7 & 5.974e-01 & 5.821e-01 & 4.624e-01 & 2.941e-01 & 8.227e-02 & 3.438e-02 \\
		\hline
		9 & 4.569e-01 & 4.050e-01 & 2.292e-01 & 4.534e-02 & 2.232e-02 & 1.051e-02 \\
		\hline
		11 & 2.630e-01 & 1.741e-01 & 6.575e-02 & 2.855e-02 & 1.371e-02 & 6.767e-03 \\
		\hline
		13 & 2.186e-01 & 1.263e-01 & 5.728e-02 & 2.392e-02 & 1.138e-02 & 5.427e-03 \\
		\hline
	\end{tabular}
	\caption{LAE, Test 2: Errors in the $L^1$-norm obtained with WENO--DeC with final time $T_f:=2000$ and several levels of mesh refinement}
	\label{tab:LAE_1d_composite_wave_longer_time_WENO_DeC}
\end{table}

\begin{table}[h]
	\centering
	\begin{tabular}{|c|c|c|c|}
		\hline
		Order & WENO--DeC & WENO--SSPRK3 & WENO--SSPRK4 \\
		\hline
		3 & 5.942e-01 & 5.942e-01 & 5.950e-01 \\
		\hline
		5 & 5.526e-01 & 3.250e-01 & 5.491e-01 \\
		\hline
		7 & 2.941e-01 & 3.161e-01 & 1.878e-01 \\
		\hline
		9 & 4.534e-02 & 3.163e-01 & 1.526e-01 \\
		\hline
		11& 2.855e-02 & 3.149e-01 & 1.559e-01 \\
		\hline
		13& 2.392e-02 & 3.149e-01 & 1.577e-01 \\
		\hline
	\end{tabular}
	\caption{LAE, Test 2: Errors in the $L^1$-norm obtained with WENO--DeC, WENO--SSPRK3 and WENO--SSPRK4 with final time $T_f:=2000$ and 400 elements}
	\label{tab:LAE_1d_composite_wave_longer_time_WENO_SSPRK}
\end{table}

\begin{remark}[On smaller final times]
	It is useful to remark that the described accuracy losses of schemes obtained coupling very high order space discretizations with lower order time discretizations may not be visible for small final time.
	Results obtained using WENO--SSPRK3, WENO--SSPRK4 and WENO--DeC with $400$ elements for a final time $T_f:=8$ are comparable, and they are omitted. 
	This may be critical in the context of the many physical applications involving long time simulations, where users could realize only a posteriori to have obtained unsatisfactory results, after huge waste of computational resources.
\end{remark}

The adoption of mSSPRK3 and mSSPRK4 is clearly not a viable solution for long time--simulations, requiring computational times which are too long. Hence, results obtained through such time discretizations are not reported. Let us also remark that they could not be generated at all for order higher than 7 within reasonable time.


\subsection{Euler equations}\label{sec:Euler_1d}
In this section, we perform tests on the one--dimensional Euler equations, obtained by~\eqref{eq:sys} with
\begin{align}
\uvec{u}:=\begin{pmatrix}
	\rho\\
	\rho u\\
	E
\end{pmatrix}, \quad \uvec{f}(\uvec{u}):=\begin{pmatrix}
\rho u\\
\rho u^2 + p\\
(E+p)u
\end{pmatrix},
\end{align}
where $\rho$ is the density of the fluid, $u$ is the velocity, $p$ is the pressure, and $E$ is the total energy.
The system is closed by an equation of state. 
Here, we consider ideal fluids and thus the following equation of state
\begin{equation}
	E=\rho e+\frac{1}{2}\rho u^2,
\end{equation}
with $e:=\frac{p}{\rho(\gamma-1)}$ being the internal energy and $\gamma$, herein set equal to 1.4, adiabatic coefficient given by the ratio between specific heats at constant pressure and volume.

For this system, we consider several different tests. We first investigate a smooth problem in Section~\ref{sec:Euler_1d_sin4} for verifying the order of accuracy.
The remaining tests are, instead, characterized by discontinuous solutions.
More in detail, we consider, five selected Riemann problems from~\cite[Section 12.4]{ToroBook} in Sections \ref{sec:Euler_1d_RP1}, \ref{sec:Euler_1d_RP2}, \ref{sec:Euler_1d_RP3}, \ref{sec:Euler_1d_RP4} and \ref{sec:Euler_1d_RP5}. 
Such tests are extremely challenging, and their numerical resolution is commonly associated with spurious oscillations and/or simulation crashes. 
%
%
The related test informations are reported in Table~\ref{tab:Euler_1d_RP}, where $x_d$ represents the location of the discontinuity in the initial condition.
For all such Riemann problems, the domain is $[0,1]$, and transmissive boundary conditions are adopted at domain extrema.
Let us emphasize that, since we could not run the second Riemann problem with our methods due to simulation crashes caused by the tough character of the problem, we have focused on a relaxed version.

\begin{table}
	\centering
	\begin{tabular}{|c||c|c|c||c|c|c||c||c|}
		\hline
		Test& $\rho_L$ & $u_L$ & $p_L$& $\rho_R$& $u_R$ & $p_R$ & $x_d$ & $T_f$\\\hline\hline
		1   &      1.0 &  0.75 &   1.0&    0.125&   0.0 &  0.1  &  0.3  & 0.2 \\
		\hline
		2   &      1.0 &  -2.0 &   0.4&      1.0&   2.0 &  0.4  &  0.5 & 0.15 \\
		\hline
		 Relaxed 2  &      1.0 &  -1.0 &   0.4&      1.0&   1.0 &  0.4  & 0.5  &  0.15 \\
		\hline
		3   &      1.0 &   0.0 &1000.0&      1.0&   0.0 & 0.01  & 0.5  & 0.012 \\
		\hline
		4   &  5.99924 &19.5975&460.894& 5.99242&-6.19633&46.0950 & 0.4  & 0.035 \\
   		\hline
		5   &     1.0  &-19.59745&1000.0& 1.0 & -19.59745 & 0.01 &  0.8 &  0.012\\
		\hline
	\end{tabular}
	\caption{Test informations for Riemann problems 1, 2 (and relaxed 2), 3, 4, 5. A relaxed version of Riemann problem 2 has been investigated in this work}\label{tab:Euler_1d_RP}
\end{table}

In Section~\ref{sec:Euler_1d_titarev_toro}, we test the methods on a shock--turbulence interaction problem introduced in~\cite{titarev2004finite} as a tough modification of the original version proposed in~\cite{shu1989efficient}.
Finally, in Section~\ref{sec:Euler_1d_critical_aspects}, we point out some critical aspects experienced in the context of these simulations.

As for the LAE, we test methods from order 3 to 13. 
The exact reference solutions to the problems discussed in Sections \ref{sec:Euler_1d_RP1}, \ref{sec:Euler_1d_RP2}, \ref{sec:Euler_1d_RP3}, \ref{sec:Euler_1d_RP4} and \ref{sec:Euler_1d_RP5} have been computed using the library NUMERICA~\cite{toro1999numerica}.
Since for the shock--turbulence interaction problem of Section~\ref{sec:Euler_1d_titarev_toro} no exact solution is available, the reference solution for such a test has been computed through a FV scheme with second order accurate van Leer's minmod spatial discretization~\cite{AbgrallMishranotes} and SSPRK2 time discretization on a very refined mesh of 200,000 cells, with
$C_{CFL}:=0.5$, reconstruction of characteristic variables, and exact Riemann solver numerical flux.
Overall, the results confirm a general advantage in
\begin{itemize}
	\item adoption of high order discretizations with respect to lower order ones;
	\item reconstruction of characteristics variables with respect to conserved ones;
	\item employment of exact Riemann solver with respect to Rusanov, with the difference between the two numerical fluxes getting smaller as the order of accuracy increases.
\end{itemize}
In particular, as also shown in~\cite{qiu2002construction,miyoshi2020short,peng2019adaptive,ghosh2012compact}, reconstruction of characteristic variables results essential in the context of discontinuous solutions to avoid highly oscillatory results and, in some cases, to prevent simulation crashes.
The employment of exact Riemann solver determines, in general, a better capturing of discontinuities and other solution features. On the other hand, the more diffusive character of Rusanov is, sometimes, able to substantially reduce over- and undershoots.
In this context, we confirm what conjectured in~\cite{micalizzi2024impact}: increasing the order of accuracy leads to smaller differences between the results obtained through different numerical fluxes. This is evident for the shock--turbolence interaction problem of Section~\ref{sec:Euler_1d_titarev_toro}.
However, depending on the test and on the mesh refinement, such differences may still be evident even for order 13.
A more systematic discussion on this aspect, considering many other numerical fluxes, is left for future works. Here, we mainly focus on the impact of the order of accuracy.

Also in this case, on smooth problems, comparisons between WENO--DeC, WENO--SSPRK3, WENO--SSPRK4, WENO--mSSPRK3 and WENO--mSSPRK4 confirm order degradation and efficiency loss in the mesh refinement for couplings involving high order space discretizations and lower order time integrations, even for reduced time step. 
Let us remark that the adoption of extra strategies, e.g., a posteriori limiting, would be helpful to reduce spurious oscillations and prevent simulation crashes. This will be matter of investigation of future works, while we only focus on the basic WENO--DeC approach here.

\subsubsection{Advection of smooth density}\label{sec:Euler_1d_sin4}
We consider here, as in~\cite{toro2005tvd}, the advection of a smooth density profile in order to verify the order of accuracy of the investigated methods.
The initial condition in terms of primitive variables reads
\begin{align}
	\begin{cases}
		\rho_0(x)&:=2+\sin^4{\left(\pi x\right)},\\
		u_0(x)&:=u_\infty, \\
		p_0(x)&:=p_\infty,
	\end{cases}
\end{align}
where $u_\infty:=1$ and $p_\infty:=1.$
The computational domain is $\Omega:=[-1,1]$, over which we assume periodic boundary conditions.
The exact solution reads $\uvec{u}(x,t):=\uvec{u}_0(x-u_{\infty}t)$.
We ran the convergence analysis with $C_{CFL}:=0.95$ until the final time $T_f:=2$. 
We only report the errors for the density. Analogous results have been obtained for the other variables and they are omitted to save space.

The results are qualitatively analogous to what obtained in the smooth test considered for the LAE in Section~\ref{sec:LAE_1d_sin4}.
For the sake of compactness, we report the convergence tables of WENO--DeC with reconstruction of characteristic variables and exact Riemann solver only in Table~\ref{tab:Euler_1d_convergence_tables_WENO_DeC_char_exact}.
All other convergence tables of WENO--DeC with reconstruction of characteristic variables and Rusanov numerical flux and with reconstruction of conserved variables are available in the supplementary material, as well as the convergence tables for WENO--SSPRK3, WENO--SSPRK4, WENO--mSSPRK3, and  WENO--mSSPRK4.
Let us remark that the results are qualitatively similar amongst the considered four different settings, coming from the combination of reconstruction of conserved or characteristic variables, and adoption of Rusanov or exact Riemann solver numerical flux.

Looking at Table~\ref{tab:Euler_1d_convergence_tables_WENO_DeC_char_exact}, one can see that the expected order of accuracy has been obtained in all the three considered norms from order 5 on.
As for the LAE, WENO3--DeC3 is characterized by a slight superconvergent character in the mesh refinement.

\begin{table}[htbp]
	\centering
	\caption{Euler equations, Advection of smooth density: convergence tables for WENO--DeC with reconstruction of characteristic variables and exact Riemann solver}
	\label{tab:Euler_1d_convergence_tables_WENO_DeC_char_exact}
	\scalebox{0.65}{ 
		\begin{tabular}{c c c c c c c c}
			\toprule
			\multirow{2}{*}{$N$} & \multicolumn{2}{c}{$L^1$ error $\rho$} & \multicolumn{2}{c}{$L^2$ error $\rho$} & \multicolumn{2}{c}{$L^{\infty}$ error $\rho$} & \multirow{2}{*}{CPU Time} \\
			\cmidrule(lr){2-3} \cmidrule(lr){4-5} \cmidrule(lr){6-7}
			& Error & Order & Error & Order & Error & Order & \\
			\midrule
			
			\multicolumn{8}{c}{\textbf{WENO3--DeC3}} \\ 
			\midrule
			160  &   1.908e-02  &  $-$  &   2.233e-02  &  $-$  &   4.937e-02  &  $-$  &   3.234e-01 \\ 
			320  &   3.706e-03  &  2.364  &   5.755e-03  &  1.957  &   1.641e-02  &  1.590  &   1.191e+00 \\ 
			640  &   4.969e-04  &  2.899  &   9.711e-04  &  2.567  &   3.729e-03  &  2.137  &   4.711e+00 \\ 
			1280  &   3.881e-05  &  3.678  &   7.751e-05  &  3.647  &   3.909e-04  &  3.254  &   1.884e+01 \\ 
			2560  &   2.391e-06  &  4.021  &   3.723e-06  &  4.380  &   1.611e-05  &  4.601  &   7.455e+01 \\ 
			5120  &   1.407e-07  &  4.087  &   1.706e-07  &  4.448  &   5.387e-07  &  4.902  &   2.979e+02 \\ 
			Average order    &      &  3.410  &      &  3.400  &      &  3.297  &    \\ 
			\midrule

			\multicolumn{8}{c}{\textbf{WENO5--DeC5}} \\ 
			\midrule
			80  &   1.586e-03  &  $-$  &   1.371e-03  &  $-$  &   2.391e-03  &  $-$  &   2.742e-01 \\ 
			160  &   5.795e-05  &  4.775  &   5.782e-05  &  4.567  &   1.214e-04  &  4.299  &   9.700e-01 \\ 
			320  &   1.545e-06  &  5.230  &   1.500e-06  &  5.268  &   3.435e-06  &  5.144  &   3.577e+00 \\ 
			640  &   3.636e-08  &  5.409  &   3.144e-08  &  5.576  &   5.727e-08  &  5.906  &   1.417e+01 \\ 
			1280  &   8.323e-10  &  5.449  &   6.821e-10  &  5.527  &   9.233e-10  &  5.955  &   5.653e+01 \\ 
			2560  &   1.972e-11  &  5.399  &   1.727e-11  &  5.304  &   2.579e-11  &  5.162  &   2.259e+02 \\ 
			Average order    &      &  5.252  &      &  5.249  &      &  5.293  &    \\ 
			\midrule

			\multicolumn{8}{c}{\textbf{WENO7--DeC7}} \\ 
			\midrule
			80  &   1.751e-04  &  $-$  &   1.919e-04  &  $-$  &   3.411e-04  &  $-$  &   1.005e+00 \\ 
			160  &   6.198e-07  &  8.142  &   7.444e-07  &  8.010  &   1.736e-06  &  7.618  &   3.485e+00 \\ 
			320  &   4.359e-09  &  7.152  &   8.156e-09  &  6.512  &   4.220e-08  &  5.362  &   1.341e+01 \\ 
			640  &   2.897e-11  &  7.233  &   7.285e-11  &  6.807  &   4.985e-10  &  6.404  &   5.376e+01 \\ 
			1280  &   1.333e-13  &  7.763  &   2.501e-13  &  8.186  &   1.541e-12  &  8.338  &   2.138e+02 \\ 
			Average order    &      &  7.573  &      &  7.379  &      &  6.931  &    \\ 
			\midrule

			\multicolumn{8}{c}{\textbf{WENO9--DeC9}} \\ 
			\midrule
			40  &   1.026e-03  &  $-$  &   1.215e-03  &  $-$  &   2.504e-03  &  $-$  &   7.076e-01 \\ 
			80  &   1.816e-06  &  9.142  &   2.216e-06  &  9.099  &   5.786e-06  &  8.758  &   2.603e+00 \\ 
			160  &   2.123e-09  &  9.740  &   2.394e-09  &  9.854  &   6.370e-09  &  9.827  &   9.466e+00 \\ 
			320  &   4.087e-12  &  9.021  &   4.668e-12  &  9.002  &   1.212e-11  &  9.038  &   3.700e+01 \\ 
			Average order    &      &  9.301  &      &  9.318  &      &  9.207  &    \\ 
			\midrule

			\multicolumn{8}{c}{\textbf{WENO11--DeC11}} \\ 
			\midrule
			40  &   1.179e-04  &  $-$  &   1.260e-04  &  $-$  &   2.437e-04  &  $-$  &   1.584e+00 \\ 
			80  &   1.584e-07  &  9.540  &   2.122e-07  &  9.214  &   5.373e-07  &  8.825  &   5.877e+00 \\ 
			160  &   3.167e-11  &  12.289  &   4.248e-11  &  12.286  &   9.911e-11  &  12.404  &   2.206e+01 \\ 
			Average order    &      &  10.914  &      &  10.750  &      &  10.615  &    \\ 
			\midrule

			\multicolumn{8}{c}{\textbf{WENO13--DeC13}} \\ 
			\midrule
			40  &   1.186e-05  &  $-$  &   1.261e-05  &  $-$  &   2.377e-05  &  $-$  &   3.310e+00 \\ 
			80  &   3.170e-09  &  11.870  &   4.343e-09  &  11.503  &   1.169e-08  &  10.990  &   1.238e+01 \\ 
			160  &   8.733e-14  &  15.148  &   1.000e-13  &  15.406  &   2.580e-13  &  15.467  &   4.663e+01 \\ 
			Average order    &      &  13.509  &      &  13.454  &      &  13.229  &    \\ 
			\midrule

			\bottomrule
	\end{tabular}}
\end{table}

\begin{figure}[htbp]
	\centering
	\begin{subfigure}[t]{0.53\textwidth}
		\centering
		\includegraphics[width=\textwidth]{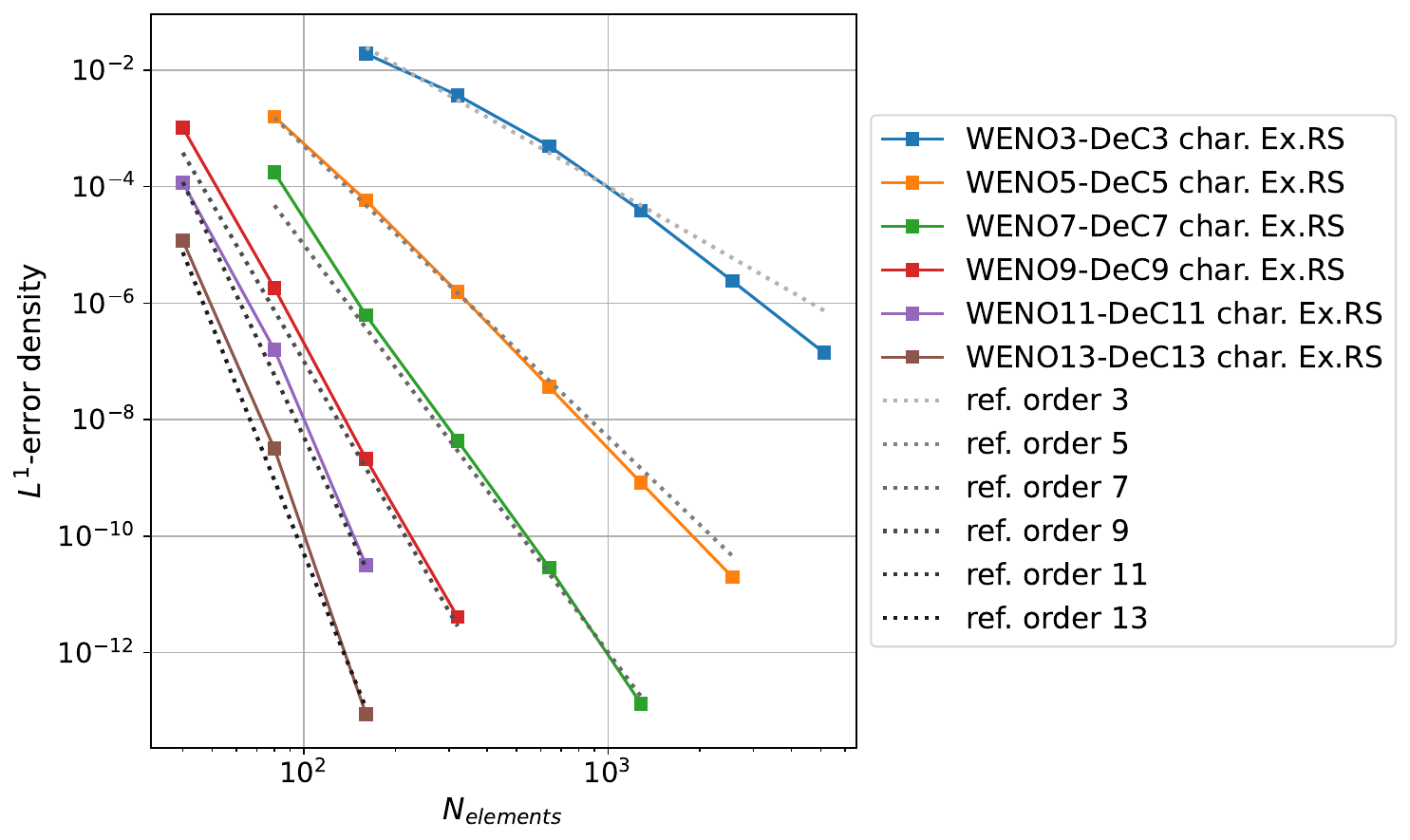}
		\caption{Convergence analysis on the density for WENO--DeC with reconstruction of characteristic variables and exact Riemann solver}
	\end{subfigure}
	\quad
	\begin{subfigure}[t]{0.44\textwidth}
		\centering
		\includegraphics[width=\textwidth]{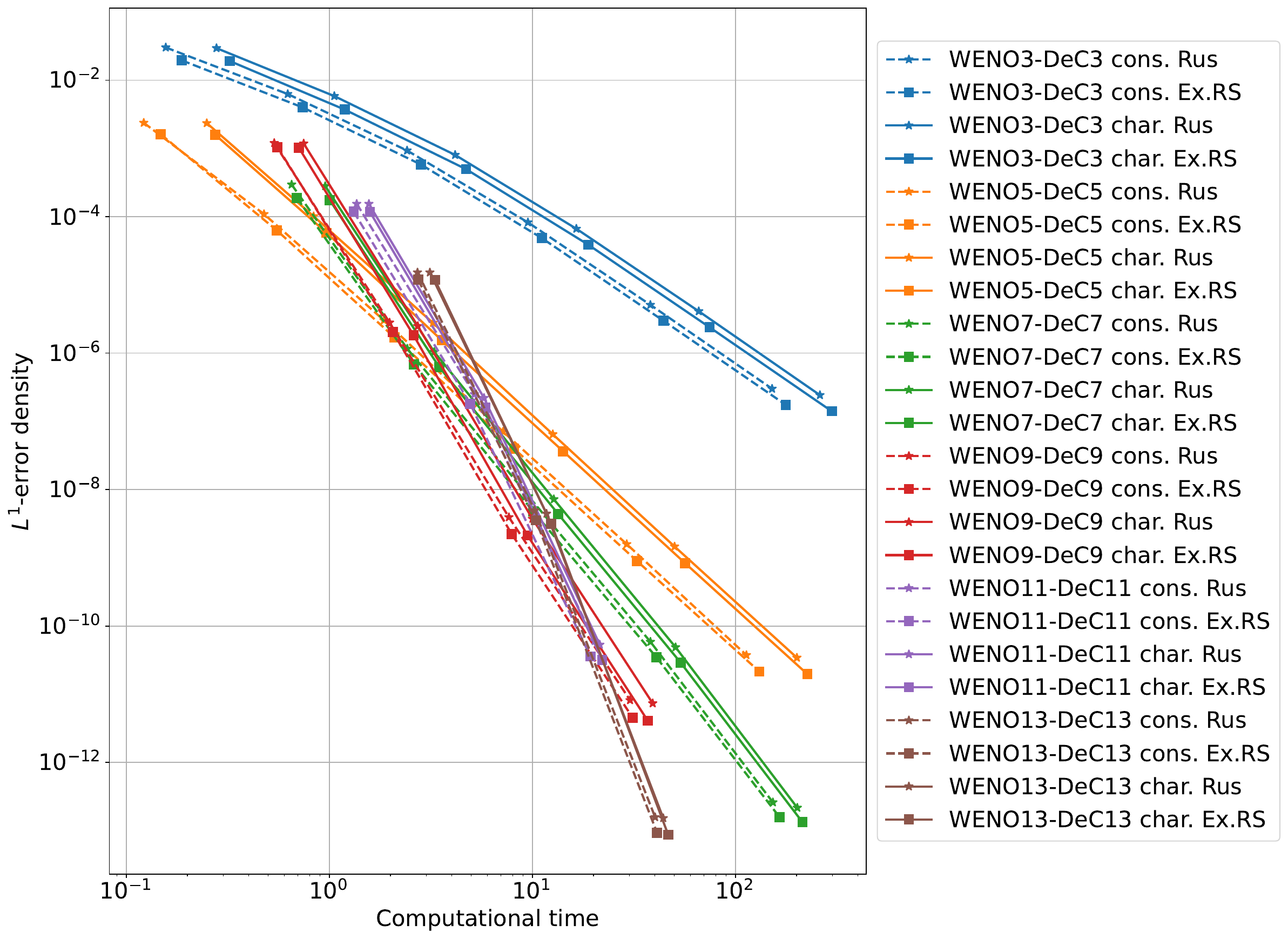}
		\caption{Efficiency analysis for WENO--DeC with all settings}
	\end{subfigure}
	\caption{Euler equations, Advection of smooth density: Results obtained with WENO--DeC}
	\label{fig:Euler_1d_sin4_WENODeC}
\end{figure}

In Figure~\ref{fig:Euler_1d_sin4_WENODeC}, the results of the convergence analysis in the $L^1$-norm obtained with WENO--DeC, reconstruction of characteristic variables and exact Riemann solver are reported, together with the efficiency analysis of all four settings.
We do not report the plots of the convergence analyses for the other three settings, however, the similar trend of the error with respect to the displayed setting can be inferred from the efficiency analysis plot.
Furthermore, we do not display results for the other norms in order to save space, since in fact they are qualitatively analogous.

As already anticipated, the left plot in Figure~\ref{fig:Euler_1d_sin4_WENODeC} confirms the achievement of the expected rates of convergence.
From the right plot, instead, we see the same trend obtained for LAE: high order methods are better for small errors.
Here, however, we have the possibility to make further comparisons concerning the efficiency of different settings corresponding to the adopted choices of reconstructed variables and numerical flux.
For this smooth test, reconstruction of conserved variables is more computationally convenient than reconstruction of characteristic ones; moreover, the exact Riemann solver is slightly more computationally convenient than Rusanov.
We will see in the following that, despite being more computationally expensive (at least in this case), the reconstruction of characteristic variables is essential in the reduction of spurious oscillations in tests involving discontinuities.

\begin{figure}[htbp]
	\centering
	\begin{subfigure}[b]{0.45\textwidth}
		\centering
		\includegraphics[width=\textwidth]{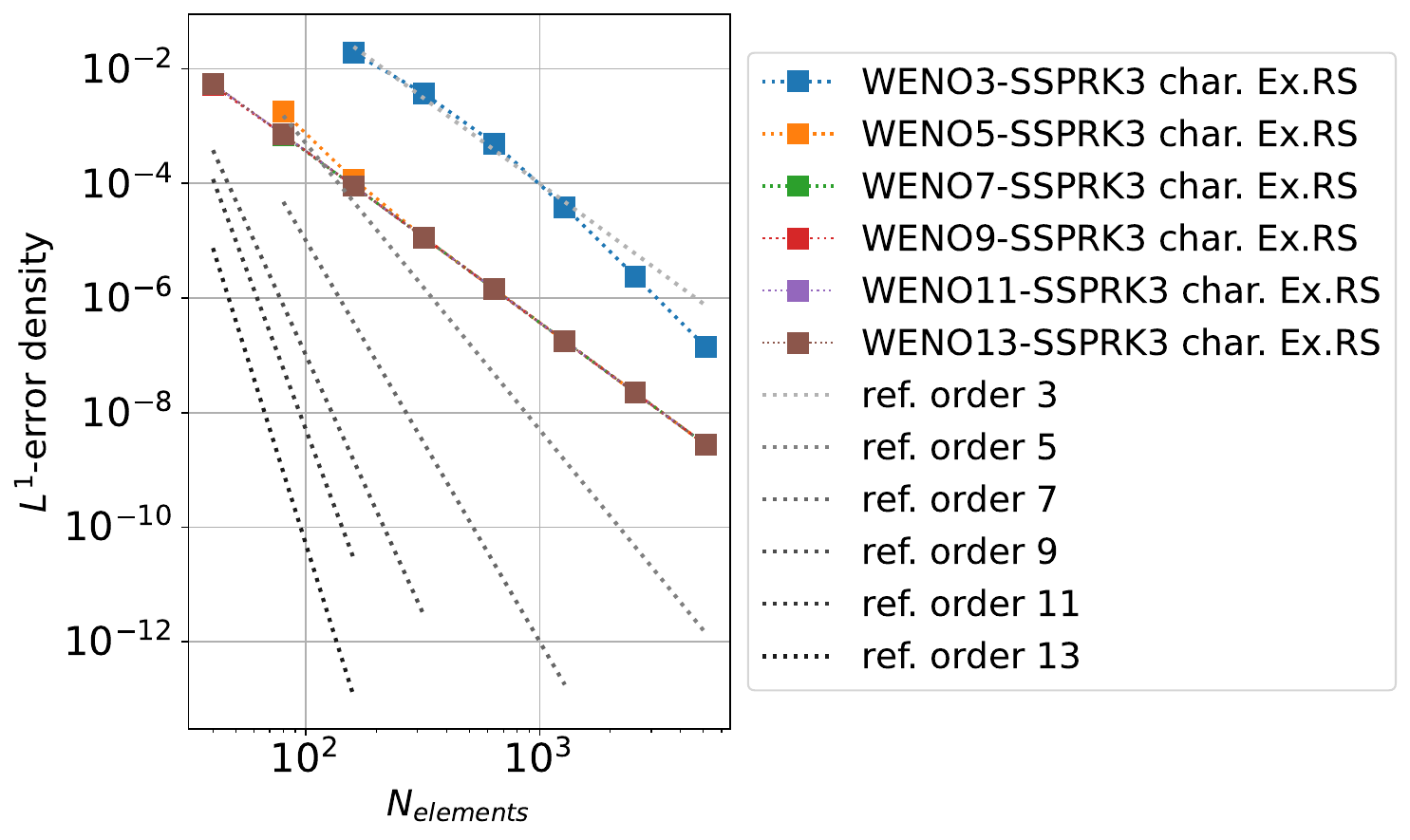}
		\caption{Convergence analysis on the density for WENO--SSPRK3 with reconstruction of characteristic variables and exact Riemann solver}
		\label{fig:Euler_1d_sin4_WENOSSPRK3_convergence}
	\end{subfigure}
	\quad
	\begin{subfigure}[b]{0.45\textwidth}
		\centering
		\includegraphics[width=\textwidth]{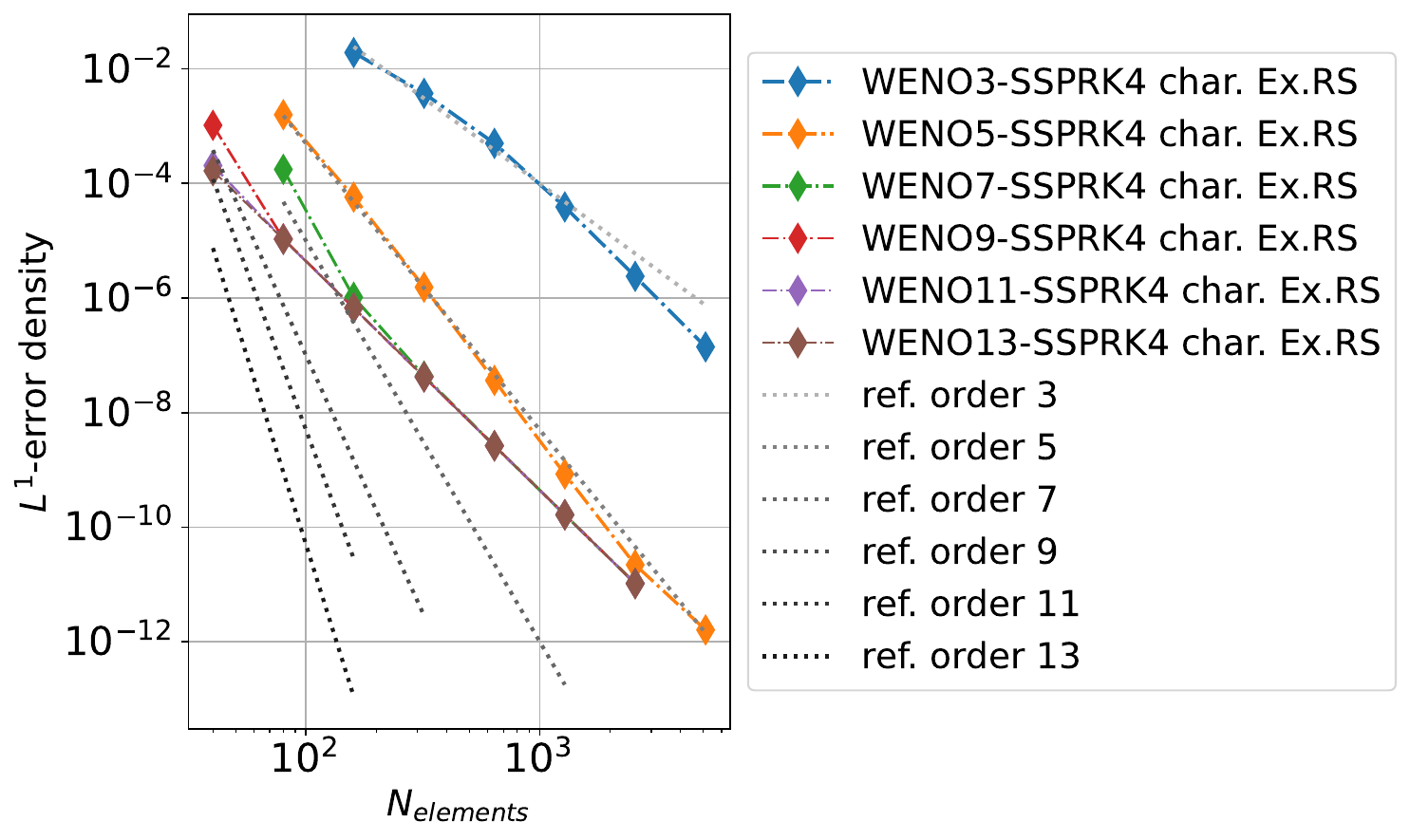}
		\caption{Convergence analysis on the density for WENO--SSPRK4 with reconstruction of characteristic variables and exact Riemann solver}
		\label{fig:Euler_1d_sin4_WENOSSPRK4_convergence}
	\end{subfigure}
	
	\vspace{0.5cm}
	
	\begin{subfigure}[b]{0.47\textwidth}
		\centering
		\includegraphics[width=\textwidth]{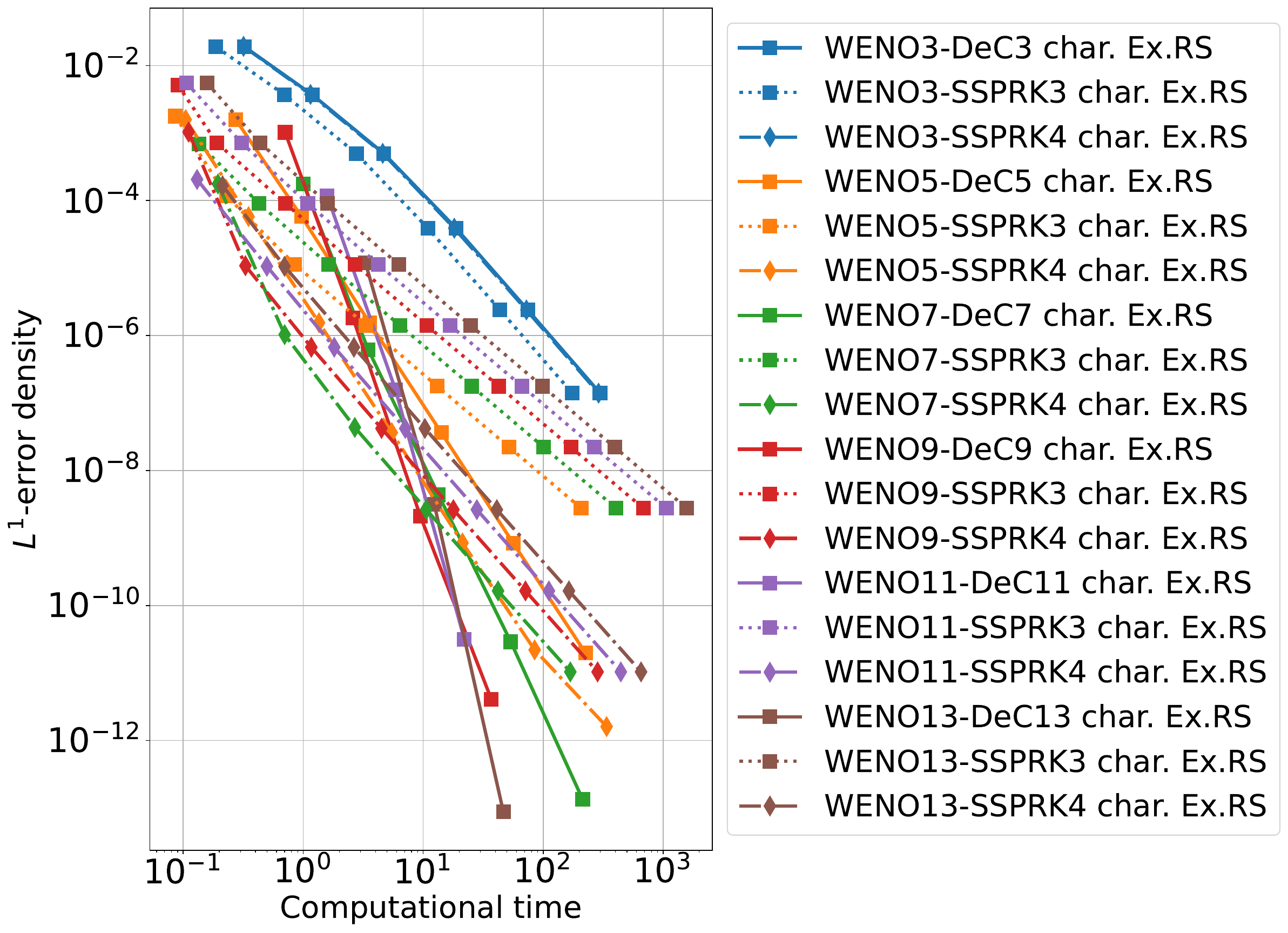}
		\caption{Efficiency analysis for WENO--SSPRK3, WENO--SSPRK4, and WENO--DeC with reconstruction of characteristic variables and exact Riemann solver}
		\label{fig:Euler_1d_sin4_WENOSSPRK_efficiency_char_exact}
	\end{subfigure}
	\quad
	\begin{subfigure}[b]{0.47\textwidth}
		\centering
		\includegraphics[width=\textwidth]{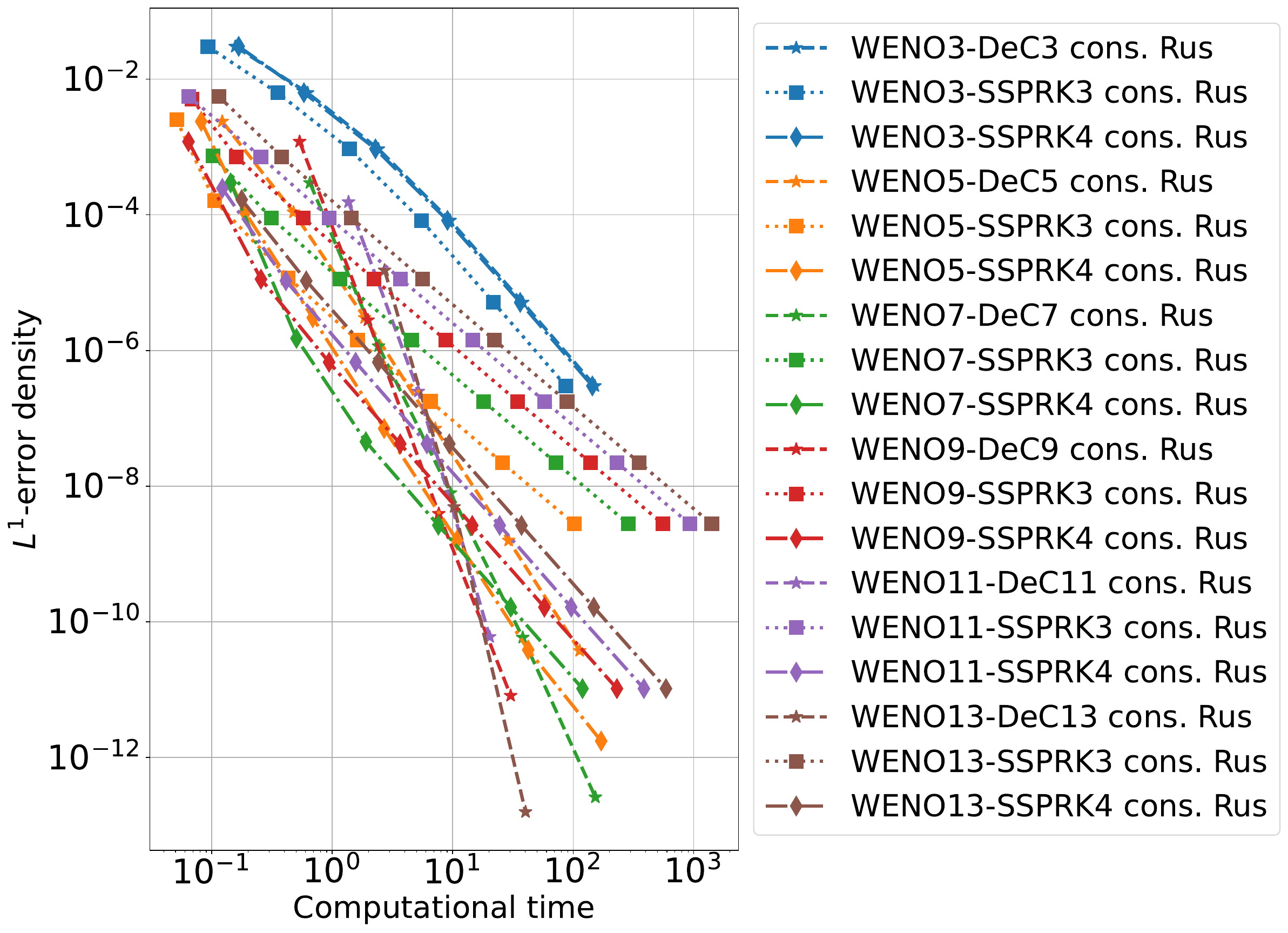}
		\caption{Efficiency analysis for WENO--SSPRK3, WENO--SSPRK4, and WENO--DeC with reconstruction of conserved variables and Rusanov}
		\label{fig:Euler_1d_sin4_WENOSSPRK_efficiency_cons_rusanov}
	\end{subfigure}
	
	\caption{Euler equations, Advection of smooth density: Results obtained with WENO--SSPRK3 and WENO--SSPRK4 on top and efficiency comparison with WENO--DeC on the bottom}
	\label{fig:Euler_1d_sin4_WENOSSPRK}
\end{figure}

As done for the LAE, we compare the results obtained through WENO--DeC with the ones obtained through WENO--SSPRK3 and WENO--SSPRK4. 
The convergence analyses in the $L^1$-norm for WENO--SSPRK3 and WENO--SSPRK4 with reconstruction of characteristic variables and exact Riemann solver are reported in Figures~\ref{fig:Euler_1d_sin4_WENOSSPRK3_convergence} and~\ref{fig:Euler_1d_sin4_WENOSSPRK4_convergence} respectively, along with efficiency comparisons with WENO--DeC for the same setting in~\ref{fig:Euler_1d_sin4_WENOSSPRK_efficiency_char_exact} and for reconstruction of conserved variables and Rusanov in~\ref{fig:Euler_1d_sin4_WENOSSPRK_efficiency_cons_rusanov}.
We clearly see, from Figure~\ref{fig:Euler_1d_sin4_WENOSSPRK3_convergence}, that with SSPRK3 third order of accuracy is obtained for any order of the space discretization from 5 on.
Instead, with SSPRK4, the asymptotic rate of convergence of convergence is clearly 4 for any space discretization from order 7 on, as one can see in Figure~\ref{fig:Euler_1d_sin4_WENOSSPRK4_convergence}.
Also in this context, we observe the usual unexpected superconvergence from WENO3.
WENO5--SSPRK4 achieves order 5 until very high levels of refinement, and only in the last one the convergence rate decreases to approximately four.
We are not sure, whether this is caused by the (expected) order degradation due to the choice of a lower order time discretization or by machine precision effects.
In any case, this is in line with what expressed in Remark~\ref{rmk:asymptotic_convergence}: some levels of refinements are required for the error in space to dominate over the error in time.
This is also the case for WENO7--SSPRK4, which achieves order 7 in the first refinement before becoming fourth order accurate.
The efficiency analysis comparisons reported in Figure~\ref{fig:Euler_1d_sin4_WENOSSPRK_efficiency_char_exact}, for reconstruction of characteristic variables and exact Riemann solver, and in Figure~\ref{fig:Euler_1d_sin4_WENOSSPRK_efficiency_cons_rusanov}, for reconstruction of conserved variables and Rusanov, confirm what observed for the LAE: for small errors, lower order time integration causes efficiency losses in very high order space discretizations.
Similar results have been obtained for reconstruction of conserved variables and exact Riemann solver and reconstruction of characteristic variables and Rusanov.

\begin{figure}[htbp]
	\centering
	\begin{subfigure}[b]{0.55\textwidth}
		\centering
		\includegraphics[width=1.3\linewidth, trim={0 710 0 0}, clip]{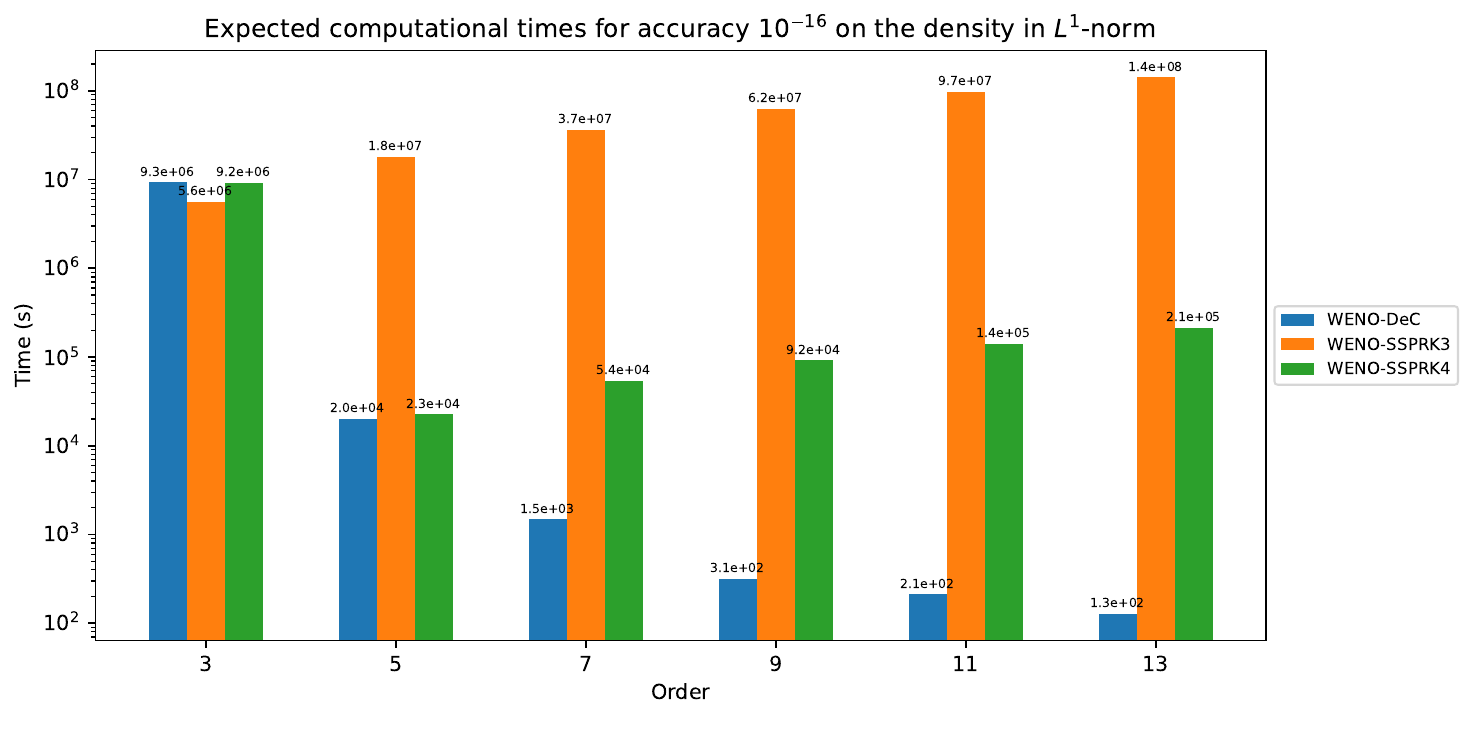}
		\caption{Results for $L^1$-norm}
	\end{subfigure}
	\\
	\begin{subfigure}[b]{0.55\textwidth}
		\centering
		\includegraphics[width=1.3\linewidth, trim={0 360 0 360}, clip]{figures/Euler_1D/sin4/advection_sin4_longer_domain_paper_Euler_1D_char_exact_standard_bar_plots_expected_times.pdf}
		\caption{Results for $L^2$-norm}
	\end{subfigure}
	\\
	\begin{subfigure}[b]{0.55\textwidth}
		\centering
		\includegraphics[width=1.3\linewidth, trim={0 0 0 700}, clip]{figures/Euler_1D/sin4/advection_sin4_longer_domain_paper_Euler_1D_char_exact_standard_bar_plots_expected_times.pdf}
		\caption{Results for $L^{\infty}$-norm}
	\end{subfigure}
	\caption{Euler equations, Advection of smooth density: Expected computational times in seconds to reach an accuracy level equal to $10^{-16}$ on the density with reconstruction of characteristic variables and exact Riemann solver}
	\label{fig:expected_time_Euler_1d}
\end{figure}

Again, we would like to give a concrete meaning to the waste of computational resources obtained adopting SSPRK3 and SSPRK4 in combination with higher order WENO reconstructions.
As a quantitative example, we report in Figure~\ref{fig:expected_time_Euler_1d}, the estimated computational times needed in order to reach an accuracy level equal to $10^{-16}$ on the density in all norms, for reconstruction of characteristic variables and exact Riemann solver.
As for the LAE, these expected values are obtained considering the regression of the last three points of the curves error versus time in logarithmic scale until reaching the desired tolerance.
Results confirm that with DeC, i.e., with order of temporal accuracy matching the one of the spatial discretization, increasing the order always leads to computational advantages.
Instead, for SSPRK3 and SSPRK4, from order 5 on, increasing the order of the space discretization determines an increase in the expected computational time.
As a matter of fact, this means that the adoption of lower order time discretizations completely spoils the advantages coming from high order space discretizations.
%
%
Let us notice that the estimated times in the three different norms are rather similar, which is a further confirmation of the fact that the expected convergence trends have been obtained in all the considered norms.
Let us also remark that results obtained with all other settings are qualitatively analogous, and therefore we omit them in order to save space.

\begin{figure}[htbp]
	\centering
	\begin{subfigure}[b]{0.45\textwidth}
		\centering
		\includegraphics[width=\textwidth]{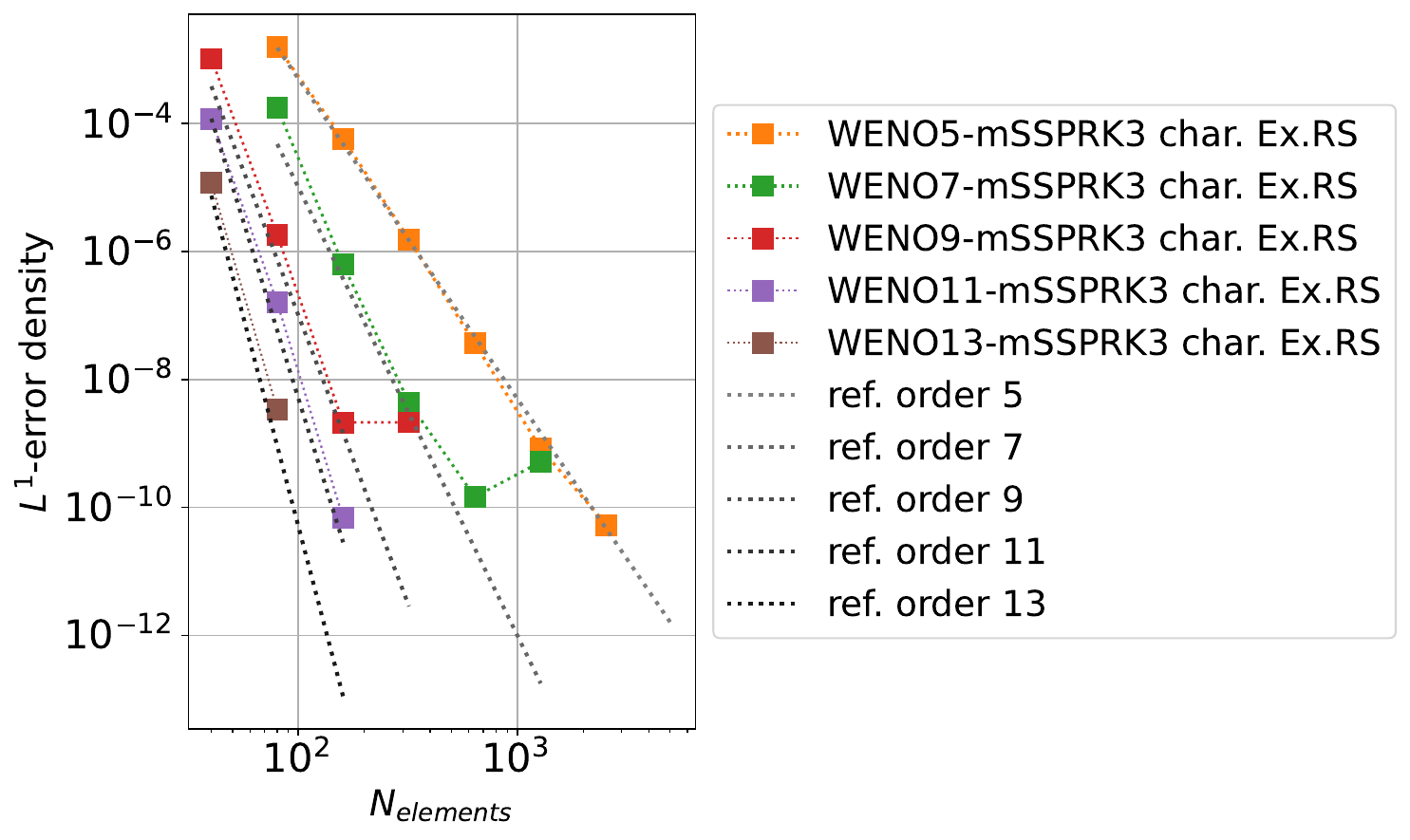}
		\caption{Convergence analysis on the density for WENO--mSSPRK3 with reconstruction of characteristic variables and exact Riemann solver}
		\label{fig:Euler_1d_sin4_WENOmSSPRK3_convergence}
	\end{subfigure}
	\quad
	\begin{subfigure}[b]{0.45\textwidth}
		\centering
		\includegraphics[width=\textwidth]{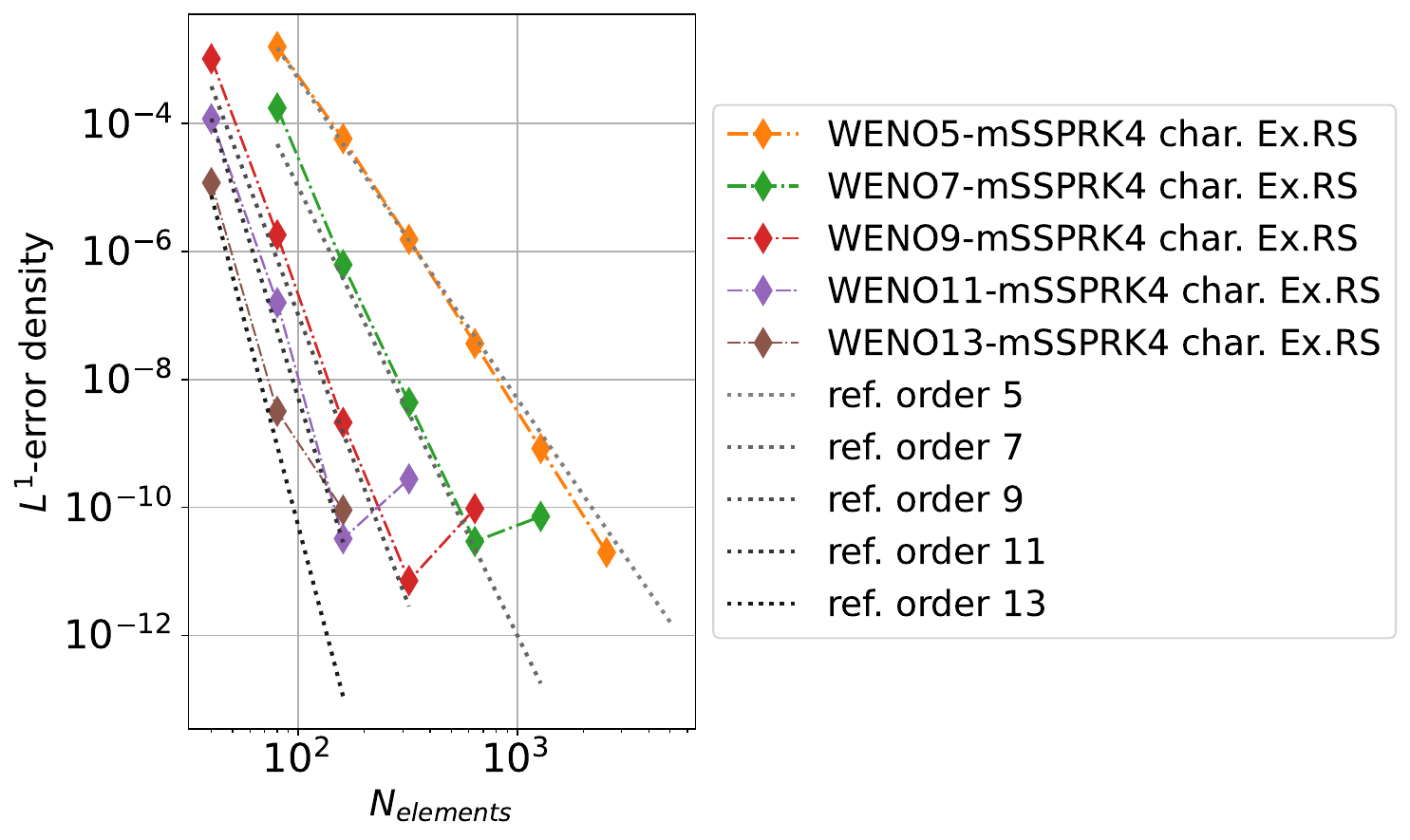}
		\caption{Convergence analysis on the density for WENO--mSSPRK4 with reconstruction of characteristic variables and exact Riemann solver}
		\label{fig:Euler_1d_sin4_WENOmSSPRK4_convergence}
	\end{subfigure}
	
	\vspace{0.5cm}
	
	\begin{subfigure}[b]{0.47\textwidth}
		\centering
		\includegraphics[width=\textwidth]{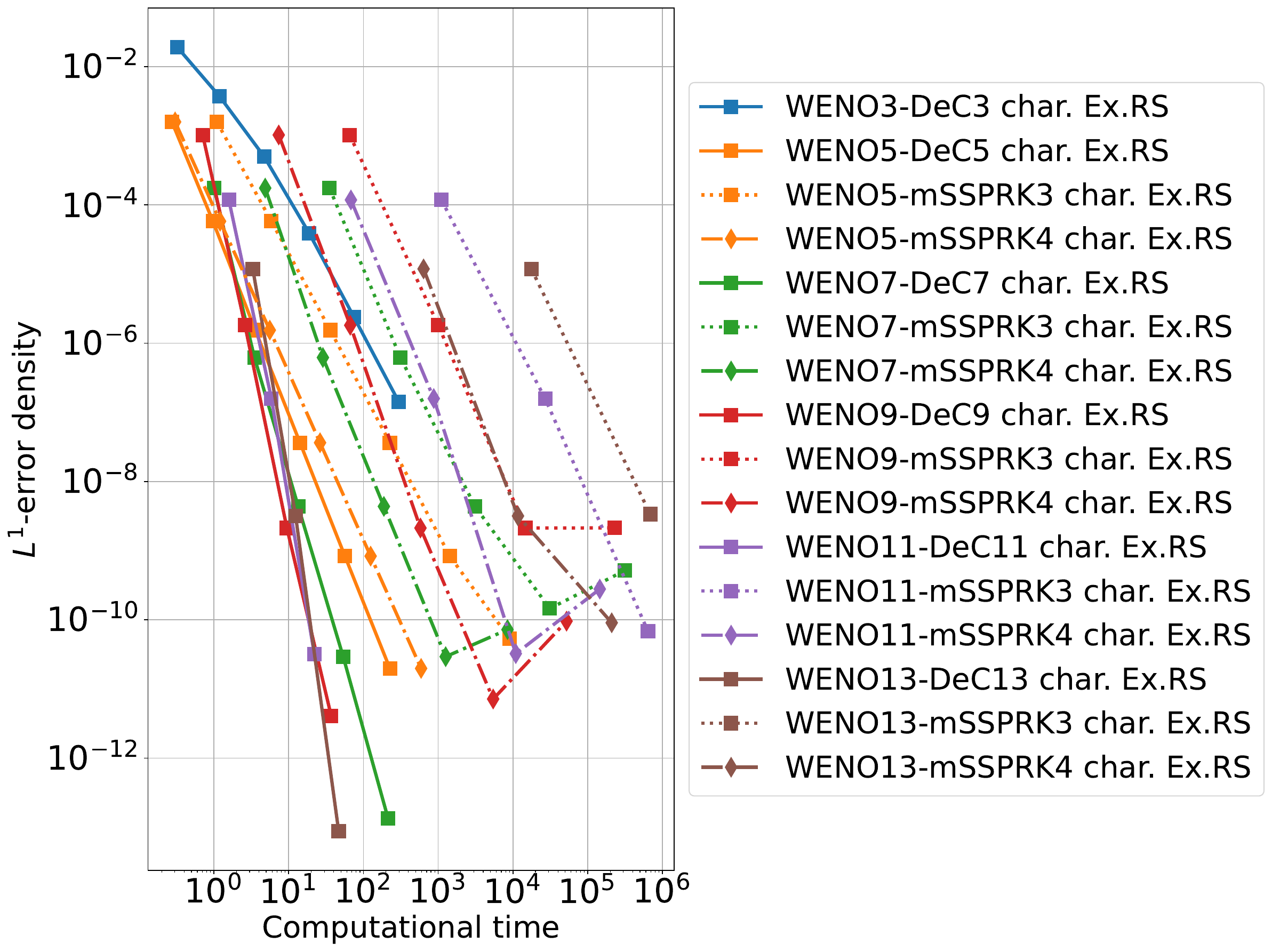}
		\caption{Efficiency analysis for WENO--mSSPRK3, WENO--mSSPRK4, and WENO--DeC with reconstruction of characteristic variables and exact Riemann solver}
		\label{fig:Euler_1d_sin4_WENOmSSPRK_efficiency_char_exact}
	\end{subfigure}
	\quad
	\begin{subfigure}[b]{0.47\textwidth}
		\centering
		\includegraphics[width=\textwidth]{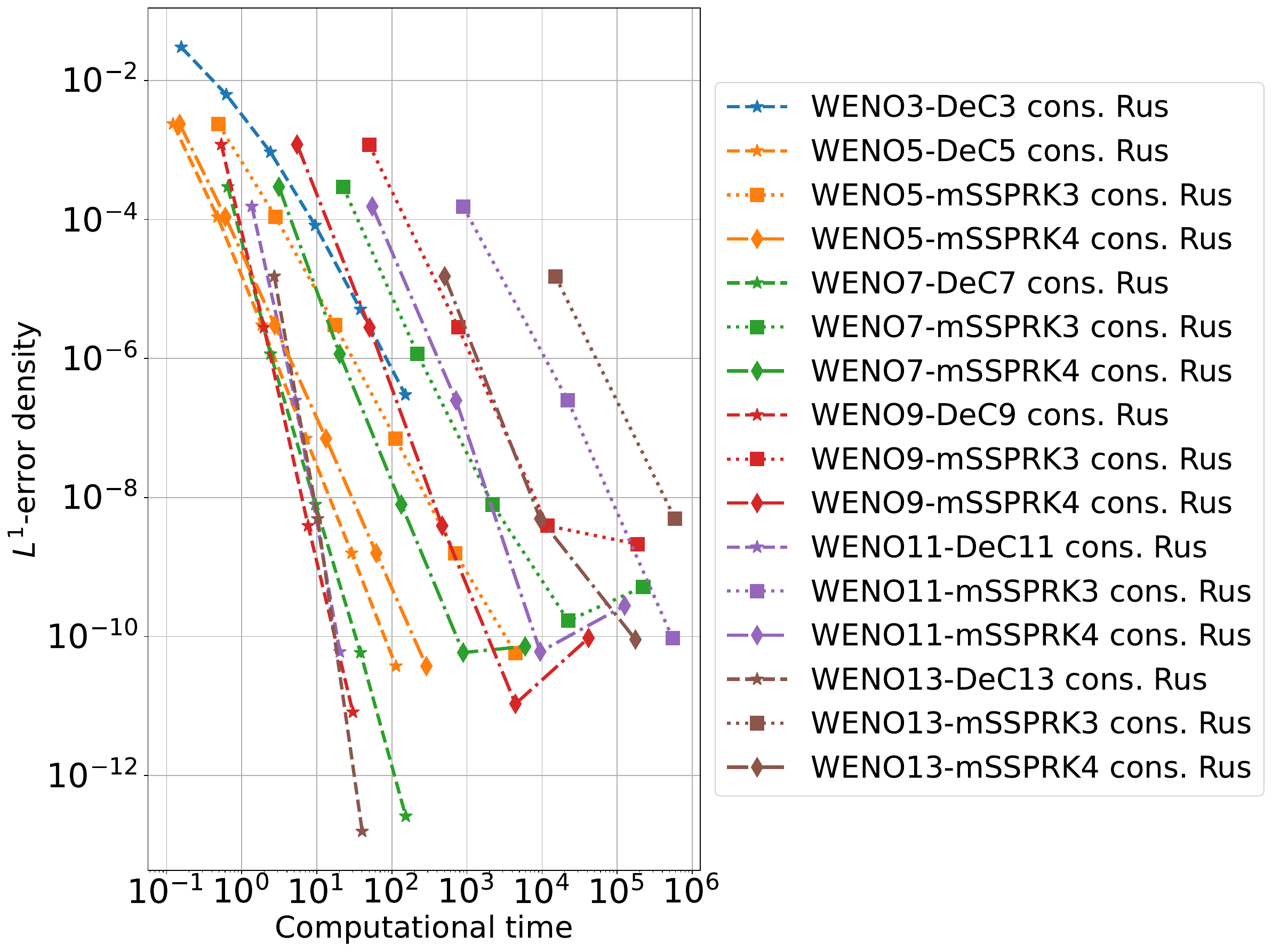}
		\caption{Efficiency analysis for WENO--mSSPRK3, WENO--mSSPRK4, and WENO--DeC with reconstruction of conserved variables and Rusanov}
		\label{fig:Euler_1d_sin4_WENOmSSPRK_efficiency_cons_rusanov}
	\end{subfigure}
	
	\caption{Euler equations, Advection of smooth density: Results obtained with WENO--mSSPRK3 and WENO--mSSPRK4 on top and efficiency comparison with WENO--DeC on the bottom}
	\label{fig:Euler_1d_sin4_WENOmSSPRK}
\end{figure}

Finally, let us compare the results obtained through WENO--DeC with the ones obtained through WENO--mSSPRK3 and WENO--mSSPRK4. 
The convergence analyses plots of WENO--mSSPRK3 and WENO--mSSPRK4 for the $L^1$-norm obtained with reconstruction of characteristic variables and exact Riemann solver are displayed in Figures~\ref{fig:Euler_1d_sin4_WENOmSSPRK3_convergence} and~\ref{fig:Euler_1d_sin4_WENOmSSPRK4_convergence} respectively, while, the efficiency comparisons with WENO--DeC are reported in~\ref{fig:Euler_1d_sin4_WENOmSSPRK_efficiency_char_exact} for the same setting and in~\ref{fig:Euler_1d_sin4_WENOmSSPRK_efficiency_cons_rusanov} for reconstruction of conserved variables and Rusanov.
Just like for LAE, we can appreciate that machine precision effects tend to arise far above usual values, as a result of the strong time step reduction. 
Moreover, such schemes are extremely more computationally costly than WENO--DeC, as can be seen from Figures~\ref{fig:Euler_1d_sin4_WENOmSSPRK_efficiency_char_exact} and~\ref{fig:Euler_1d_sin4_WENOmSSPRK_efficiency_cons_rusanov}.

Also in this context, we conclude that WENO--DeC is the most desirable option amongst the schemes considered.

\subsubsection{Riemann problem 1 (Modified Sod)}\label{sec:Euler_1d_RP1}
This problem is a variant of the Sod shock tube~\cite{sod1978survey}, characterized by nonzero initial velocity of the left state. The solution presents a sonic point, causing troubles in many numerical approaches usually due to the adoption of linearized Riemann solvers~\cite{ToroBook,tang2005sonic}.

\begin{figure}[htbp]
	\centering
	\begin{subfigure}[b]{0.7\textwidth}
		\centering
		\includegraphics[width=\textwidth]{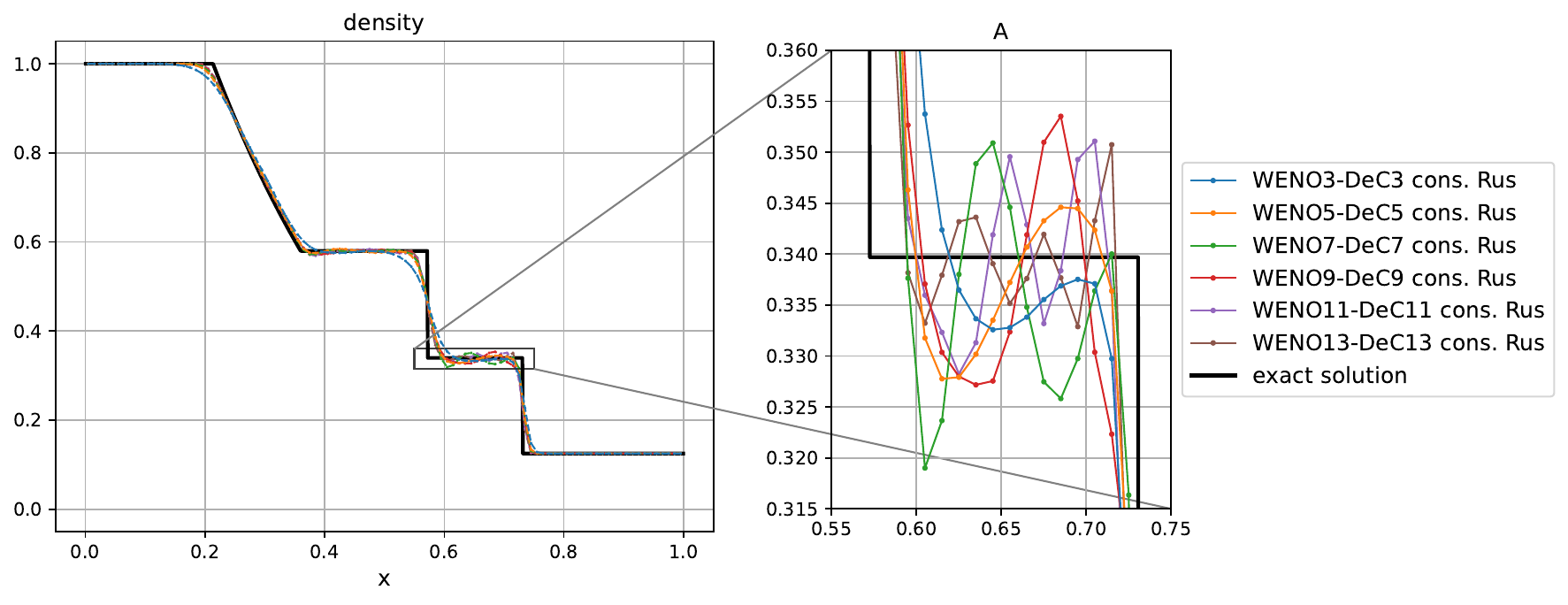}
		\caption{Reconstruction of conserved variables and Rusanov}
	\end{subfigure}
	\\
	\begin{subfigure}[b]{0.7\textwidth}
		\centering
		\includegraphics[width=\textwidth]{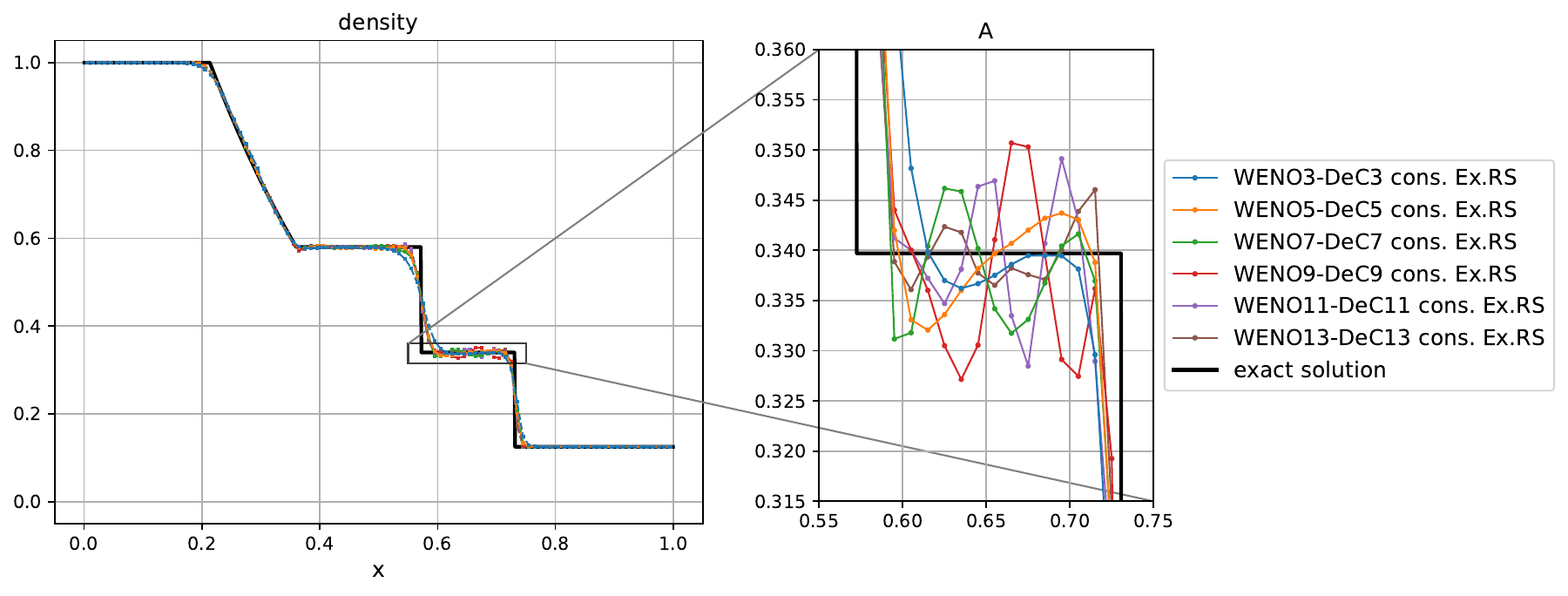}
		\caption{Reconstruction of conserved variables and exact Riemann solver}
	\end{subfigure}\\
	\begin{subfigure}[b]{0.7\textwidth}
		\centering
		\includegraphics[width=\textwidth]{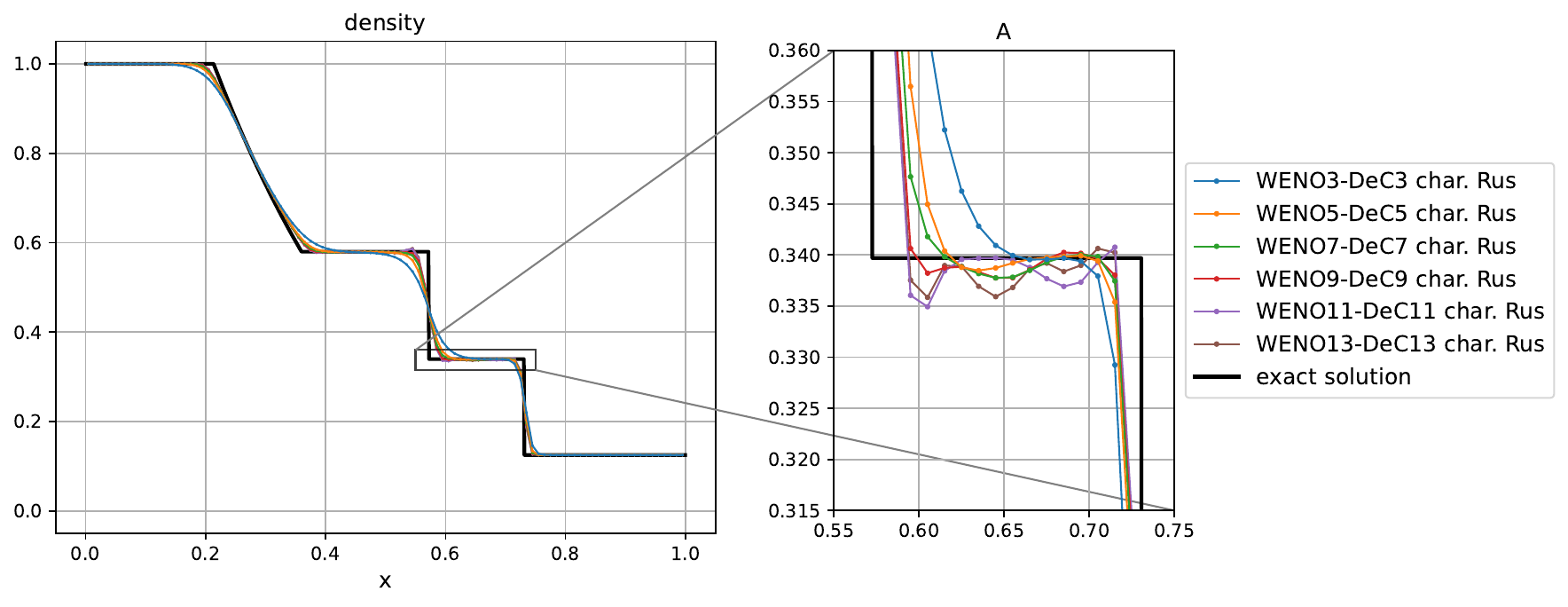}
		\caption{Reconstruction of characteristic variables and Rusanov}
	\end{subfigure}
	\\
	\begin{subfigure}[b]{0.7\textwidth}
		\centering
		\includegraphics[width=\textwidth]{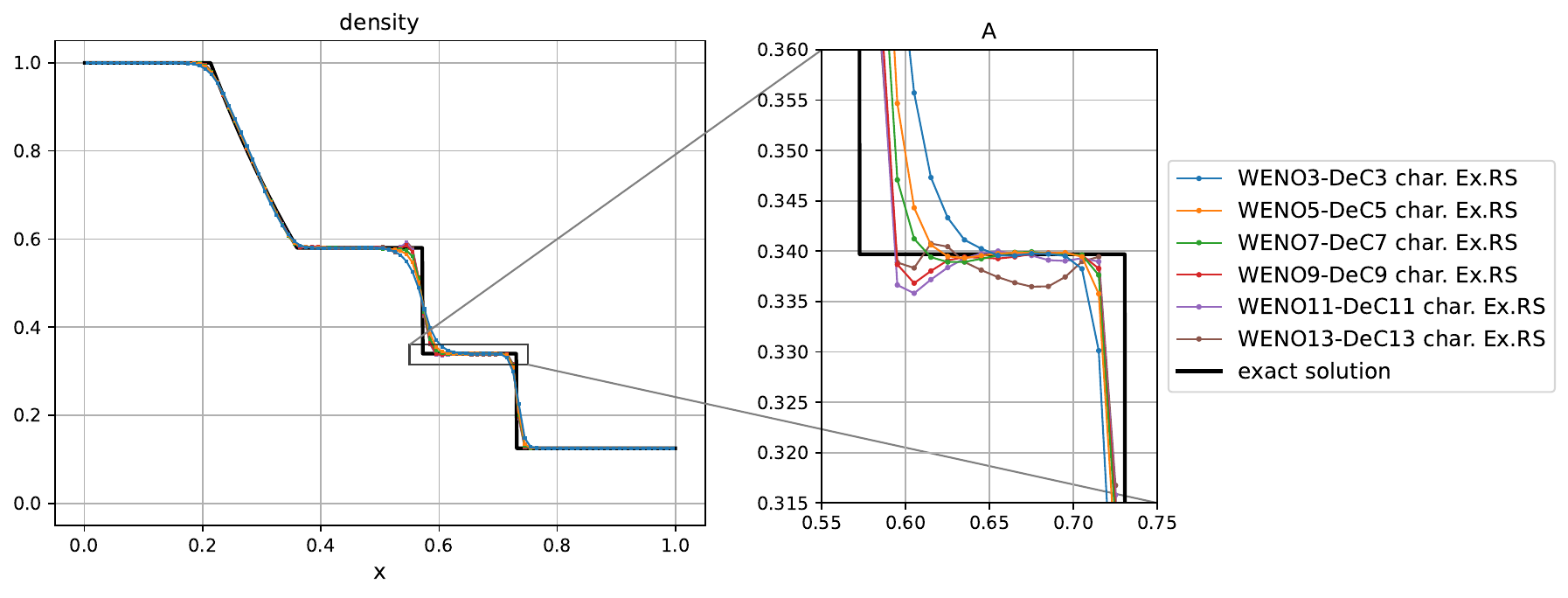}
		\caption{Reconstruction of characteristic variables and exact Riemann solver}
	\end{subfigure}
	\caption{Euler equations, Riemann problem 1: Results for the density obtained through WENO--DeC for reconstruction of conserved and characteristic variables with $C_{CFL}:=0.95$}
	\label{fig:Euler_1d_RP1_all_settings_zoom_density_plateau}
\end{figure}

We did not experience any problem in running this test with all settings and orders for $C_{CFL}:=0.95$ with 100 elements.
In Figure~\ref{fig:Euler_1d_RP1_all_settings_zoom_density_plateau}, we report the results for the density obtained with all settings and orders.
On this test, all schemes produce rather nice results, with sonic point handled perfectly.
However, we can see that the reconstruction of conserved variables leads to spurious oscillations for very high orders, see the zoom on the plateau of the density between the shock wave and the contact discontinuity.
Results are much cleaner when reconstruction of characteristic variables is adopted.
This is in line with what shown in~\cite{qiu2002construction,miyoshi2020short,peng2019adaptive,ghosh2012compact}.
In order to save space, in the following, we will not provide results with reconstruction of conserved variables in the main manuscript.
However, they can be found in the supplementary material.


\begin{figure}[htbp]
	\centering
	\begin{subfigure}[b]{0.9\textwidth}
		\centering
		\includegraphics[width=\textwidth]{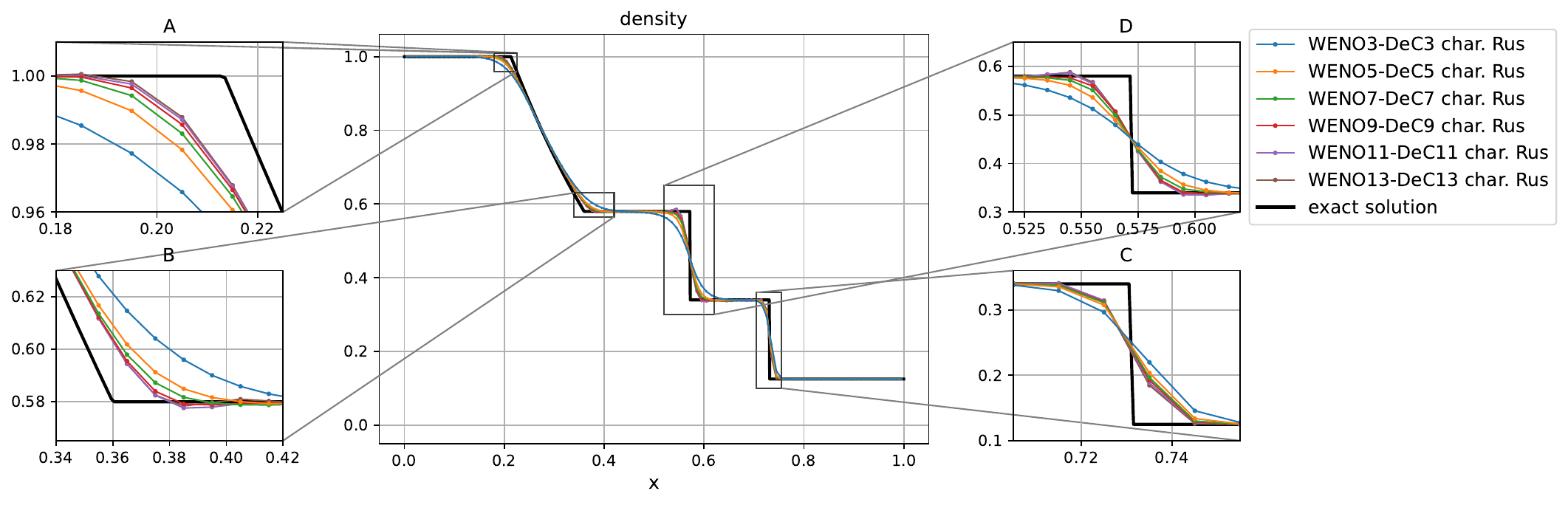}
		\caption{Reconstruction of characteristic variables and Rusanov}
	\end{subfigure}
	\\
	\begin{subfigure}[b]{0.9\textwidth}
		\centering
		\includegraphics[width=\textwidth]{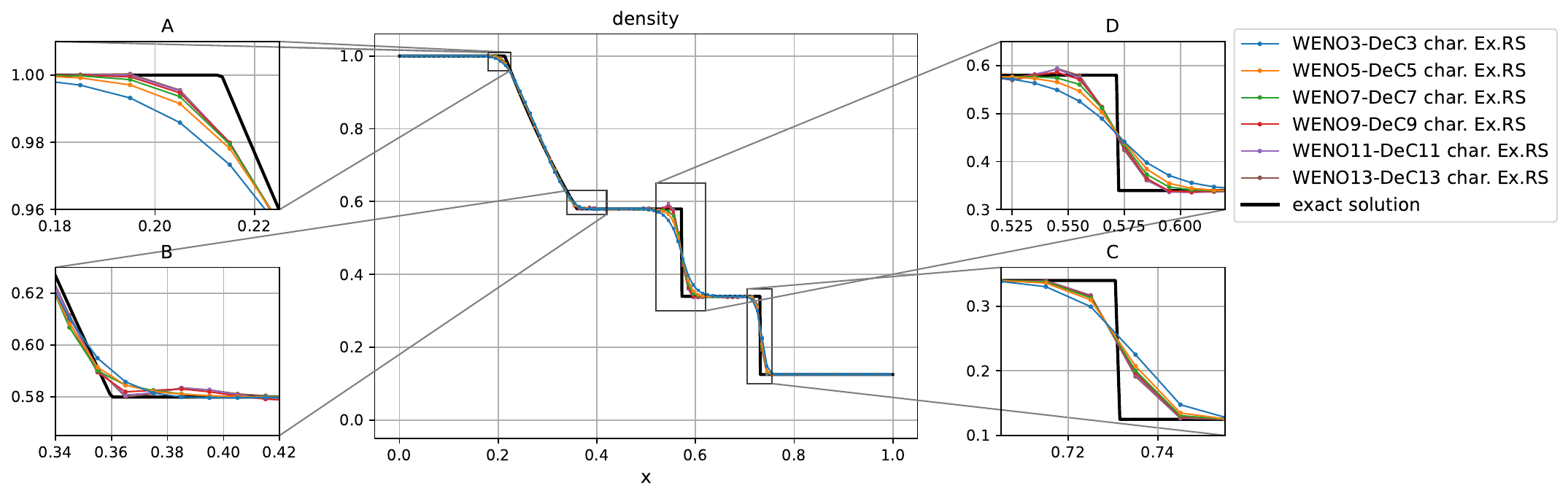}
		\caption{Reconstruction of characteristic variables and exact Riemann solver}
	\end{subfigure}
	\caption{Euler equations, Riemann problem 1: Results for the density obtained through WENO--DeC for reconstruction of characteristic variables with $C_{CFL}:=0.95$}
	\label{fig:Euler_1d_RP1_zoom_density}
\end{figure}

In Figure~\ref{fig:Euler_1d_RP1_zoom_density}, we report some zooms on the density profiles obtained with reconstruction of characteristic variables and both numerical fluxes.
The employment of the exact Riemann solver provides a better capturing of the rarefaction wave with respect to Rusanov, as can be seen in panels A and B.
For very high orders, little over- and undershoots can be observed in panel D and, only for Rusanov numerical flux, also in panel B.
Zoom on other variables are omitted to save space. However, same considerations hold for them: the reconstruction of conserved variables is associated with spurious oscillations for very high orders, with the results obtained for reconstruction of characteristic variables being much cleaner; the adoption of the exact Riemann solver provides a better capturing of some solution details with respect to Rusanov.
Overall, the best results are obtained for reconstruction of characteristic variables and exact Riemann solver.
Let us notice that higher order methods produce more reliable results with respect to lower order ones, and tend to better capture the details of the solution, including discontinuities.

The analysis of the results obtained for the original Sod shock tube test lead to similar conclusions, hence, we do not report them for the sake of compactness.

\subsubsection{Riemann problem 2 (Double rarefaction)}\label{sec:Euler_1d_RP2}
This test, also known as ``123 problem'', is one of the toughest Riemann problems proposed in~\cite{ToroBook}. Two diverging streams determine the formation of a near--vacuum--area.
Simulation crashes due to negative density or pressure are very frequent for this test.
Just like for the previous Riemann problem, we consider 100 cells.


For $0.05 \leq C_{CFL}\leq 0.85$ (we did not test for smaller values), the test works only for order 3 and only if characteristic variables are reconstructed and Rusanov numerical flux is employed.
All other orders and settings crash.
This is why we have focused on a relaxed version of the problem with smaller values of the speed of the diverging streams, according to Table~\ref{tab:Euler_1d_RP}.
Let us remark that, despite this little modification, the test still preserves the challenging features of the original version.
The following comments concern the results obtained on such a new relaxed version with 100 cells.


For $C_{CFL}:=0.95$, only order 3 does not crash for any considered reconstructed variable and numerical flux.
The value of $C_{CFL}$ must be reduced to prevent simulation crashes.
With $C_{CFL}\leq 0.7$, there are no simulation failures for all orders and both numerical fluxes with reconstruction of characteristic variables.
Instead, if conserved variables are reconstructed, lower values of $C_{CFL}$ are required.
More in detail, with exact Riemann solver, for $0.55 \leq C_{CFL}\leq 0.95$, only order 3 works; for $0.05\leq C_{CFL}\leq 0.5$ only orders 3 and 5 work. 
On the other hand, with Rusanov, for $0.8\leq C_{CFL}\leq 0.95$ only order 3 works; for $C_{CFL}\leq 0.75$ orders 3, 5 and 7 work; for $C_{CFL}\leq 0.6$ orders 3,5,7 and 9 work; for $C_{CFL}\leq 0.55$ all orders work besides order 13; for $C_{CFL}\leq 0.45$ all orders work.
We perform the comparison between the settings for $C_{CFL}:= 0.45$.

\begin{figure}[htbp]
	\centering
	\begin{subfigure}[b]{0.9\textwidth}
		\centering
		\includegraphics[width=\textwidth]{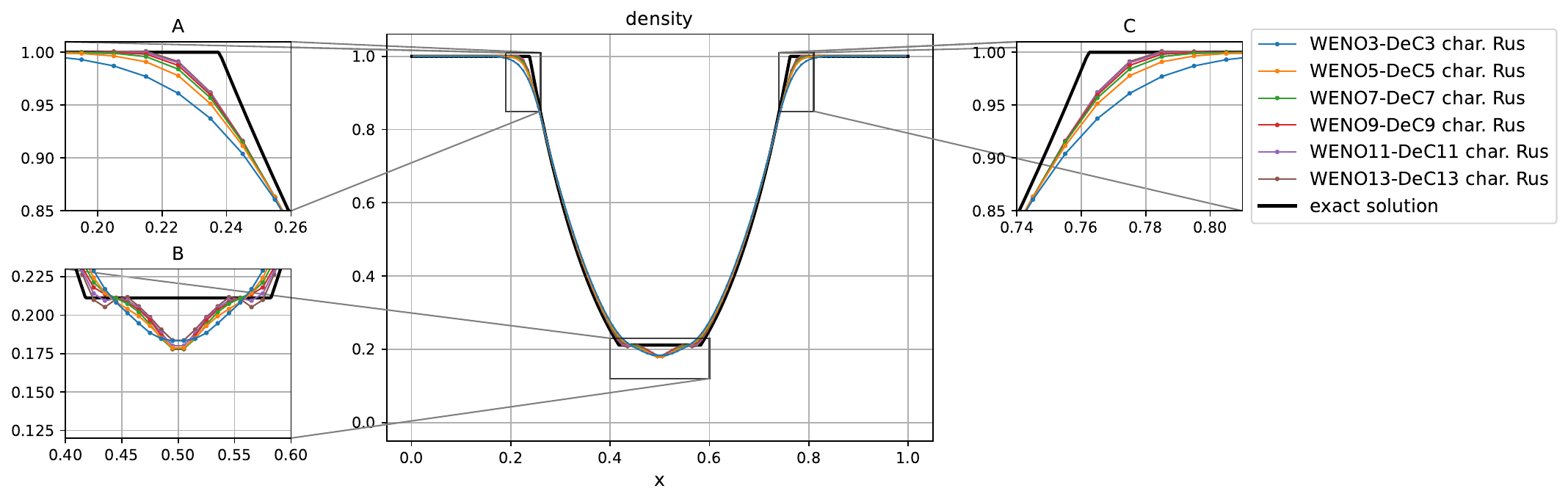}
		\caption{Reconstruction of characteristic variables and Rusanov}
	\end{subfigure}
	\\
	\begin{subfigure}[b]{0.9\textwidth}
		\centering
		\includegraphics[width=\textwidth]{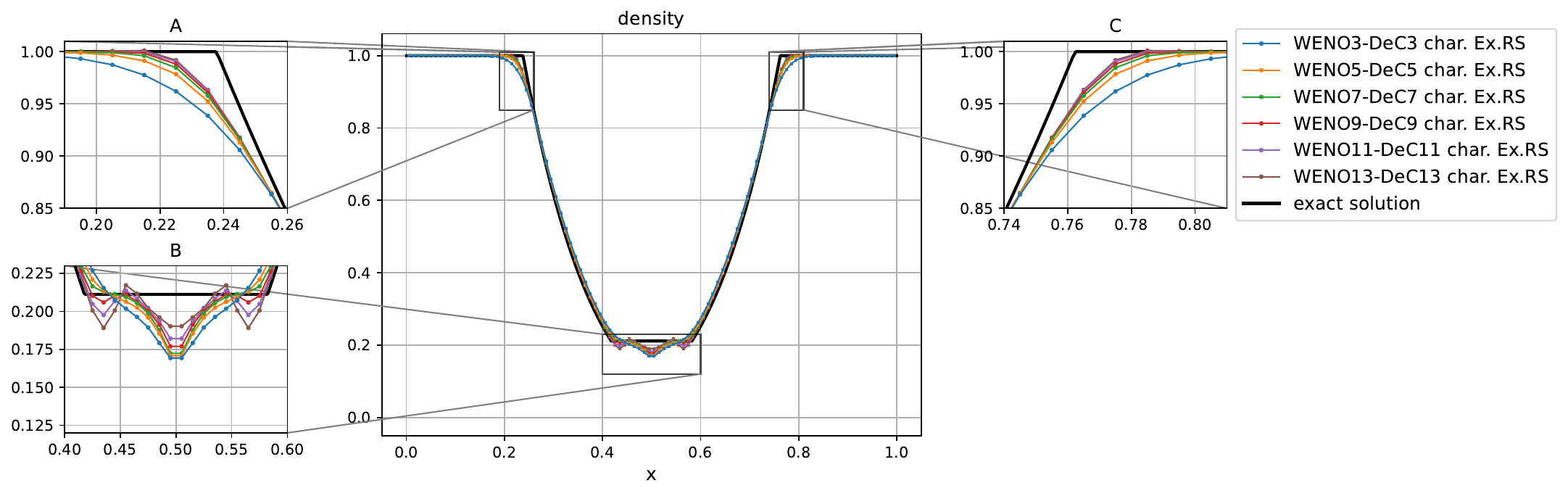}
		\caption{Reconstruction of characteristic variables and exact Riemann solver}
	\end{subfigure}
	\caption{Euler equations, relaxed Riemann problem 2: Results for the density obtained through WENO--DeC for reconstruction of characteristic variables with $C_{CFL}:=0.45$}
	\label{fig:Euler_1d_RP2_relaxed_zoom_density}
\end{figure}

In Figure~\ref{fig:Euler_1d_RP2_relaxed_zoom_density}, we report zooms on some details of the density for reconstruction of characteristic variables only; results for reconstruction of conserved variables can be found in the supplementary material.
The exact Riemann solver produces small spurious over- and undershoots at the foot of the rarefactions, which are smeared by Rusanov numerical flux.
Overall, the best results are obtained with reconstruction of characteristic variables and Rusanov.
In general, we can observe how higher order discretizations are able to provide results which are closer to the reference.

\subsubsection{Riemann problem 3 (Left Woodward--Colella)}\label{sec:Euler_1d_RP3}
This Riemann problem corresponds to the left part of the Woodward--Colella test~\cite{woodward1984numerical}.
The severe jump in the pressure, which amounts to 5 orders of magnitude, makes this test very challenging from the computational point of view.


For $C_{CFL}:=0.95$, all the simulations work up to order 9, besides when characteristic variables are reconstructed and Rusanov is used. 
In such a case all simulations crash.
Simulations work up to order 11 for $0.05 \leq C_{CFL}\leq 0.45$ with reconstruction of characteristic variables and exact Riemann solver, and up to order 9 for all the other settings. We did not test for smaller values of $C_{CFL}$.

\begin{figure}[htbp]
	\centering
	\begin{subfigure}[b]{0.9\textwidth}
		\centering
		\includegraphics[width=\textwidth]{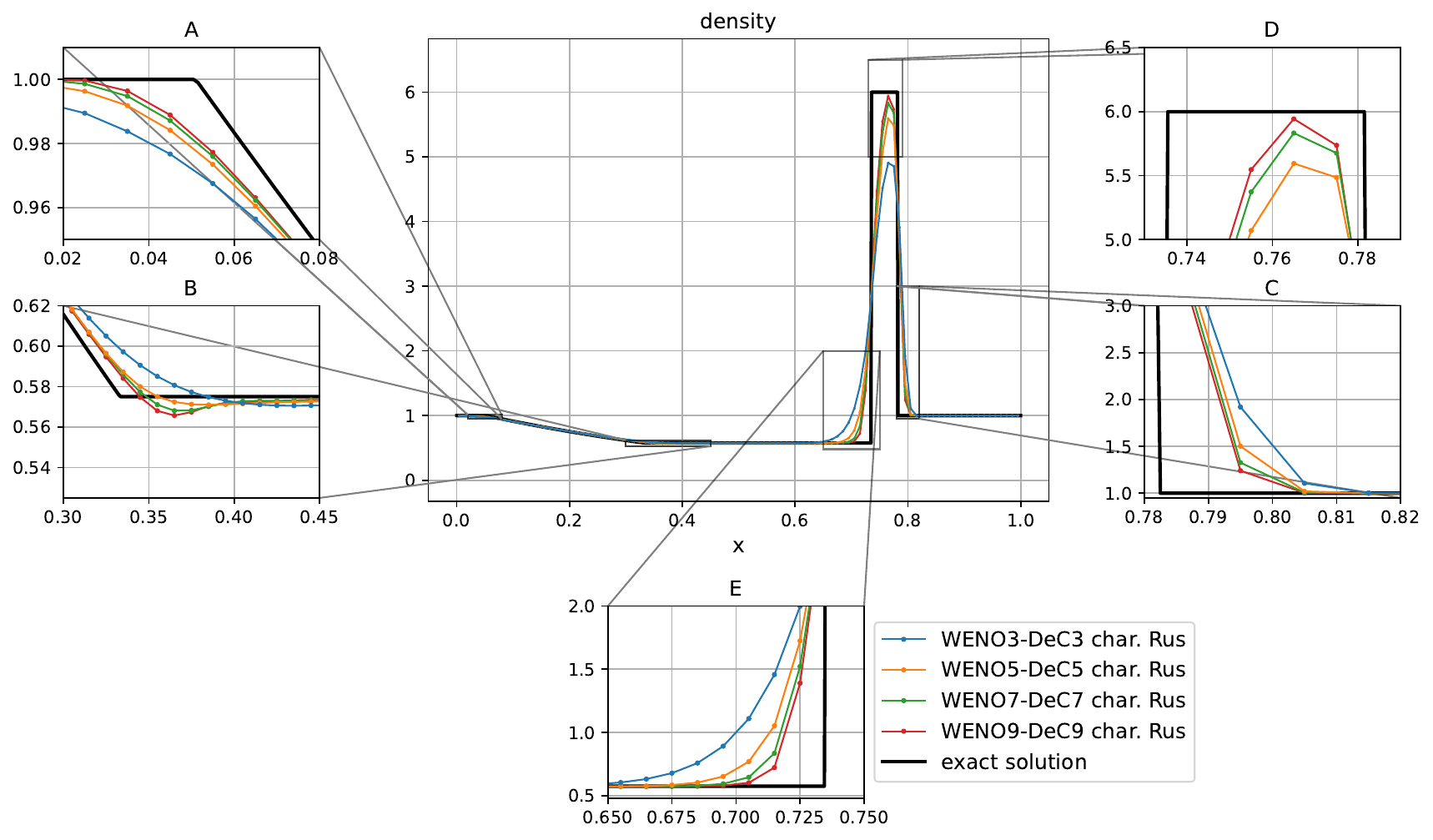}
		\caption{Reconstruction of characteristic variables and Rusanov}
	\end{subfigure}
	\\
	\begin{subfigure}[b]{0.9\textwidth}
		\centering
		\includegraphics[width=\textwidth]{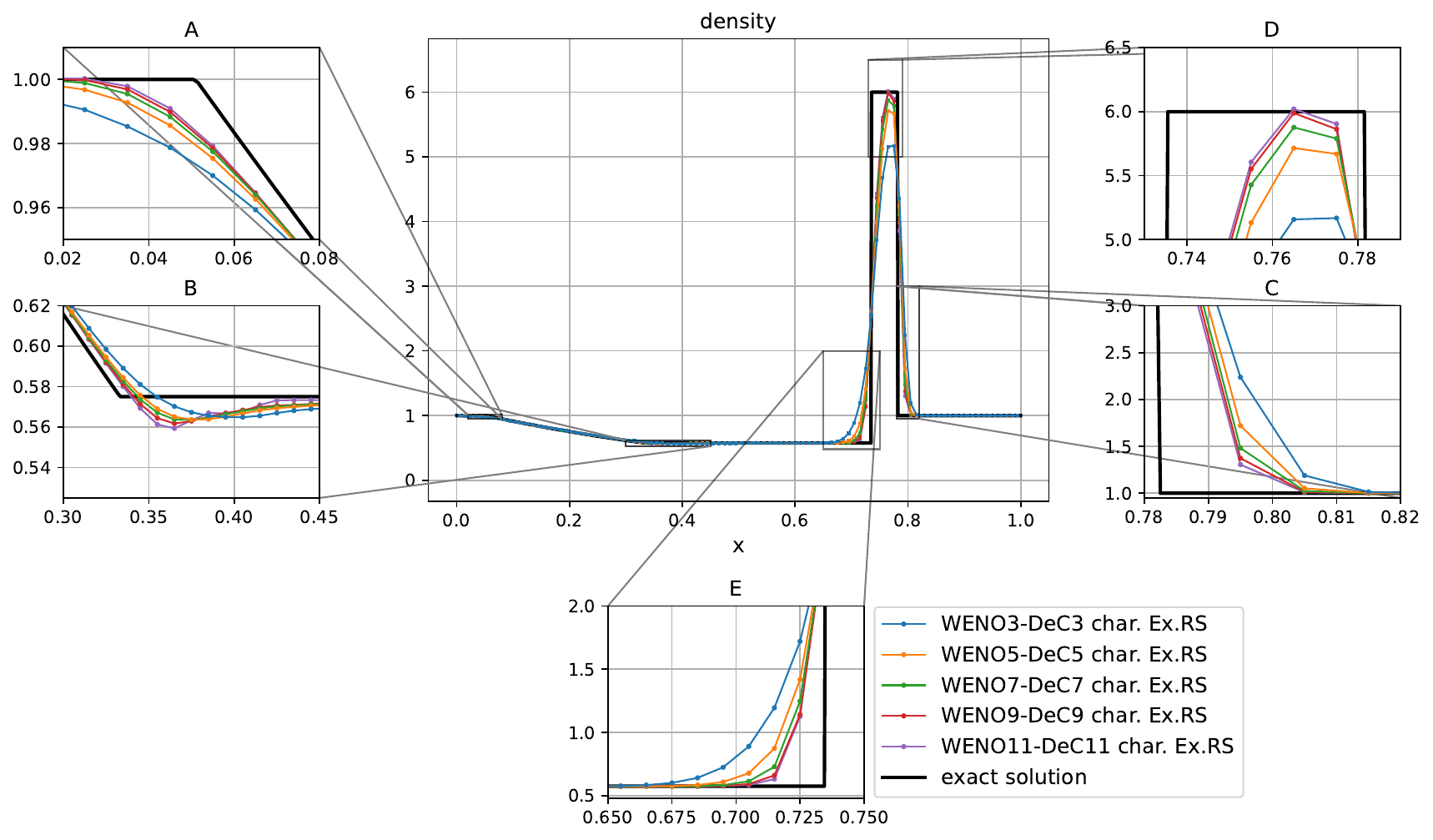}
		\caption{Reconstruction of characteristic variables and exact Riemann solver}
	\end{subfigure}
	\caption{Euler equations, Riemann problem 3: Results for the density obtained through WENO--DeC for reconstruction of characteristic variables with $C_{CFL}:=0.45$}
	\label{fig:Euler_1d_RP3_zoom_density}
\end{figure}

We perform the comparison for $C_{CFL}:= 0.45$ and 100 cells.
In Figure~\ref{fig:Euler_1d_RP3_zoom_density}, we report the density details for reconstruction of characteristic variables with both numerical fluxes; results obtained with reconstruction of conserved variables are available in the supplementary material.
Looking at the results in Figure~\ref{fig:Euler_1d_RP3_zoom_density}, in particular in panel D, we can see that the exact Riemann solver captures a bit better the peak of the density, especially for order 3. Differences among the numerical fluxes tend to disappear as the order of accuracy increases.
With no exceptions, increasing the order of the method is associated with higher quality of the obtained results.

\subsubsection{Riemann problem 4 (Collision of two shocks)}\label{sec:Euler_1d_RP4}
In this test, two streams collide creating a very high--pressure--area in the center of the domain.
Despite this, we experienced almost no problem in running this test.

\begin{figure}[htbp]
	\centering
	\begin{subfigure}[b]{0.9\textwidth}
		\centering
		\includegraphics[width=\textwidth]{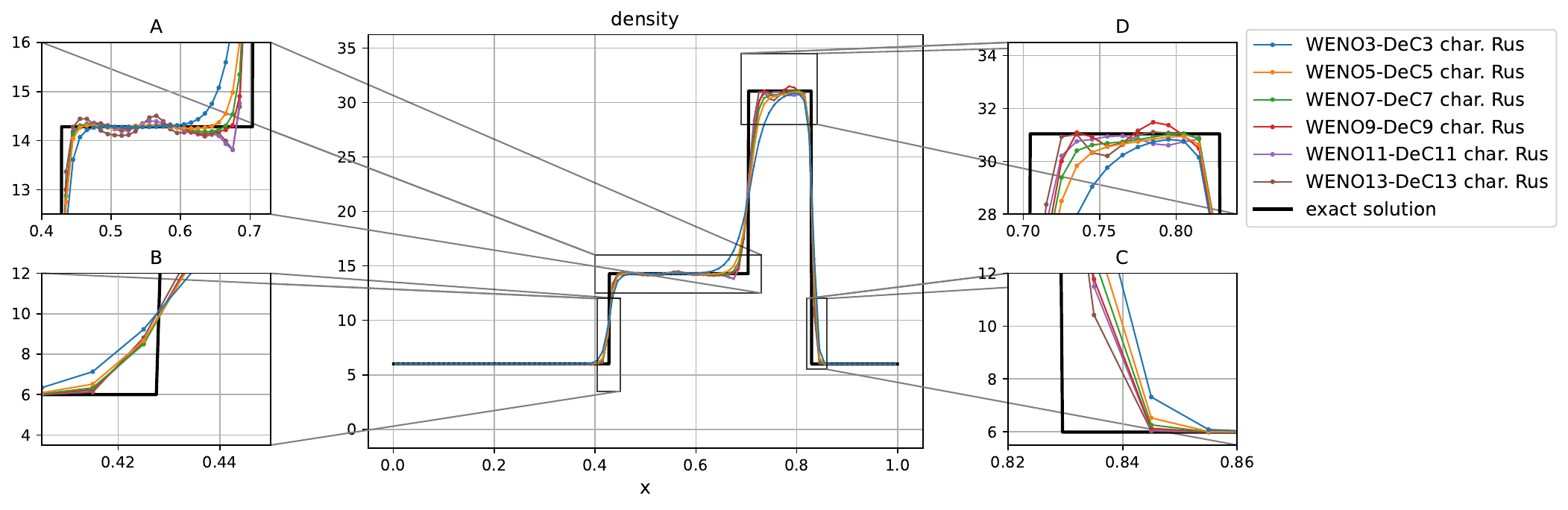}
		\caption{Reconstruction of characteristic variables and Rusanov}
	\end{subfigure}
	\\
	\begin{subfigure}[b]{0.9\textwidth}
		\centering
		\includegraphics[width=\textwidth]{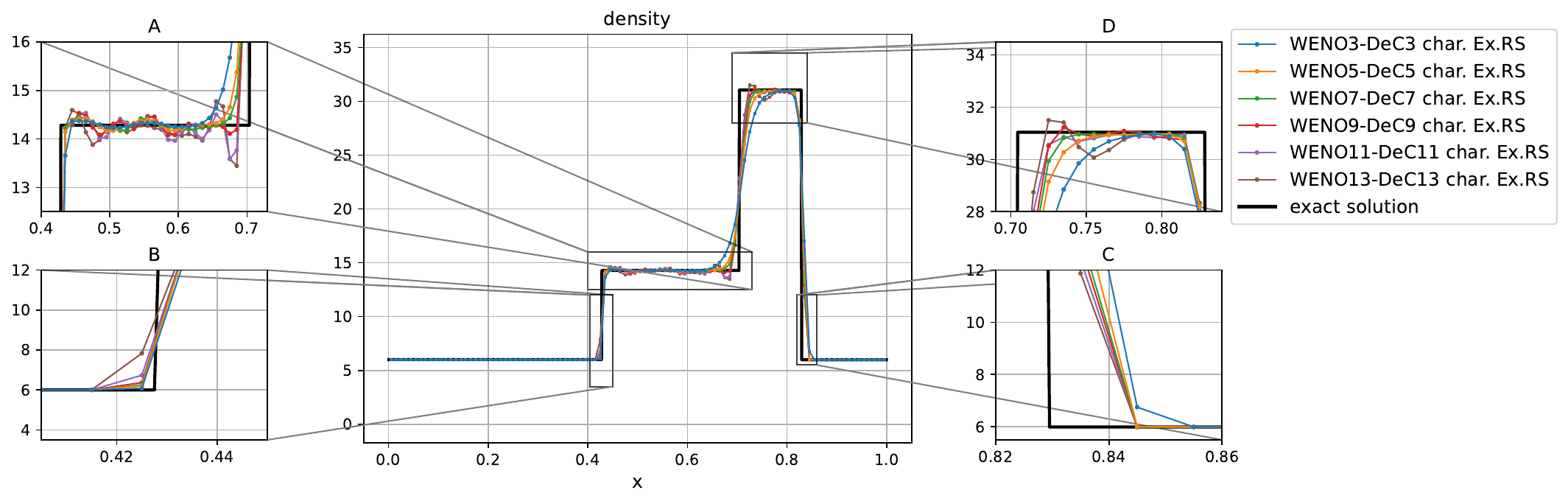}
		\caption{Reconstruction of characteristic variables and exact Riemann solver}
	\end{subfigure}
	\caption{Euler equations, Riemann problem 4: Results for the density obtained through WENO--DeC for reconstruction of characteristic variables with $C_{CFL}:=0.95$}
	\label{fig:Euler_1d_RP4_zoom_density}
\end{figure}

On computational meshes of 100 cells, all settings experience no crashes for all orders and $C_{CFL}:=0.95$, besides reconstruction of conserved variables with Rusanov numerical flux, which only works up to order 11.
We could not find a value of $C_{CFL}$ for which order 13 works with such setting. We tested up to $C_{CFL}:=0.05$.
We perform a comparison between the settings for $C_{CFL}:=0.95$, and we only report results for reconstruction of characteristic variables in Figure~\ref{fig:Euler_1d_RP4_zoom_density}. Also in this case, results for reconstruction of conserved variables are available in the supplementary material.
%
%
In this test, exact Riemann solver and Rusanov perform similarly.
Little spurious oscillations are present in both cases as can be seen from the density zooms.
In particular, they are slightly smaller for Rusanov, see panel A, due to its more diffusive character. On the other hand, the exact Riemann solver is able to more sharply capture the left shock as can bee seen from panel B.

\subsubsection{Riemann problem 5 (Stationary contact)}\label{sec:Euler_1d_RP5}
This Riemann problem can be regarded as a modification of the third one, in which the two states move towards left with rather high velocity.
As usual, we have considered 100 cells for the numerical experiments.

The behavior of the methods on this test was quite irregular.
Reconstruction of conserved variables with Rusanov always crashes for all orders and any value of $C_{CFL}$ between 0.95 and 0.05.
For $C_{CFL}:=0.95$, only reconstruction of characteristic variables with Rusanov for order 3 works. All other settings crash.
For $C_{CFL}:=0.75$, reconstruction of characteristic variables with Rusanov works up to order 7. All other settings crash.
For $C_{CFL}\leq 0.65$, reconstruction of characteristic variables with Rusanov works up to order 11. All other setting crash, for $0.35 \leq C_{CFL}\leq 0.65$, besides some irregular exceptions involving exact Riemann solver for orders 9, 11 and 13 for some values of $C_{CFL}$, for both variable reconstructions.
For $C_{CFL}:=0.3$, reconstruction of conserved variables with exact Riemann solver works for order higher than 5,
while, reconstruction of characteristic variables with exact Riemann solver works from order 7 to 11.
For $C_{CFL}:=0.2$, reconstruction of characteristic variables with exact Riemann solver works from order 5 to 11.
For the same $C_{CFL}$, reconstruction of conserved variables works for all orders with exact Riemann solver.
For $C_{CFL}:=0.1$, all settings work but reconstruction of characteristic variables with Rusanov for order 13 and reconstruction of conserved variables with Rusanov, which always crashes for all orders, as already anticipated.
Issues do not occur at the beginning of the simulations but rather at later times, and they are due to the reconstructions producing negative values of density or pressure. As already stated, fixing this problem is not in the goals of this paper and is left for future works.

\begin{figure}[htbp]
	\centering
	\begin{subfigure}[b]{0.9\textwidth}
		\centering
		\includegraphics[width=\textwidth]{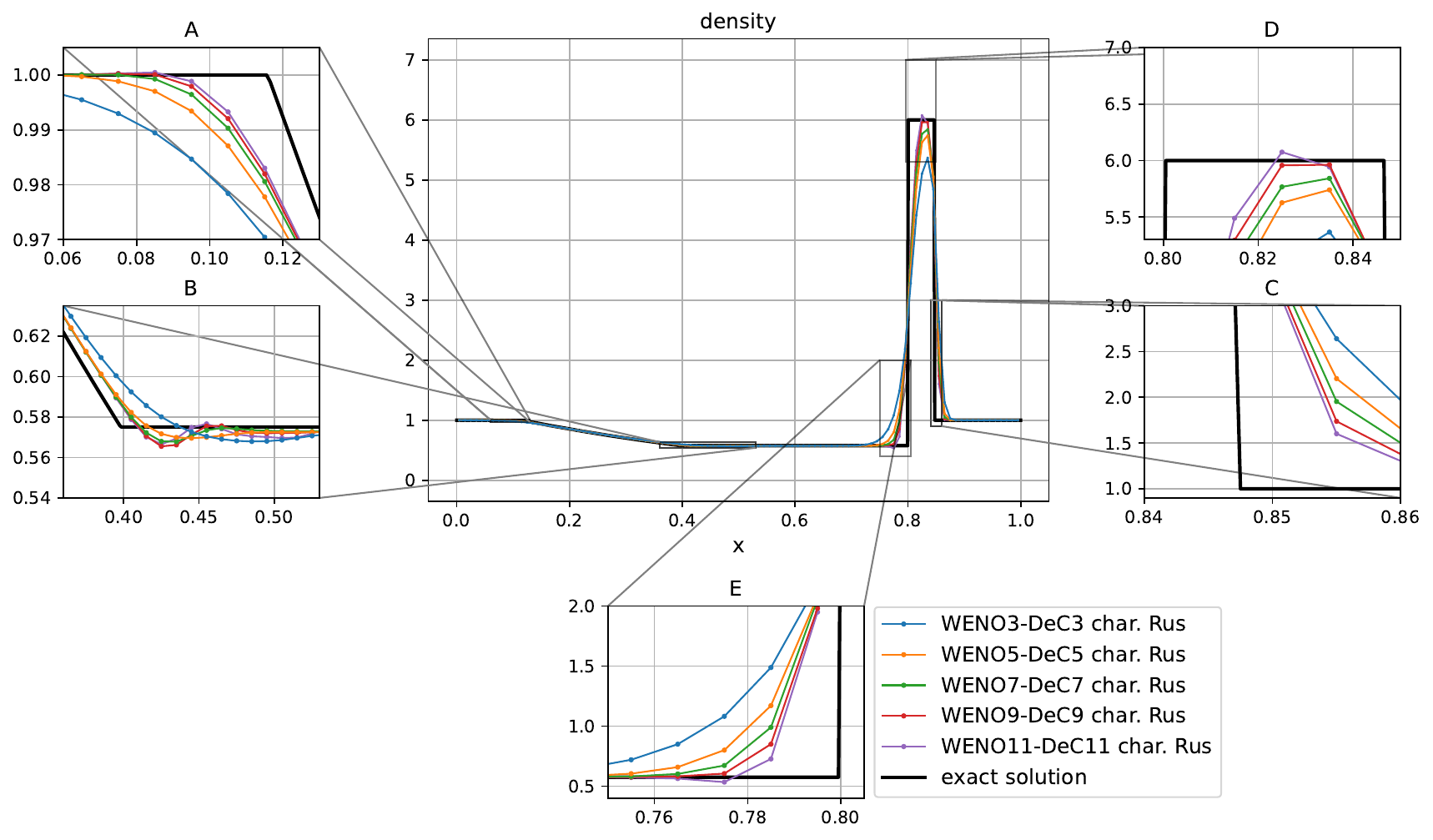}
		\caption{Reconstruction of characteristic variables and Rusanov}
	\end{subfigure}
	\\
	\begin{subfigure}[b]{0.9\textwidth}
		\centering
		\includegraphics[width=\textwidth]{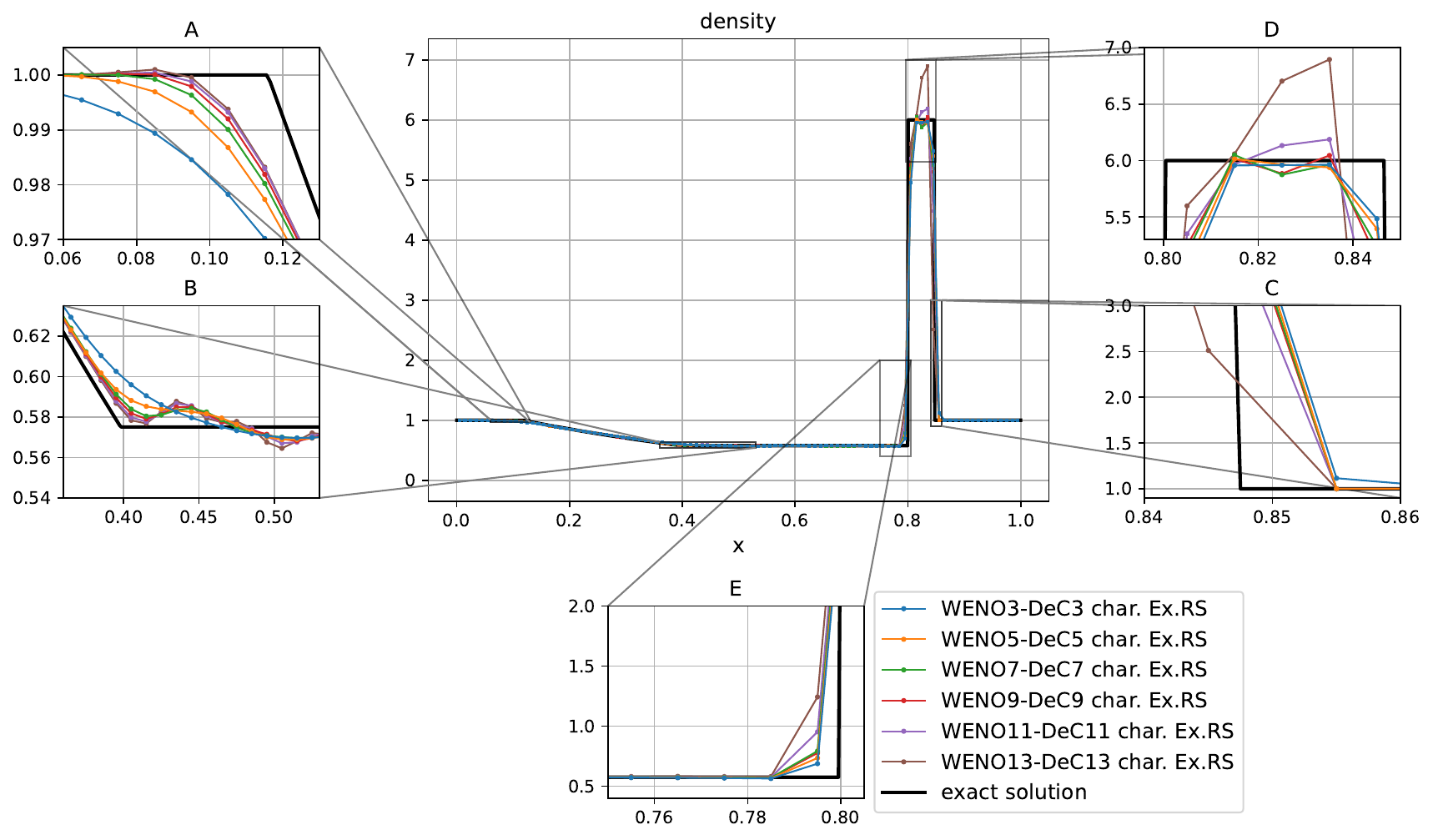}
		\caption{Reconstruction of characteristic variables and exact Riemann solver}
	\end{subfigure}
	\caption{Euler equations, Riemann problem 5: Results for the density obtained through WENO--DeC for reconstruction of characteristic variables with $C_{CFL}:=0.1$}
	\label{fig:Euler_1d_RP5_zoom_density}
\end{figure}

Let us comment the results obtained for $C_{CFL}:=0.1$ and reconstruction of characteristic variables, reported in Figure~\ref{fig:Euler_1d_RP5_zoom_density}; results for reconstruction of conserved variables are instead reported in the supplementary material.
We can appreciate a much higher resolution when the exact Riemann solver is employed.
The density peak is smeared by Rusanov, but is perfectly captured even by order 3 with reconstruction of conserved variables, when the exact Riemann solver is employed.
Results get better as the order increase, excluded order 13 which is characterized by ugly spurious overshoots not present for other orders.

\subsubsection{Modified shock--turbulence interaction}\label{sec:Euler_1d_titarev_toro}
This last test, introduced in~\cite{titarev2004finite}, consists in a challenging modification of the original problem presented in~\cite{shu1989efficient} by Shu and Osher.
On the computational domain $\Omega:=[-5,5]$, the initial conditions read 
\begin{align}
	\begin{pmatrix}
		\rho\\
		u\\
		p
	\end{pmatrix}(x,0):=\begin{cases}
		(1.515695,0.523346,1.80500)^T, \quad  &\text{if}~x< -4.5,\\
		(1.0+0.1\sin{(20 \pi x)},0,1)^T, \quad &\text{otherwise}.
	\end{cases} 
	\label{eq:Euler_1d_titarev_toro_IC}
\end{align}
Inflow and transmissive boundary conditions are prescribed at the left and right boundaries respectively. The final time is $T_f:=5.$

\begin{figure}[htbp]
	\centering
	\begin{subfigure}[b]{0.7\textwidth}
		\centering
		\includegraphics[width=\textwidth]{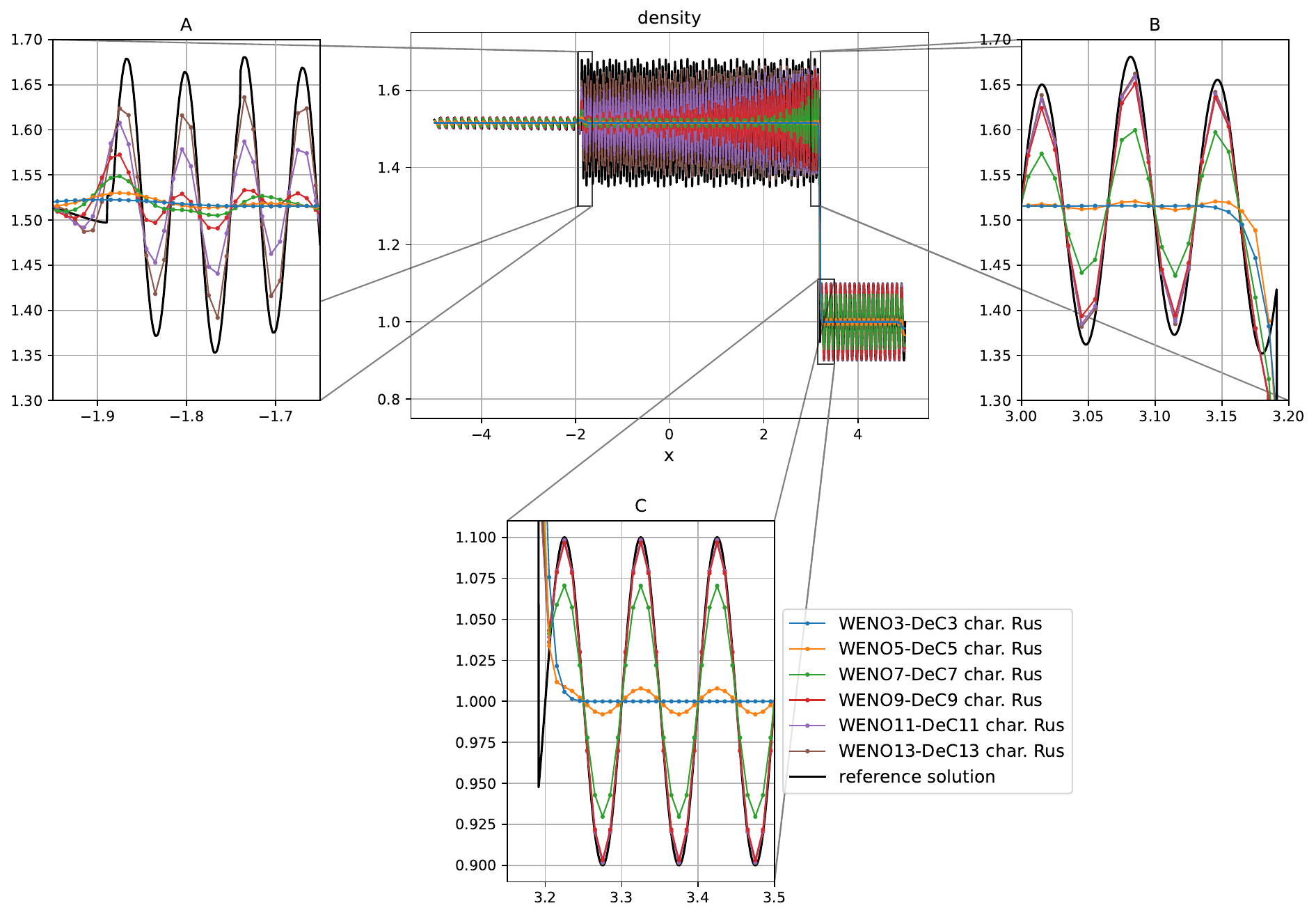}
		\caption{Reconstruction of characteristic variables and Rusanov}
	\end{subfigure}
	\\
	\begin{subfigure}[b]{0.7\textwidth}
		\centering
		\includegraphics[width=\textwidth]{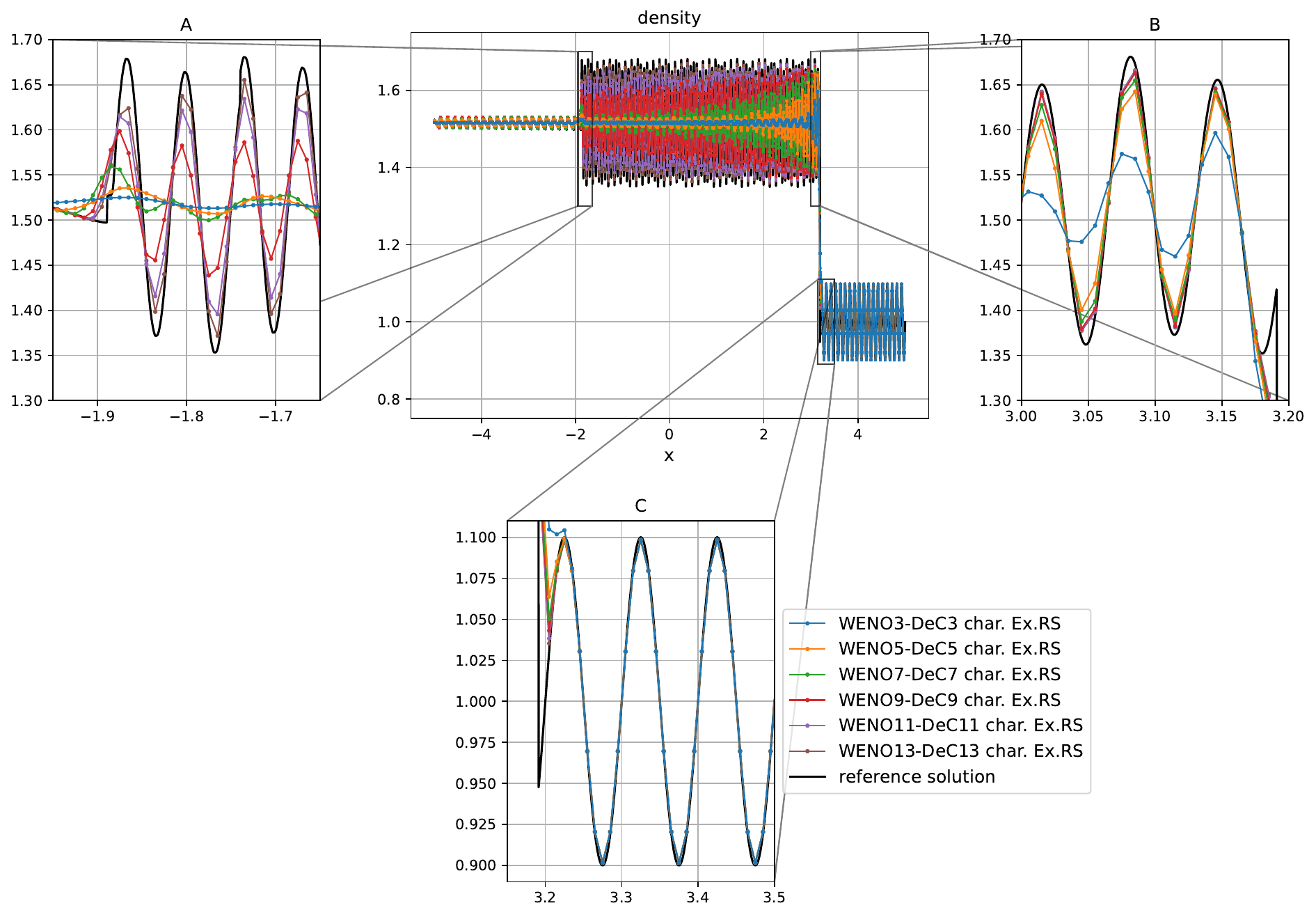}
		\caption{Reconstruction of characteristic variables and exact Riemann solver}
	\end{subfigure}
	\caption{Euler equations, Modified shock--turbulence interaction: Results for the density obtained through WENO--DeC for reconstruction of characteristic variables with $C_{CFL}:=0.95$}
	\label{fig:Euler_1d_titarev_toro_zoom_density}
\end{figure}

\begin{figure}[htbp]
	\centering
	\begin{subfigure}[b]{0.7\textwidth}
		\centering
		\includegraphics[width=\textwidth]{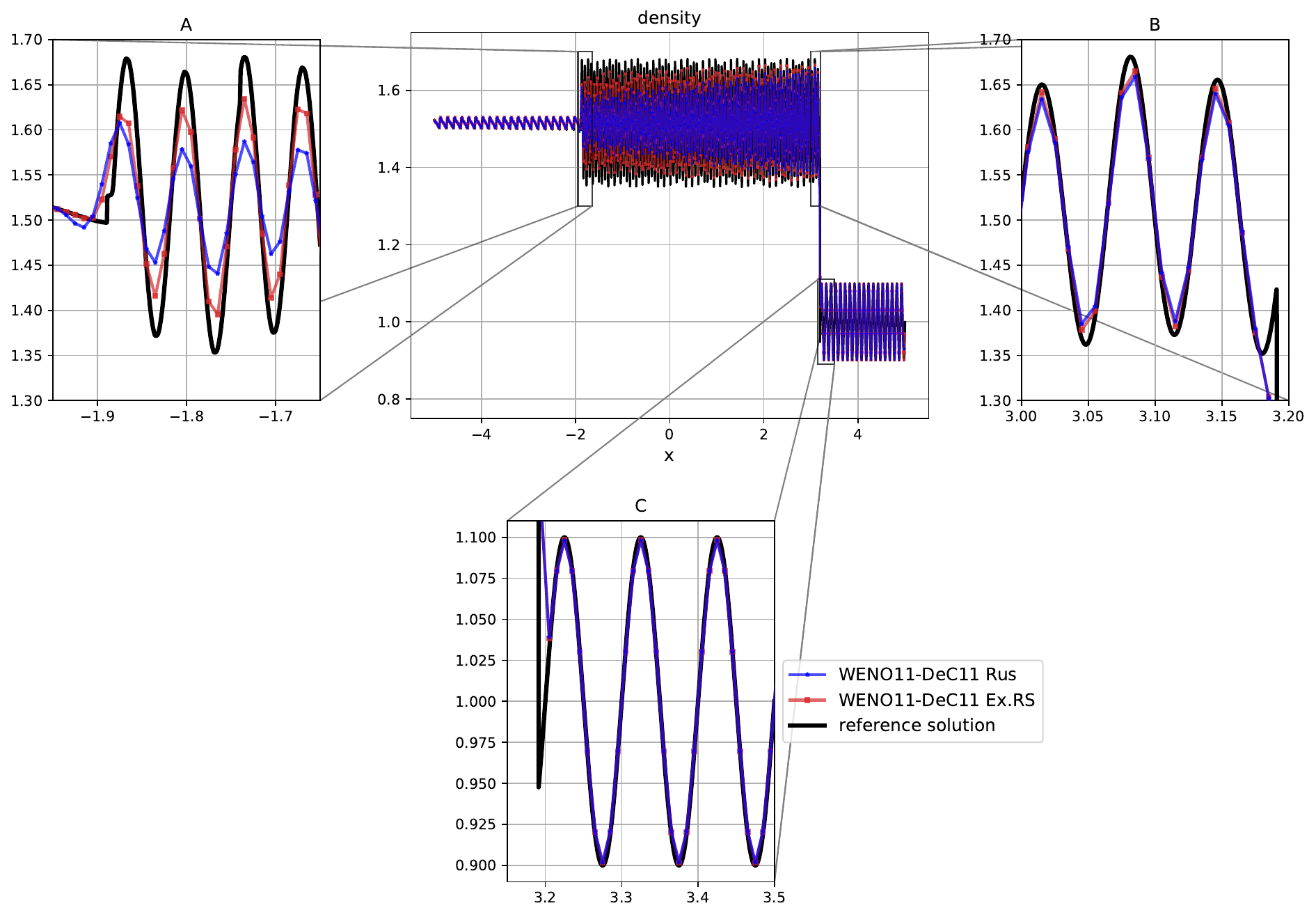}
		\caption{Order 11}
	\end{subfigure}
	\\
	\begin{subfigure}[b]{0.7\textwidth}
		\centering
		\includegraphics[width=\textwidth]{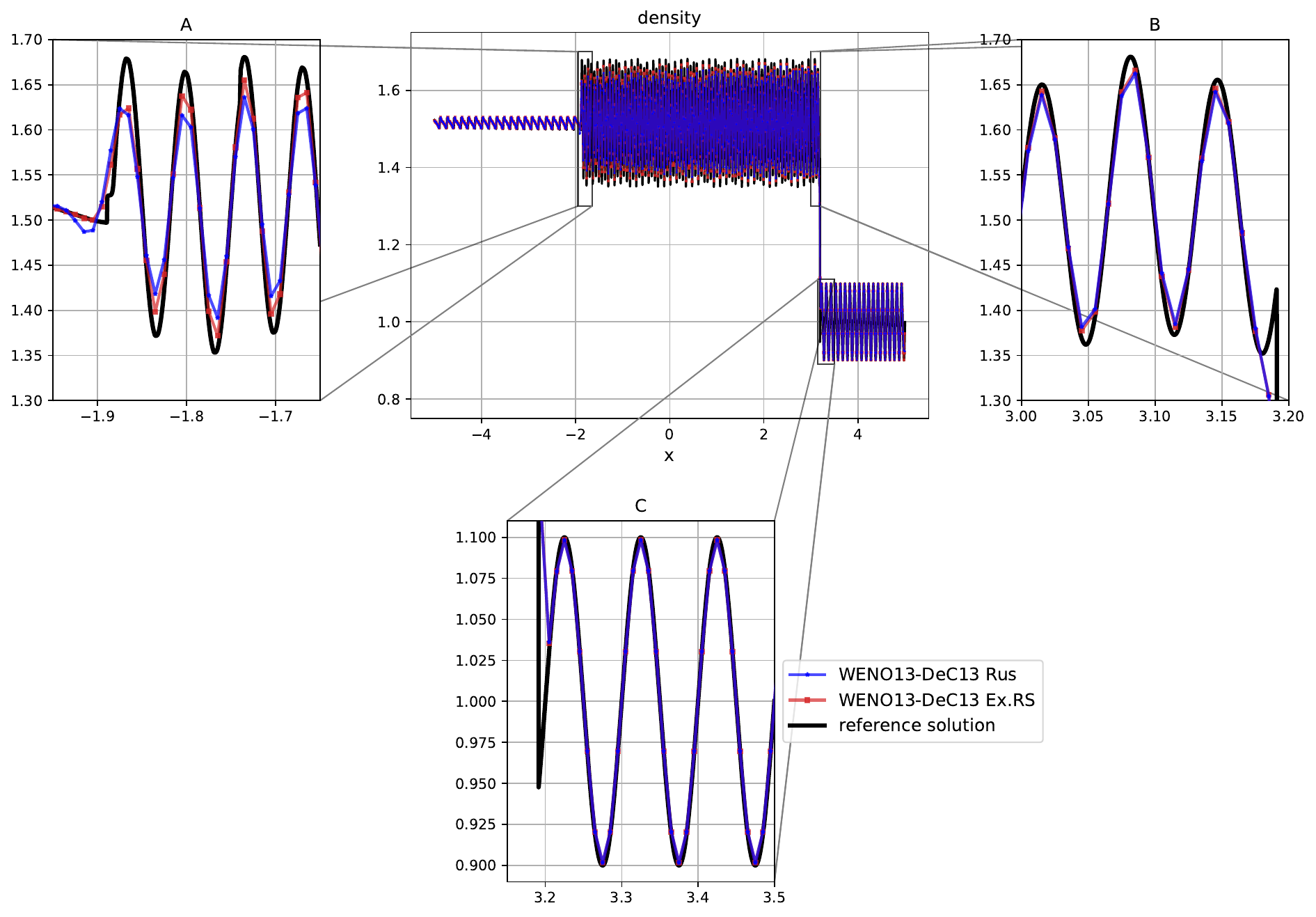}
		\caption{Order 13}
	\end{subfigure}
	\caption{Euler equations, Modified shock--turbulence interaction: Results for the density obtained through WENO--DeC for reconstruction of characteristic variables with $C_{CFL}:=0.95$ for Rusanov and exact Riemann solver}
	\label{fig:Euler_1d_titarev_toro_numerical_fluxes}
\end{figure}

Such a test is very challenging, as it consists of a shock interacting with a turbulent flow with very high frequency oscillations, hence, requiring at the same time a good capturing of smooth structures and a robust handling of discontinuities.
We considered 1000 cells and $C_{CFL}:=0.95$, and we experienced no simulation failures for any reconstructed variable, numerical flux and order.
The results obtained for the density with reconstruction of characteristic variables are reported in Figure~\ref{fig:Euler_1d_titarev_toro_zoom_density}, whereas, the ones obtained for reconstruction of conserved variables can be found in the supplementary material.
The quality of the results improves as the order of the method increases.
We can appreciate huge differences between Rusanov and exact Riemann solver, being the latter more able to capture many solution structures even for order 3, see panels B and, in particular, C.
Differences tend to decrease for increasing order but are still quite evident for order 9.
This is, in fact, the test in which the most evident differences between different numerical fluxes have been recorded in the investigation conducted in~\cite{micalizzi2024impact} considering 8 numerical fluxes up to order 7.
There, it was conjectured that such differences may eventually disappear increasing the order of accuracy to very high numbers.
This seems to be the case, as can be appreciated from Figure~\ref{fig:Euler_1d_titarev_toro_numerical_fluxes}, where we report the results for order 11 and 13 with reconstruction of characteristic variables.
The difference between the numerical fluxes, absolutely evident up to order 7 in~\cite{micalizzi2024impact}, becomes much smaller for such high orders.
In particular, there are no differences at all in panel C, and barely visible ones in panel B.
One can still appreciate some little differences in panel A, i.e., in the tail of the turbulence, with exact Riemann solver being more accurate than Rusanov.

Let us remark that, for the considered mesh refinement, even the 7-th order versions of the investigated schemes are dramatically far from convergence. In such a test, increasing the order to very high values is mandatory for a satisfactory description of the flow.


\subsubsection{Critical aspects}\label{sec:Euler_1d_critical_aspects}
The results obtained in the whole section put in evidence that WENO--DeC schemes, even for very high order, are able to tackle tough simulations involving strong discontinuities.
Nonetheless, spurious oscillations and simulation crashes can be experienced, as we showed in the context of Riemann problems 2 and 6. They are due to the fact that we have investigated the basic WENO--DeC framework without any further modification. Nothing at all is made in order to ensure, for example, positivity of density and pressure of the reconstructed values.
Special care under this point of view, along with other adaptive strategies such as a posteriori limiting, could prevent the failure of many simulations and reduce spurious oscillations, providing better results.
Investigations in this direction, e.g., with a posteriori limiting, are planned for future works.
Let us notice that we also tried adopting rigorous speed estimates for $s^x$ and $s^y$ for the computation of the time step in Equation~\ref{eq:CFL}. In particular, we used the rigorous bounds from~\cite{toro2020bounds} computed as detailed in~\cite[Section 9.5]{ToroBook}. As already stated, the simple direct estimates adopted herein do not bound the actual wave speeds, which may be the cause of simulation crashes for (apparently) safe values of $C_{CFL}$ even for low order schemes. We saw no improvement in robustness switching to rigorous speed estimates, which confirms the fact that simulation crashes are really due to negative values of density and pressure provided by the reconstruction.

Reconstructing characteristic variables is almost mandatory to avoid oscillatory results, whereas the preferability of a given numerical flux depends on the specific test.
Generally speaking, the exact Riemann solver provides sharper results, which can, however, present small spurious oscillations.
On the other hand, results obtained with Rusanov are in general slightly less accurate but also less oscillatory.
As conjectured in~\cite{micalizzi2024impact}, differences among numerical fluxes tend to decrease for increasing order, however, depending on the test, they may not disappear even for order 13.
Other tests not reported here, involving stationary and moving contact discontinuities, put in evidence that also for very high order the choice of the numerical flux still influences the quality of the obtained results.
We leave for future works the investigation of this aspect.

\section{Conclusions and future perspectives}\label{sec:conclusions}
In this work, we have investigated the WENO--DeC arbitrary high order framework for the numerical solution of hyperbolic PDEs.
In particular, we have focused on very high order realizations of such a framework, i.e., with order of spatial and temporal discretization higher than 5.
We tested the methods on several benchmarks for one--dimensional LAE and Euler equations.
Numerical results confirm the advantages of adopting very high order methods both on smooth and non--smooth problems.
In the first case, they are able to achieve smaller errors within coarser discretizations and shorter computational times; in the latter case, they provide a sharper capturing of various structures of the analytical solution.
%
%
We have also shown how the common practice of coupling high order space discretizations with lower order time integrations results in efficiency loss with respect to ``truly'' high order frameworks, such as WENO--DeC, where the order of the time discretization matches the one of the space discretization.

We plan to extend this investigation in several directions, for example, with multidimensional studies and coupling with other limiting techniques. 
In fact, on tough tests, spurious oscillations and simulation crashes can occur and further work is required to make the schemes fail--safe.

\subsubsection*{Acknowledgments}
LM has been funded by the LeRoy B. Martin, Jr. Distinguished Professorship Foundation. 
EFT gratefully acknowledges the support from the ``Special Programme on Numerical Methods for Non-linear Hyperbolic Partial Differential Equations'', sponsored by the Department of Mathematics and the Shenzhen International Center for Mathematics, Southern University of Science and Technology, Shenzhen, People Republic of China.
Both the authors would like to express gratitude to Prof. Alina Chertock from North Carolina State University, Prof. Alexander Kurganov from Southern University of Science and Technology and Prof. Chi--Wang Shu from Brown University, the organizers of the Workshop on Numerical Methods for Shallow Water Models, held in Shenzhen, May 11-15 2024. This investigation is the result of interesting discussions between the authors held during such event.
Furthermore, the authors acknowledge Prof. Rémi Abgrall from the University of Zurich for granting access to useful computational resources, and Prof. David Ketcheson from King Abdullah University of Science and Technology for the computation of more accurate coefficients of SSPRK4.

\appendix
\section{Accurate coefficients of SSPRK(5,4)}\label{app:SSPRK4}
The following coefficients for SSPRK(5,4) have been computed up to accuracy $10^{-30}$.

\begin{align}
\scalebox{0.5}{$
	A := \begin{pmatrix}
0 & 0 & 0 & 0 & 0 \\
0.391752226869253785640632115627 & 0 & 0 & 0 & 0  \\
0.217669096357834985920253802915 & 0.368410592709066783214662112772 & 0 & 0 & 0  \\
0.0826920866830935842609242437786 & 0.139958502107426395108400626025 & 0.251891774371960822884363746140 & 0 & 0  \\
0.0679662835740483884329695316049 & 0.115034698453668419467815057942 & 0.207034898772936576352392025561 & 0.544974750295139481064416383368 & 0 
	\end{pmatrix},
	$}
\end{align}

\begin{align}
	\scalebox{0.5}{$
		\uvec{b} := \begin{pmatrix}
			0.146811876157875933686947006683\\ 
			0.248482909391317264243714136087\\ 
			0.104258830279481225354037031167\\ 
			0.274438901048480694917546480567\\ 
			0.226007483122844881797755345495
		\end{pmatrix},\quad
		\uvec{c} := \begin{pmatrix}
		0\\
		0.391752226869253785640632115627\\\textbf{}
		0.586079689066901769134915915687\\
		0.4745423631624808022536886159436\\
		0.9350106310957928653175929984759
		\end{pmatrix}
		$}
\end{align}

\printbibliography

@article{huang2014re,
  title={A re-averaged WENO reconstruction and a third order CWENO scheme for hyperbolic conservation laws},
  author={Huang, Chieh-Sen and Arbogast, Todd and Hung, Chen-Hui},
  journal={Journal of Computational Physics},
  volume={262},
  pages={291--312},
  year={2014},
  publisher={Elsevier}
}

@article{zhu2017new,
  title={A new type of finite volume WENO schemes for hyperbolic conservation laws},
  author={Zhu, Jun and Qiu, Jianxian},
  journal={Journal of Scientific Computing},
  volume={73},
  pages={1338--1359},
  year={2017},
  publisher={Springer}
}

@incollection{zhang2016eno,
  title={ENO and WENO schemes},
  author={Zhang, Y-T and Shu, C-W},
  booktitle={Handbook of numerical analysis},
  volume={17},
  pages={103--122},
  year={2016},
  publisher={Elsevier}
}

@article{zhang2024application,
  title={Application of a Novel High-Order WENO Scheme in LES Simulations},
  author={Zhang, Shuo and Zhong, Dongdong and Wang, Hao and Wu, Xingshuang and Ge, Ning},
  journal={Applied Sciences},
  volume={14},
  number={17},
  pages={7875},
  year={2024},
  publisher={MDPI}
}

@article{wang2025extremum,
  title={An extremum properties (EP)-based discontinuous sensor and hybrid weighted essentially non-oscillatory scheme on tetrahedral meshes},
  author={Wang, Zhenming and Zhu, Jun and Tan, Yan and Tian, Linlin and Zhao, Ning},
  journal={Journal of Computational Physics},
  pages={113979},
  year={2025},
  publisher={Elsevier}
}

@article{tang2005sonic,
  title={On the sonic point glitch},
  author={Tang, Huazhong},
  journal={Journal of Computational Physics},
  volume={202},
  number={2},
  pages={507--532},
  year={2005},
  publisher={Elsevier}
}

@article{zhu2008runge,
  title={Runge--Kutta discontinuous Galerkin method using WENO limiters II: unstructured meshes},
  author={Zhu, Jun and Qiu, Jianxian and Shu, Chi-Wang and Dumbser, Michael},
  journal={Journal of Computational Physics},
  volume={227},
  number={9},
  pages={4330--4353},
  year={2008},
  publisher={Elsevier}
}

@article{wolf2007high,
  title={High-order ENO and WENO schemes for unstructured grids},
  author={Wolf, WR and Azevedo, JLF},
  journal={International Journal for Numerical Methods in Fluids},
  volume={55},
  number={10},
  pages={917--943},
  year={2007},
  publisher={Wiley Online Library}
}

@article{micalizzi2024impact,
  title={Impact of Numerical Fluxes on High Order Semidiscrete WENO--DeC Finite Volume Schemes},
  author={Micalizzi, Lorenzo and Toro, Eleuterio Francisco},
  journal={arXiv e-prints},
  pages={arXiv--2411},
  year={2024}
}

@article{toro2020bounds,
  title={Bounds for wave speeds in the Riemann problem: Direct theoretical estimates},
  author={Toro, Eleuterio F and M{\"u}ller, Lucas O and Siviglia, Annunziato},
  journal={Computers \& Fluids},
  volume={209},
  pages={104640},
  year={2020},
  publisher={Elsevier}
}

@article{jiang1996efficient,
  title={Efficient implementation of weighted ENO schemes},
  author={Jiang, Guang-Shan and Shu, Chi-Wang},
  journal={Journal of computational physics},
  volume={126},
  number={1},
  pages={202--228},
  year={1996},
  publisher={Elsevier}
}

@article{harten1987uniformly1,
  title={Uniformly high-order accurate nonoscillatory schemes. I},
  author={Harten, Ami and Osher, Stanley},
  journal={SIAM Journal on Numerical Analysis},
  volume={24},
  number={2},
  pages={279--309},
  year={1987},
  publisher={SIAM}
}

@article{harten1987uniformly2,
title = {Uniformly high order accurate essentially non-oscillatory schemes, III},
journal = {Journal of Computational Physics},
volume = {71},
number = {2},
pages = {231-303},
year = {1987},
issn = {0021-9991},
doi = {https://doi.org/10.1016/0021-9991(87)90031-3},
url = {https://www.sciencedirect.com/science/article/pii/0021999187900313},
author = {Ami Harten and Bjorn Engquist and Stanley Osher and Sukumar R Chakravarthy}
}

@misc{godlewski2021numerical,
  title={Numerical Approximation of Hyperbolic Systems of Conservation Laws},
  author={Godlewski, Edwige and Raviart, Pierre-Arnaud},
  year={2021},
  publisher={Springer New York, 2nd edition}
}

@article{toro1989fast,
  title={A fast Riemann solver with constant covolume applied to the random choice method},
  author={Toro, Eleuterio Francisco},
  journal={International journal for numerical methods in fluids},
  volume={9},
  number={9},
  pages={1145--1164},
  year={1989},
  publisher={Wiley Online Library}
}

@book{hirsch2007numerical,
  title={Numerical computation of internal and external flows: The fundamentals of computational fluid dynamics},
  author={Hirsch, Charles},
  year={2007},
  publisher={Elsevier}
}

@misc{toro2024computational,
  title={Computational Algorithms for Shallow Water Equations},
  author={Toro, Eleuterio F},
  year={2024},
  publisher={Springer}
}

@article{qiu2002construction,
  title={On the construction, comparison, and local characteristic decomposition for high-order central WENO schemes},
  author={Qiu, Jianxian and Shu, Chi-Wang},
  journal={Journal of Computational Physics},
  volume={183},
  number={1},
  pages={187--209},
  year={2002},
  publisher={Elsevier}
}

@article{miyoshi2020short,
  title={A short note on reconstruction variables in shock capturing schemes for magnetohydrodynamics},
  author={Miyoshi, Takahiro and Minoshima, Takashi},
  journal={Journal of Computational Physics},
  volume={423},
  pages={109804},
  year={2020},
  publisher={Elsevier}
}

@article{peng2019adaptive,
  title={An adaptive characteristic-wise reconstruction WENO-Z scheme for gas dynamic Euler equations},
  author={Peng, Jun and Zhai, Chuanlei and Ni, Guoxi and Yong, Heng and Shen, Yiqing},
  journal={Computers \& Fluids},
  volume={179},
  pages={34--51},
  year={2019},
  publisher={Elsevier}
}

@article{ghosh2012compact,
  title={Compact reconstruction schemes with weighted ENO limiting for hyperbolic conservation laws},
  author={Ghosh, Debojyoti and Baeder, James D},
  journal={SIAM Journal on Scientific Computing},
  volume={34},
  number={3},
  pages={A1678--A1706},
  year={2012},
  publisher={SIAM}
}

@article{qiu2006numerical,
  title={A numerical study for the performance of the Runge--Kutta discontinuous Galerkin method based on different numerical fluxes},
  author={Qiu, Jianxian and Khoo, Boo Cheong and Shu, Chi-Wang},
  journal={Journal of Computational Physics},
  volume={212},
  number={2},
  pages={540--565},
  year={2006},
  publisher={Elsevier}
}

@article{leidi2024performance,
  title={Performance of high-order Godunov-type methods in simulations of astrophysical low Mach number flows},
  author={Leidi, Giovanni and Andrassy, R and Barsukow, Wasilij and Higl, J and Edelmann, Philipp VF and R{\"o}pke, Friedrich K},
  journal={Astronomy \& Astrophysics},
  volume={686},
  pages={A34},
  year={2024},
  publisher={EDP Sciences}
}

@article{qiu2007numerical,
  title={A numerical comparison of the Lax--Wendroff discontinuous Galerkin method based on different numerical fluxes},
  author={Qiu, Jianxian},
  journal={Journal of Scientific Computing},
  volume={30},
  pages={345--367},
  year={2007},
  publisher={Springer}
}

@article{qiu2008development,
  title={Development and comparison of numerical fluxes for LWDG methods},
  author={Qiu, Jianxian},
  journal={Numerical Mathematics, Theory, Methods and Application},
  volume={1},
  pages={1--32},
  year={2008}
}

@article{hongxia2020numerical,
  title={A Numerical Comparison of the HWENO Method Based on Different Numerical Fluxes},
  author={Hongxia, Liu and Yiming, Liu and He, Juan},
  journal={Authorea Preprints},
  year={2020},
  publisher={Authorea}
}

@article{san2015evaluation,
  title={Evaluation of Riemann flux solvers for WENO reconstruction schemes: Kelvin--Helmholtz instability},
  author={San, Omer and Kara, Kursat},
  journal={Computers \& Fluids},
  volume={117},
  pages={24--41},
  year={2015},
  publisher={Elsevier}
}

@article{toro2000centred,
  title={Centred TVD schemes for hyperbolic conservation laws},
  author={Toro, Eleuterio Francisco and Billett, SJ},
  journal={IMA Journal of Numerical Analysis},
  volume={20},
  number={1},
  pages={47--79},
  year={2000},
  publisher={OUP}
}

@incollection{chen2003centred,
  title={Centred Schemes for Non-Linear Hyperbolic Equations},
  author={Chen, GQ and Toro, Eleuterio Francisco},
  booktitle={Isaac Newton Institute for Mathematical Sciences Preprint Series NI 03046-NPA},
  year={2003}
}

@article{toro1992restoration,
  title={Restoration of the contact surface in the HLL-Riemann solver},
  author={Toro, Eleuterio F and Spruce, Michael and Speares, William},
  journal={Technical Report CoA 9204, Department of Aerospace Science, College of Aeronautics, Cranfield Institute of Technology, UK},
  year={1992}
}

@article{toro1994restoration,
  title={Restoration of the contact surface in the HLL-Riemann solver},
  author={Toro, Eleuterio F and Spruce, Michael and Speares, William},
  journal={Shock waves},
  volume={4},
  pages={25--34},
  year={1994},
  publisher={Springer}
}

@article{kurganov2023new,
  title={New low-dissipation central-upwind schemes},
  author={Kurganov, Alexander and Xin, Ruixiao},
  journal={Journal of Scientific Computing},
  volume={96},
  number={2},
  pages={56},
  year={2023},
  publisher={Springer}
}

@article{kurganov2001semidiscrete,
  title={Semidiscrete central-upwind schemes for hyperbolic conservation laws and Hamilton--Jacobi equations},
  author={Kurganov, Alexander and Noelle, Sebastian and Petrova, Guergana},
  journal={SIAM Journal on Scientific Computing},
  volume={23},
  number={3},
  pages={707--740},
  year={2001},
  publisher={SIAM}
}

@article{kurganov2000new,
  title={New high-resolution central schemes for nonlinear conservation laws and convection--diffusion equations},
  author={Kurganov, Alexander and Tadmor, Eitan},
  journal={Journal of computational physics},
  volume={160},
  number={1},
  pages={241--282},
  year={2000},
  publisher={Elsevier}
}

@article{harten1983upstream,
  title={On upstream differencing and Godunov-type schemes for hyperbolic conservation laws},
  author={Harten, Amiram and Lax, Peter D and Leer, Bram van},
  journal={SIAM review},
  volume={25},
  number={1},
  pages={35--61},
  year={1983},
  publisher={SIAM}
}

@techreport{Toro1996,
  author = {E. F. Toro},
  title = {On Glimm--Related Schemes for Conservation Laws},
  institution = {Department of Mathematics and Physics, Manchester Metropolitan University},
  year = {1996},
  number = {MMU--9602},
  address = {UK}
}

@article{Rusanov1961,
  author = {V. V. Rusanov},
  title = {Calculation of Interaction of Non–Steady Shock Waves with Obstacles},
  journal = {J. Comput. Math. Phys. USSR},
  volume = {1},
  pages = {267--279},
  year = {1961}
}

@article{popov2024space,
  title={Space--time adaptive ADER-DG finite element method with LST-DG predictor and a posteriori sub-cell ADER-WENO finite-volume limiting for multidimensional detonation waves simulation},
  author={Popov, IS},
  journal={Computers \& Fluids},
  pages={106425},
  year={2024},
  publisher={Elsevier}
}

@article{spiteri2002new,
  title={A new class of optimal high-order strong-stability-preserving time discretization methods},
  author={Spiteri, Raymond J and Ruuth, Steven J},
  journal={SIAM Journal on Numerical Analysis},
  volume={40},
  number={2},
  pages={469--491},
  year={2002},
  publisher={SIAM}
}

@article{rogerson1990numerical,
  title={A numerical study of the convergence properties of ENO schemes},
  author={Rogerson, AM and Meiburg, E},
  journal={Journal of Scientific Computing},
  volume={5},
  pages={151--167},
  year={1990},
  publisher={Springer}
}

@article{shu1990numerical,
  title={Numerical experiments on the accuracy of ENO and modified ENO schemes},
  author={Shu, Chi-Wang},
  journal={Journal of Scientific Computing},
  volume={5},
  number={2},
  pages={127--149},
  year={1990},
  publisher={Springer}
}

@article{woodward1984numerical,
  title={The numerical simulation of two-dimensional fluid flow with strong shocks},
  author={Woodward, Paul and Colella, Phillip},
  journal={Journal of computational physics},
  volume={54},
  number={1},
  pages={115--173},
  year={1984},
  publisher={Elsevier}
}

@misc{AbgrallMishranotes,
author = {Mishra, S. and Fjordholm, U. S. and Abgrall, R.},
title = {Numerical methods for conservation laws and related equations},
url = {https://metaphor.ethz.ch/x/2019/hs/401-4671-00L/literature/mishra_hyperbolic_pdes.pdf}
}

@article{titarev2004finite,
  title={Finite-volume WENO schemes for three-dimensional conservation laws},
  author={Titarev, Vladimir A and Toro, Eleuterio F},
  journal={Journal of Computational Physics},
  volume={201},
  number={1},
  pages={238--260},
  year={2004},
  publisher={Elsevier}
}

@article{toro1999numerica,
  title={NUMERICA: A Library of Source Codes for Teaching, Research and Applications},
  author={Toro, Eleuterio Francisco},
  year={1999}
}

@article{abgrall2024staggered,
  title={Staggered Schemes for Compressible Flow: A General Construction},
  author={Abgrall, Remi},
  journal={SIAM Journal on Scientific Computing},
  volume={46},
  number={1},
  pages={A399--A428},
  year={2024},
  publisher={SIAM}
}

@article{liu1994weighted,
  title={Weighted essentially non-oscillatory schemes},
  author={Liu, Xu-Dong and Osher, Stanley and Chan, Tony},
  journal={Journal of computational physics},
  volume={115},
  number={1},
  pages={200--212},
  year={1994},
  publisher={Elsevier}
}

@article{shi2002technique,
  title={A technique of treating negative weights in WENO schemes},
  author={Shi, Jing and Hu, Changqing and Shu, Chi-Wang},
  journal={Journal of Computational Physics},
  volume={175},
  number={1},
  pages={108--127},
  year={2002},
  publisher={Elsevier}
}

@article{shu1998essentially,
  title={Essentially non-oscillatory and weighted essentially non-oscillatory schemes for hyperbolic conservation laws},
author={Shu, Chi-Wang},
journal={Advanced numerical approximation of nonlinear hyperbolic equations},
pages={325--432},
year={1998},
publisher={Springer}
}

@article{shu1989efficient,
  title={Efficient implementation of essentially non-oscillatory shock-capturing schemes, II},
  author={Shu, Chi-Wang and Osher, Stanley},
  journal={Journal of computational physics},
  volume={83},
  number={1},
  pages={32--78},
  year={1989},
  publisher={Elsevier}
}

@article{shu1988efficient,
  title={Efficient implementation of essentially non-oscillatory shock-capturing schemes},
  author={Shu, Chi-Wang and Osher, Stanley},
  journal={Journal of computational physics},
  volume={77},
  number={2},
  pages={439--471},
  year={1988},
  publisher={Elsevier}
}

@article{dumbser2014posteriori,
  title={A posteriori subcell limiting of the discontinuous Galerkin finite element method for hyperbolic conservation laws},
  author={Dumbser, Michael and Zanotti, Olindo and Loub{\`e}re, Rapha{\"e}l and Diot, Steven},
  journal={Journal of Computational Physics},
  volume={278},
  pages={47--75},
  year={2014},
  publisher={Elsevier}
}

@article{gerolymos2009very,
  title={Very-high-order WENO schemes},
  author={Gerolymos, Georges A and S{\'e}n{\'e}chal, David and Vallet, Isabelle},
  journal={Journal of Computational Physics},
  volume={228},
  number={23},
  pages={8481--8524},
  year={2009},
  publisher={Elsevier}
}

@article{velasco2023spectral,
  title={Spectral difference method with a posteriori limiting: application to the Euler equations in one and two space dimensions},
  author={Velasco Romero, David A and Han-Veiga, Maria and Teyssier, Romain},
  journal={Monthly Notices of the Royal Astronomical Society},
  volume={520},
  number={3},
  pages={3591--3608},
  year={2023},
  publisher={Oxford University Press}
}

@article{ruuth2002two,
  title={Two barriers on strong-stability-preserving time discretization methods},
  author={Ruuth, Steven J and Spiteri, Raymond J},
  journal={Journal of Scientific Computing},
  volume={17},
  pages={211--220},
  year={2002},
  publisher={Springer}
}

@article{shu1988total,
  title={Total-variation-diminishing time discretizations},
  author={Shu, Chi-Wang},
  journal={SIAM Journal on Scientific and Statistical Computing},
  volume={9},
  number={6},
  pages={1073--1084},
  year={1988},
  publisher={SIAM}
}

@article{gottlieb2001strong,
  title={Strong stability-preserving high-order time discretization methods},
  author={Gottlieb, Sigal and Shu, Chi-Wang and Tadmor, Eitan},
  journal={SIAM review},
  volume={43},
  number={1},
  pages={89--112},
  year={2001},
  publisher={SIAM}
}

@article{evstigneev2016construction,
  title={On the construction and properties of WENO schemes order five, seven, nine, eleven and thirteen. Part 1. Construction and stability},
  author={Evstigneev, Nikolai Mikhailovich},
  journal={Computer research and modeling},
  volume={8},
  number={5},
  pages={721--753},
  year={2016}
}

@inproceedings{Evstigneev2016OnTC,
  title={On the construction and properties of WENO schemes order five, seven, nine, eleven and thirteen. Part 2. Numerical examples},
  author={Evstigneev, Nikolai Mikhailovich},
  year={2016},
  url={https://api.semanticscholar.org/CorpusID:202408991}
}

@article{sod1978survey,
  title={A survey of several finite difference methods for systems of nonlinear hyperbolic conservation laws},
  author={Sod, Gary A},
  journal={Journal of computational physics},
  volume={27},
  number={1},
  pages={1--31},
  year={1978},
  publisher={Elsevier}
}

@article{van1974towards,
  title={Towards the ultimate conservative difference scheme. II. Monotonicity and conservation combined in a second-order scheme},
  author={Van Leer, Bram},
  journal={Journal of computational physics},
  volume={14},
  number={4},
  pages={361--370},
  year={1974},
  publisher={Elsevier}
}

@article{van1982comparative,
  title={A comparative study of computational methods in cosmic gas dynamics},
  author={van Albada, Geert Dick and Van Leer, Bram and Roberts Jr, WW},
  journal={Astronomy and Astrophysics, vol. 108, no. 1, Apr. 1982, p. 76-84.},
  volume={108},
  pages={76--84},
  year={1982}
}

@article{sweby1984high,
  title={High resolution schemes using flux limiters for hyperbolic conservation laws},
  author={Sweby, Peter K},
  journal={SIAM journal on numerical analysis},
  volume={21},
  number={5},
  pages={995--1011},
  year={1984},
  publisher={SIAM}
}

@article{roe1986characteristic,
  title={Characteristic-based schemes for the Euler equations},
  author={Roe, Philip L},
  journal={Annual review of fluid mechanics},
  volume={18},
  number={1},
  pages={337--365},
  year={1986},
  publisher={Annual Reviews 4139 El Camino Way, PO Box 10139, Palo Alto, CA 94303-0139, USA}
}

@article{vanleer1977towardsiii,
title = {Towards the ultimate conservative difference scheme III. Upstream-centered finite-difference schemes for ideal compressible flow},
journal = {Journal of Computational Physics},
volume = {23},
number = {3},
pages = {263-275},
year = {1977},
issn = {0021-9991},
doi = {https://doi.org/10.1016/0021-9991(77)90094-8},
url = {https://www.sciencedirect.com/science/article/pii/0021999177900948},
author = {Bram {Van Leer}}
}

@inproceedings{chakravarthy1983high,
  title={High resolution applications of the Osher upwind scheme for the Euler equations},
  author={Chakravarthy, S and Osher, Stanley},
  booktitle={6th Computational Fluid Dynamics Conference Danvers},
  pages={1943},
  year={1983}
}

@Article{Godunov,
 Author = {Godunov, S. K.},
 Title = {A difference method for numerical calculation of discontinuous solutions of the equations of hydrodynamics},
 FJournal = {Matematicheski{\u{\i}} Sbornik. Novaya Seriya},
 Journal = {Mat. Sb., Nov. Ser.},
 Volume = {47},
 Pages = {271--306},
 Year = {1959},
 Language = {Russian},
 zbMATH = {3273813},
 Zbl = {0171.46204}
}

@article{abgrall1994essentially,
  title={On essentially non-oscillatory schemes on unstructured meshes: analysis and implementation},
  author={Abgrall, R{\'e}mi},
  journal={Journal of Computational Physics},
  volume={114},
  number={1},
  pages={45--58},
  year={1994},
  publisher={Elsevier}
}

@article{abgrall2023extensions,
  title={Extensions of Active Flux to arbitrary order of accuracy},
  author={Abgrall, Remi and Barsukow, Wasilij},
  journal={ESAIM: Mathematical Modelling and Numerical Analysis},
  volume={57},
  number={2},
  pages={991--1027},
  year={2023},
  publisher={EDP Sciences}
}

@phdthesis{lore_phd_thesis,
         address = {Z{\"u}rich},
           title = {Efficient Iterative Arbitrary High Order Methods: Adaptivity and Structure Preservation},
           month = {March},
          school = {University of Zurich},
          author = {Lorenzo Micalizzi},
            year = {2024},
             url = {https://doi.org/10.5167/uzh-258683}
}

@article{ciallella2024high,
  title={A high-order, fully well-balanced, unconditionally positivity-preserving finite volume framework for flood simulations},
  author={Ciallella, Mirco and Micalizzi, Lorenzo and Michel-Dansac, Victor and {\"O}ffner, Philipp and Torlo, Davide},
  journal={arXiv preprint arXiv:2402.12248},
  year={2024}
}

@article{veiga2024improving,
  title={On improving the efficiency of {ADER} methods},
  author={Veiga, Maria Han and Micalizzi, Lorenzo and Torlo, Davide},
  journal={Applied Mathematics and Computation},
  volume={466},
  pages={128426},
  year={2024},
  publisher={Elsevier}
}

@article{micalizzi2023efficient,
  title={Efficient Iterative Arbitrary High-Order Methods: an Adaptive Bridge Between Low and High Order},
  author={Micalizzi, Lorenzo and Torlo, Davide and Boscheri, Walter},
  journal={Communications on Applied Mathematics and Computation},
  pages={1--38},
  year={2023},
  publisher={Springer}
}

@article{lax1954weak,
  title={Weak solutions of nonlinear hyperbolic equations and their numerical computation},
  author={Lax, Peter D},
  journal={Communications on pure and applied mathematics},
  volume={7},
  number={1},
  pages={159--193},
  year={1954},
  publisher={Wiley Online Library}
}

@article{micalizzi2023new,
  title={A new efficient explicit deferred correction framework: analysis and applications to hyperbolic pdes and adaptivity},
  author={Micalizzi, Lorenzo and Torlo, Davide},
  journal={Communications on Applied Mathematics and Computation},
  pages={1--36},
  year={2023},
  publisher={Springer}
}

@inproceedings{fox1949some,
  title={Some new methods for the numerical integration of ordinary differential equations},
  author={Fox, Leslie and Goodwin, ET},
  booktitle={Mathematical Proceedings of the Cambridge Philosophical Society},
  volume={45},
  number={3},
  pages={373--388},
  year={1949},
  organization={Cambridge University Press}
}

@article{micalizzi2024novel,
  title={Novel well-balanced continuous interior penalty stabilizations},
  author={Micalizzi, Lorenzo and Ricchiuto, Mario and Abgrall, R{\'e}mi},
  journal={Journal of Scientific Computing},
  volume={100},
  number={1},
  pages={1--45},
  year={2024},
  publisher={Springer}
}

@article{Decoriginal,
    AUTHOR = {Dutt, Alok and Greengard, Leslie and Rokhlin, Vladimir},
     TITLE = {Spectral deferred correction methods for ordinary differential
              equations},
   JOURNAL = {BIT},
  FJOURNAL = {BIT. Numerical Mathematics},
    VOLUME = {40},
      YEAR = {2000},
    NUMBER = {2},
     PAGES = {241--266},
      ISSN = {0006-3835},
   MRCLASS = {65L05 (65F20 65L20)},
  MRNUMBER = {1765736},
MRREVIEWER = {Ekkehard Wagenf\"{u}hrer},
       DOI = {10.1023/A:1022338906936}, 
}

@article{Decremi,
    AUTHOR = {Abgrall, Rémi},
     TITLE = {High order schemes for hyperbolic problems using globally
              continuous approximation and avoiding mass matrices},
   JOURNAL = {Journal of Scientific Computing},
    VOLUME = {73},
      YEAR = {2017},
    NUMBER = {2-3},
     PAGES = {461--494},
      ISSN = {0885-7474},
   MRCLASS = {65M60 (65M22)},
  MRNUMBER = {3719595},
MRREVIEWER = {Beny Neta},
       DOI = {10.1007/s10915-017-0498-4},
}

@article{han2021dec,
	title={{Dec} and {Ader}: similarities, differences and a unified framework},
	author={Han Veiga, Maria and {\"O}ffner, Philipp and Torlo, Davide},
	journal={Journal of Scientific Computing},
	volume={87},
	number={1},
	pages={1--35},
	year={2021},
	publisher={Springer}
}

@article{ketcheson2014comparison,
	title={A comparison of high-order explicit {Runge--Kutta}, extrapolation, and deferred correction methods in serial and parallel},
	author={Ketcheson, David and Bin Waheed, Umair},
	journal={Communications in Applied Mathematics and Computational Science},
	volume={9},
	number={2},
	pages={175--200},
	year={2014},
	publisher={Mathematical Sciences Publishers}
}

@article{minion2003semi,
	title={Semi-implicit spectral deferred correction methods for ordinary differential equations},
	author={Minion, Michael L},
	journal={Communications in Mathematical Sciences},
	volume={1},
	number={3},
	pages={471--500},
	year={2003},
	publisher={International Press of Boston}
}

@article{huang2006accelerating,
	title={Accelerating the convergence of spectral deferred correction methods},
	author={Huang, Jingfang and Jia, Jun and Minion, Michael},
	journal={Journal of Computational Physics},
	volume={214},
	number={2},
	pages={633--656},
	year={2006},
	publisher={Elsevier}
}

@article{michel2022spectral,
	title={Spectral analysis of high order continuous {FEM} for hyperbolic {PDEs} on triangular meshes: influence of approximation, stabilization, and time-stepping},
	author={Michel, Sixtine and Torlo, Davide and Ricchiuto, Mario and Abgrall, R{\'e}mi},
	journal={arXiv preprint arXiv:2206.06150},
	year={2022}
}

@article{michel2021spectral,
	title={Spectral analysis of continuous {FEM} for hyperbolic {PDEs}: influence of approximation, stabilization, and time-stepping},
	author={Michel, Sixtine and Torlo, Davide and Ricchiuto, Mario and Abgrall, R{\'e}mi},
	journal={Journal of Scientific Computing},
	volume={89},
	number={2},
	pages={1--41},
	year={2021},
	publisher={Springer}
}

@article{jund2007arbitrary,
	title={{Arbitrary High-Order Finite Element Schemes and High-Order Mass Lumping}},
	author={Jund, S{\'e}bastien and Salmon, St{\'e}phanie},
	journal={International Journal of Applied Mathematics \& Computer Science},
	volume={17},
	number={3},
	year={2007},
	pages={375--393}
}

@article{abgrall2019high,
	title={High-order residual distribution scheme for the time-dependent {Euler} equations of fluid dynamics},
	author={Abgrall, R{\'e}mi and Bacigaluppi, Paola and Tokareva, Svetlana},
	journal={Computers \& Mathematics with Applications},
	volume={78},
	number={2},
	pages={274--297},
	year={2019},
	publisher={Elsevier}
}

@article{abgrall2020high,
	title={High order asymptotic preserving deferred correction implicit-explicit schemes for kinetic models},
	author={Abgrall, R{\'e}mi and Torlo, Davide},
	journal={SIAM Journal on Scientific Computing},
	volume={42},
	number={3},
	pages={B816--B845},
	year={2020},
	publisher={SIAM}
}

@Article{ADERNSE,
author ={M. Dumbser},
title = {Arbitrary High Order {PNPM} Schemes on Unstructured Meshes for the Compressible {{Navier}--{Stokes}} Equations}, 
journal = {Computers \& Fluids},
volume = {39},
pages = {60--76},
year = {2010},  
}

@book{ToroBook,
	title={Riemann Solvers and Numerical Methods for Fluid Dynamics: A Practical Introduction},
	author={Toro, E.F.},
	isbn={9783540498346},
	lccn={2009921818},
	url={https://books.google.com/books?id=SqEjX0um8o0C},
	year={2009},
	publisher={Springer},
	address={Berlin Heidelberg}
}

@incollection{toro2016riemann,
  title={The Riemann problem: solvers and numerical fluxes},
  author={Toro, Eleuterio F},
  booktitle={Handbook of Numerical Analysis},
  volume={17},
  pages={19--54},
  year={2016},
  publisher={Elsevier}
}

@book{leveque2002finite,
	title={Finite volume methods for hyperbolic problems},
	author={LeVeque, Randall J},
	volume={31},
	year={2002},
	publisher={Cambridge university press}
}

@article{minion2011hybrid,
	title={A hybrid parareal spectral deferred corrections method},
	author={Minion, Michael},
	journal={Communications in Applied Mathematics and Computational Science},
	volume={5},
	number={2},
	pages={265--301},
	year={2011},
	publisher={Mathematical Sciences Publishers}
}

@article{minion2004semi,
	title={Semi-implicit projection methods for incompressible flow based on spectral deferred corrections},
	author={Minion, Michael L},
	journal={Applied numerical mathematics},
	volume={48},
	number={3-4},
	pages={369--387},
	year={2004},
	publisher={Elsevier}
}

@article{speck2015multi,
	title={A multi-level spectral deferred correction method},
	author={Speck, Robert and Ruprecht, Daniel and Emmett, Matthew and Minion, Michael and Bolten, Matthias and Krause, Rolf},
	journal={BIT Numerical Mathematics},
	volume={55},
	number={3},
	pages={843--867},
	year={2015},
	publisher={Springer}
}

@article{layton2005implications,
	title={Implications of the choice of quadrature nodes for Picard integral deferred corrections methods for ordinary differential equations},
	author={Layton, Anita T and Minion, Michael L},
	journal={BIT Numerical Mathematics},
	volume={45},
	number={2},
	pages={341--373},
	year={2005},
	publisher={Springer}
}

@article{abgrall2021relaxation,
	author = {Abgrall, R\'emi and Le M\'el\'edo, \'Elise  and \"Offner, Philipp and Torlo, Davide},
	title = {Relaxation {Deferred} {Correction} {Methods} and their {Applications} to {Residual} {Distribution} {Schemes}},
	journal = {The SMAI Journal of computational mathematics},
	pages = {125--160},
	publisher = {Soci\'et\'e de Math\'ematiques Appliqu\'ees et Industrielles},
	volume = {8},
	year = {2022},
	doi = {10.5802/smai-jcm.82},
	language = {en}
}

@article{boscarino2018implicit,
	title={Implicit-explicit integral deferred correction methods for stiff problems},
	author={Boscarino, Sebastiano and Qiu, Jing-Mei and Russo, Giovanni},
	journal={SIAM Journal on Scientific Computing},
	volume={40},
	number={2},
	pages={A787--A816},
	year={2018},
	publisher={SIAM}
}

@article{boscarino2016error,
	title={Error estimates of the integral deferred correction method for stiff problems},
	author={Boscarino, Sebastiano and Qiu, Jing-Mei},
	journal={ESAIM: Mathematical Modelling and Numerical Analysis},
	volume={50},
	number={4},
	pages={1137--1166},
	year={2016},
	publisher={EDP Sciences}
}

@article{toro2005tvd,
  title={TVD fluxes for the high-order ADER schemes},
  author={Toro, Eleuterio F and Titarev, Vladimir A},
  journal={Journal of Scientific Computing},
  volume={24},
  pages={285--309},
  year={2005},
  publisher={Springer}
}

@incollection{toro2001towards,
	title={Towards very high order {G}odunov schemes},
	author={Toro, E. F. and Millington, RC and Nejad, LAM},
	booktitle={Godunov methods},
	pages={907--940},
	year={2001},
	publisher={Springer}
}

@article{titarev2002ader,
	title={ADER: Arbitrary high order {G}odunov approach},
	author={Titarev, Vladimir A and Toro, Eleuterio F},
	journal={Journal of Scientific Computing},
	volume={17},
	number={1-4},
	pages={609--618},
	year={2002},
	publisher={Springer}
}

@InProceedings{toro2024ader,
author="Toro, Eleuterio F.
and Titarev, Vladimir
and Dumbser, Michael
and Iske, Armin
and Goetz, Claus R.
and Castro, Crist{\'o}bal E.
and Montecinos, Gino I.
and Dematt{\`e}, Riccardo",
editor="Par{\'e}s, Carlos
and Castro, Manuel J.
and Morales de Luna, Tom{\'a}s
and Mu{\~{n}}oz-Ruiz, Mar{\'i}a Luz",
title="The ADER Approach for Approximating Hyperbolic Equations to Very High Accuracy",
booktitle="Hyperbolic Problems: Theory, Numerics, Applications. Volume I",
year="2024",
publisher="Springer Nature Switzerland",
address="Cham",
pages="83--105",
isbn="978-3-031-55260-1"
}

@article{dumbser2008unified,
	title={A unified framework for the construction of one-step finite volume and discontinuous {G}alerkin schemes on unstructured meshes},
	author={Dumbser, Michael and Balsara, Dinshaw S and Toro, Eleuterio F and Munz, Claus-Dieter},
	journal={Journal of Computational Physics},
	volume={227},
	number={18},
	pages={8209--8253},
	year={2008},
	publisher={Elsevier}
}

@inproceedings{dumbser2006arbitrary,
  title={Arbitrary high order finite volume schemes for linear wave propagation},
  author={Dumbser, Michael and Schwartzkopff, T and Munz, C-D},
  booktitle={Computational Science and High Performance Computing II: The 2nd Russian-German Advanced Research Workshop, Stuttgart, Germany, March 14 to 16, 2005},
  pages={129--144},
  year={2006},
  organization={Springer}
}

@article{dumbser2009very,
	title={{Very high order PNPM schemes on unstructured meshes for the resistive relativistic MHD equations}},
	author={Dumbser, Michael and Zanotti, Olindo},
	journal={Journal of Computational Physics},
	volume={228},
	number={18},
	pages={6991--7006},
	year={2009},
	publisher={Elsevier}
}

@article{boscheri2019high,
	title={{High order direct Arbitrary-Lagrangian-Eulerian (ALE) PNPM schemes with WENO Adaptive-Order reconstruction on unstructured meshes}},
	author={Boscheri, Walter and Balsara, Dinshaw S},
	journal={Journal of Computational Physics},
	volume={398},
	pages={108899},
	year={2019},
	publisher={Elsevier}
}

@article{ciallella2022arbitrary,
	title={An arbitrary high order and positivity preserving method for the shallow water equations},
	author={Ciallella, Mirco and Micalizzi, Lorenzo and {\"O}ffner, Philipp and Torlo, Davide},
	journal={Computers \& Fluids},
	volume={247},
	pages={105630},
	year={2022},
	publisher={Elsevier},
	doi="10.1016/j.compfluid.2022.105630"
}

@article{offner2020arbitrary,
	title={Arbitrary high-order, conservative and positivity preserving {P}atankar-type deferred correction schemes},
	author={{\"O}ffner, Philipp and Torlo, Davide},
	journal={Applied Numerical Mathematics},
	volume={153},
	pages={15--34},
	year={2020},
	publisher={Elsevier}
}

@article{ciallella2023arbitrary,
  title={Arbitrary high order weno finite volume scheme with flux globalization for moving equilibria preservation},
  author={Ciallella, Mirco and Torlo, Davide and Ricchiuto, Mario},
  journal={Journal of Scientific Computing},
  volume={96},
  number={2},
  pages={53},
  year={2023},
  publisher={Springer}
}

@article{abgrall2020multidimensional,
  title={Multidimensional staggered grid residual distribution scheme for Lagrangian hydrodynamics},
  author={Abgrall, R{\'e}mi and Lipnikov, Konstantin and Morgan, Nathaniel and Tokareva, Svetlana},
  journal={SIAM Journal on Scientific Computing},
  volume={42},
  number={1},
  pages={A343--A370},
  year={2020},
  publisher={SIAM}
}

@article{balsara2000monotonicity,
	title={Monotonicity preserving weighted essentially non-oscillatory schemes with increasingly high order of accuracy},
	author={Balsara, Dinshaw S and Shu, Chi-Wang},
	journal={Journal of Computational Physics},
	volume={160},
	number={2},
	pages={405--452},
	year={2000},
	publisher={Elsevier}
}

@article{shi2003resolution,
  title={Resolution of high order WENO schemes for complicated flow structures},
  author={Shi, Jing and Zhang, Yong-Tao and Shu, Chi-Wang},
  journal={Journal of Computational Physics},
  volume={186},
  number={2},
  pages={690--696},
  year={2003},
  publisher={Elsevier}
}

@article{hermes2012linear,
  title={Linear stability of WENO schemes coupled with explicit Runge--Kutta schemes},
  author={Hermes, V and Klioutchnikov, I and Olivier, H29246581253},
  journal={International journal for numerical methods in fluids},
  volume={69},
  number={6},
  pages={1065--1095},
  year={2012},
  publisher={Wiley Online Library}
}

@article{gao2020seventh,
  title={Seventh and ninth orders characteristic-wise alternative WENO finite difference schemes for hyperbolic conservation laws},
  author={Gao, Zhen and Fang, Li-Li and Wang, Bao-Shan and Wang, Yinghua and Don, Wai Sun},
  journal={Computers \& Fluids},
  volume={202},
  pages={104519},
  year={2020},
  publisher={Elsevier}
}

\end{document}